# Some series and integrals involving the Riemann zeta function, binomial coefficients and the harmonic numbers

## Volume I

Donal F. Connon

18 February 2008

## Abstract

In this series of seven papers, predominantly by means of elementary analysis, we establish a number of identities related to the Riemann zeta function, including the following:

$$\sum_{n=0}^{\infty}\frac{1}{2^n}\sum_{k=0}^{n}\int_a^b \binom{n}{k} p(x)\cos 2kx\, dx = \int_a^b p(x)\, dx$$

$$\sum_{n=1}^{\infty}\frac{1}{2^n}\sum_{k=1}^{n}\int_a^b \binom{n}{k} p(x)\sin 2kx\, dx = \int_a^b p(x)\cot x\, dx$$

$$\sum_{n=0}^{\infty}\int_a^b p(x)\cos^n x\, \cos nx\, dx = \int_a^b p(x) dx$$

$$\sum_{n=0}^{\infty}\int_a^b p(x)\cos^n x\, \sin nx\, dx = \int_a^b p(x)\cot x\, dx$$

$$\sum_{n=1}^{\infty} t^n \sum_{k=1}^{n}\binom{n}{k}\frac{x^k}{k^s} = \frac{1}{1-t} Li_s\left(\frac{xt}{1-t}\right)$$

$$\varsigma_a(3) = \frac{1}{2}\sum_{n=1}^{\infty}\frac{1}{n2^n}\left\{(H_n^{(1)})^2 + H_n^{(2)}\right\} = \sum_{n=0}^{\infty}\frac{1}{n2^n}\sum_{k=1}^{n}\binom{n}{k}\frac{(-1)^{k+1}}{k^2}$$

$$\varsigma_a(4) = \sum_{n=1}^{\infty}\frac{1}{n2^n}\left\{\frac{1}{6}\left(H_n^{(1)}\right)^3 + \frac{1}{2}H_n^{(1)}H_n^{(2)} + \frac{1}{3}H_n^{(3)}\right\} = \sum_{n=1}^{\infty}\frac{1}{n2^n}\sum_{k=1}^{n}\binom{n}{k}\frac{(-1)^{k+1}}{k^3}$$

$$\varsigma_a(5) = \sum_{n=1}^{\infty}\frac{1}{n2^n}\left\{\left(H_n^{(1)}\right)^4 + 6\left(H_n^{(1)}\right)^2 H_n^{(2)} + 8H_n^{(1)}H_n^{(3)} + 3\left(H_n^{(2)}\right)^2 + 6H_n^{(4)}\right\}$$

where $p(x)$ is a suitably behaved continuously differentiable function and $\varsigma_a(s)$ is the alternating Riemann zeta function.

Whilst this paper is mainly expository, some of the formulae reported in it are believed to be new, and the paper may also be of interest specifically due to the fact that most of the various identities have been derived by elementary methods.

**CONTENTS OF VOLUMES I TO VI:** **Volume/page**

**SECTION:**

**APPENDICES (Volume VI):**

**GUIDE TO LOCATION OF FORMULAE**

| | Contains formulae numbered: From | To |
|---|---|---|
| VOLUME I | (1.1) | (3.266) |
| VOLUME II(a) | (4.1.1) | (4.3.184) |
| VOLUME II(b) | (4.3.200) | (4.4.44l) |
| VOLUME III | (4.4.45) | (4.4.117t) |
| VOLUME IV | (4.4.118) | (4.4.252b) |
| VOLUME V | (5.1) | (8.62) |
| VOLUME VI | Appendices, Acknowledgements and References | |

The formulae in Appendix A in Volume VI are prefixed by the letter A, and so on for the other appendices.

Some of the identities considered in the remaining six volumes are summarised below:

**Volume II(a)**

$$\varsigma_a(5)=\sum_{n=1}^{\infty}\frac{1}{n2^n}\left\{\left(H_n^{(1)}\right)^4+6\left(H_n^{(1)}\right)^2H_n^{(2)}+8H_n^{(1)}H_n^{(3)}+3\left(H_n^{(2)}\right)^2+6H_n^{(4)}\right\}$$

$$24\varsigma(5)=2\sum_{k=1}^{\infty}\frac{H_k^{(3)}}{k^2}+3\sum_{k=1}^{\infty}\frac{H_k^{(1)}H_k^{(2)}}{k^2}+\sum_{k=1}^{\infty}\frac{\left[H_k^{(1)}\right]^3}{k^2}$$

$$\psi(u)=\sum_{k=0}^{\infty}\frac{1}{k+1}\sum_{j=0}^{k}\binom{k}{j}(-1)^j\log(u+j)$$

$$\log\Gamma(u)=\sum_{n=0}^{\infty}\frac{1}{n+1}\sum_{k=0}^{n}\binom{n}{k}(-1)^k(u+k)\log(u+k)+\frac{1}{2}-u+\frac{1}{2}\log(2\pi)$$

$$\sum_{k=1}^{\infty}\left[\frac{(2k)!}{2^{2k}(k!)^2}\right]^3\sum_{j=0}^{k-1}\frac{1}{2j+1}=\frac{\pi}{24}\frac{6\pi-\sqrt{2}\Gamma(1/4)}{\Gamma^4(3/4)}$$

$$2\log G(1+u)=-\sum_{n=0}^{\infty}\frac{1}{n+1}\sum_{k=0}^{n}\binom{n}{k}(-1)^k(u+k)^2\log(u+k)+2u\log\Gamma(u)+2u^2-u$$

$$+2\varsigma'(-1)-\frac{11}{12}-\frac{1}{2}(u-1)(3u+2)$$

$$\frac{1}{2}\int_0^x u\log\sin\pi u\,du=[\varsigma'(-2,x)+\varsigma'(-2,1-x)]-2x[\varsigma'(-1,x)-\varsigma'(-1,1-x)]+\frac{\varsigma(3)}{2\pi^2}$$

$$\int_0^1\log\Gamma_3(x+1)dx=-\frac{1}{24}\log(2\pi)+\frac{3\varsigma(3)}{8\pi^2}$$

where $\varsigma'(s,x)=\frac{\partial}{\partial s}\varsigma(s,x)$, $G(x)$ and $\Gamma_3(x)$ are the Barnes double and triple gamma functions.

**Volume II(b)**

$$\gamma_p(u)=-\frac{1}{p+1}\sum_{n=0}^{\infty}\frac{1}{n+1}\sum_{k=0}^{n}\binom{n}{k}(-1)^k\log^{p+1}(u+k)$$

$$\int_1^x \gamma_p(u)\,du = \frac{(-1)^{p+1}}{p+1}[\varsigma^{(p+1)}(0,x) - \varsigma^{(p+1)}(0)]$$

$$\gamma_1\left(\frac{1}{4}\right) = \frac{1}{2}[2\gamma_1 - 15\log^2 2 - 6\gamma\log 2] - \frac{1}{2}\pi\left[\gamma + 4\log 2 + 3\log\pi - 4\log\Gamma\left(\frac{1}{4}\right)\right]$$

$$\int_1^\infty \left[\frac{1}{1-t} - \frac{1}{t\log(1/t)}\right]\frac{\log\log(1/t)}{t^u}dt = \gamma\gamma_0(u) + 2\gamma_1(u) + \gamma\log u + \log^2 u$$

$$\sum_{n=0}^\infty \frac{1}{n+1}\sum_{k=0}^n \binom{n}{k}(-1)^k k\log^2(k+1) = \gamma_1 - \log(2\pi) - \frac{1}{2}\gamma^2 + \frac{1}{24}\pi^2 + \frac{1}{2}\log^2(2\pi)$$

$$\sum_{n=1}^\infty t^n \sum_{k=1}^n \binom{n}{k}\frac{x^k}{(k+y)^s} = \frac{xt}{(1-t)\Gamma(s)}\int_0^\infty \frac{u^{s-1}e^{-(y+1)u}}{\left[1-\left(1+xe^{-u}\right)t\right]}du = \frac{xt}{(1-t)^2}\Phi\left(\frac{xt}{(1-t)}, s, y+1\right)$$

where $\gamma_p(u)$ is the Stieltjes constant and $\Phi(z,s,u)$ is the Hurwitz-Lerch zeta function.

## Volume III

$$Li_{s+1}(x) = \sum_{n=1}^\infty \frac{1}{n2^n}\sum_{k=1}^n \binom{n}{k}\frac{x^k}{k^s}$$

$$\sum_{n=1}^\infty \frac{1}{n^2}\sum_{k=1}^n \binom{n}{k}(-1)^k \frac{x^k}{k^s} = \log x\, Li_{s+1}(x) - (s+1)Li_{s+2}(x)$$

$$\sum_{n=1}^\infty \frac{w^n}{n^p}\sum_{k=1}^n \binom{n}{k}\frac{(-1)^k x^k}{k^s} = \frac{(-1)^{s+1}x}{\Gamma(s+1)}\int_0^1 \frac{\log^s t\, Li_{p-1}\left[w(1-xt)\right]}{1-xt}dt$$

$$\sum_{n=0}^\infty \frac{1}{2^{n+1}}\sum_{k=0}^n \binom{n}{k}\frac{(-1)^k}{(k+x)^s} = \sum_{k=0}^\infty \frac{(-1)^k}{(k+x)^s}$$

$$\frac{x}{s-1}\sum_{n=0}^\infty \frac{1}{n+1}\sum_{k=0}^n \binom{n}{k}(-1)^k \frac{x^k}{(k+1)^{s-1}} = \frac{\log x}{s-1}Li_{s-1}(x) + Li_s(x)$$

$$\int_0^1 \frac{t^{\alpha-1} - t^{\beta-1}}{(1+t)\log t}dt = \sum_{n=0}^\infty \frac{1}{2^{n+1}}\sum_{k=0}^n \binom{n}{k}(-1)^k \log\frac{k+\alpha}{k+\beta} = \log\frac{\Gamma\left(\frac{1+\alpha}{2}\right)\Gamma\left(\frac{\beta}{2}\right)}{\Gamma\left(\frac{1+\beta}{2}\right)\Gamma\left(\frac{\alpha}{2}\right)}$$

$$\sum_{n=0}^{\infty}\frac{w^n}{n+1}\sum_{k=0}^{n}\binom{n}{k}\frac{x^k}{(k+1)^{s-1}}=\frac{(-1)^{s-1}}{w\Gamma(s-1)}\int_0^1\frac{\log^{s-2}v\log\left[1-w(1+xv)\right]}{1+xv}dv$$

$$\int_o^1\left\{t-\frac{Li_p\left[(1-y)t\right]}{1-y}\right\}\frac{dy}{\log y}=\sum_{n=1}^{\infty}\frac{t^{n+1}}{(n+1)^p}\sum_{k=0}^{n}\binom{n}{k}(-1)^{k+1}\log(1+k)$$

where $Li_s(x)$ is the polylogarithm function.

## Volume IV

$$\sum_{n=1}^{\infty}\frac{H_n^{(1)}}{n^s}x^n=\int_0^x\frac{Li_s(x)-Li_s(y)}{x-y}dy$$

$$\sum_{n=1}^{\infty}\frac{H_n^{(3)}}{n2^n}=\varsigma(2)\log^2 2-\frac{7}{8}\varsigma(3)\log 2+Li_4(1/2)-\frac{1}{6}\log^4 2$$

$$\int_0^{\infty}u^{x-1}e^{-ku}\log^n u\,du=\frac{1}{k^x}\sum_{j=0}^{n}(-1)^{n-j}\binom{n}{j}\Gamma^{(j)}(x)\log^{n-j}k$$

$$\sum_{n=1}^{\infty}\frac{1}{2^{n+1}}\sum_{k=1}^{n}\binom{n}{k}\frac{H_k^{(2)}}{k^4}=\varsigma(2)\varsigma(4)+\varsigma^2(3)-\frac{25}{12}\varsigma(6)$$

$$\sum_{n=1}^{\infty}\frac{1}{2^{n+1}}\sum_{k=1}^{n}\binom{n}{k}\frac{H_k^{(3)}}{k^3}=\frac{1}{2}\varsigma^2(3)+\frac{1}{2}\varsigma(6)$$

$$\frac{1}{3}\sum_{n=0}^{\infty}\frac{1}{n+1}\sum_{k=0}^{n}\binom{n}{k}\frac{(-1)^k H_{k+1}^{(2)}}{(k+1)^3}=\frac{4}{3}\varsigma^2(3)-\frac{29}{12}\varsigma(6)+\varsigma(2)\varsigma(4)$$

$$\frac{1}{2}\int_0^1\frac{\log^2(1-x)\log\left[1-u(1-x)\right]}{1-x}dx=-Li_4(u)$$

$$\sum_{n=1}^{\infty}\frac{\log n\cos(n\pi/2)}{(\pi n)^2}=\frac{1}{48}\left[\log(2\pi)+\gamma-1\right]+\frac{1}{4}\varsigma'(-1)$$

## Volume V

$$\varsigma(2n+1)=(-1)^{n+1}\frac{(2\pi)^{2n+1}}{(2n+1)!}\int_0^{1/2}B_{2n+1}(x)\cot(\pi x)dx$$

$$\frac{1}{2}\int_a^b p(x)\,dx = \sum_{n=0}^{\infty}\int_a^b p(x)\cos\alpha nx\,dx$$

$$\frac{1}{2}\int_a^b p(x)\cot(\alpha x/2)\,dx = \sum_{n=0}^{\infty}\int_a^b p(x)\sin\alpha nx\,dx$$

$$\int_0^{\pi/8} x\cot x\,dx = \frac{\pi}{16}\log\left[2-\sqrt{2}\right]+\frac{1}{8}\left[1-\sqrt{2}\right]G+\frac{1}{64}\left[\sqrt{2}\varsigma\left(2,\frac{1}{8}\right)-2\left(\sqrt{2}+1\right)\pi^2\right]$$

$$\sum_{n=1}^{\infty}\frac{Ci(n\pi)}{n^4} = \varsigma(4)\left[\gamma+\log\pi-\frac{55}{32}\right]-\varsigma'(4)$$

$$\frac{\pi^2}{240} = \sum_{n=0}^{\infty}\frac{(n-1)\varsigma(2n)}{(n+1)(n+2)(n+3)2^n}$$

$$\sum_{n=1}^{\infty}\frac{H_n^{(1)}H_n^{(2)}}{2^n} = \frac{11}{6}\pi^2\log 2-\frac{179}{20}\varsigma(3)-\frac{1}{3}\log^3 2$$

$$\frac{\pi}{8}\log 2+\frac{G}{2} = \frac{1}{\sqrt{2}}\sum_{n=0}^{\infty}\frac{1}{2^{3n}(2n+1)^2}\binom{2n}{n}$$

where $p(x)$ is a suitably behaved continuously differentiable function and $Ci(x)$ is the cosine integral.

## Volume VI

$$B_n = \sum_{k=0}^{n}\frac{(-1)^k k!}{k+1}S(n,k)$$

$$\int_0^1 \frac{x^{n-1}}{1+x^n}\log\log\left(\frac{1}{x}\right)dx = -\frac{\log(2)\log(2n^2)}{2n}$$

$$\int_0^1 \frac{1}{1+x}\log^{q-1}\left(\frac{1}{x}\right)\left[\log\log\left(\frac{1}{x}\right)\right]^2 dx = \Gamma''(q)\varsigma_a(q)+2\Gamma'(q)\varsigma_a'(q)+\Gamma(q)\varsigma_a''(q)$$

$$\gamma+\frac{1}{2n}-\sum_{k=1}^{2N+1}\frac{B_{2k}}{2kn^{2k}} < H_n-\log n < \gamma+\frac{1}{2n}-\sum_{k=1}^{2N}\frac{B_{2k}}{2kn^{2k}}$$

$$\lim_{n\to\infty}\left[\frac{1}{2}\left(H_n^{(1)}\right)^2-\gamma\log n-\frac{1}{2}\log^2 n\right] = \frac{1}{2}\gamma^2$$

$$\lim_{n\to\infty}\left[\sum_{k=1}^{n}\frac{\left(H_k^{(1)}\right)^2}{k}+\sum_{k=1}^{n}\frac{H_k^{(2)}}{k}-H_n^{(1)}H_n^{(2)}\right]=\frac{4}{3}\varsigma(3)$$

$$\lim_{n\to\infty}\left[\frac{1}{6}\left(H_n^{(1)}\right)^3+\frac{1}{2}H_n^{(1)}H_n^{(2)}-\frac{1}{6}\log^3 n-\frac{\gamma}{2}\log^2 n-\frac{1}{2}\left(\varsigma(2)+\gamma^2\right)\log n\right]=\frac{1}{2}\varsigma(2)\gamma+\frac{1}{6}\gamma^3$$

$$1-\frac{1}{2}\log(2\pi)=\sum_{n=1}^{\infty}\log\frac{1}{e}\left(1+\frac{1}{n}\right)^{n+\frac{1}{2}}$$

$$\log G(u+1)-\frac{1}{2}u\log(2\pi)+\frac{1}{2}u(u+1)=\sum_{n=1}^{\infty}\left[\frac{1}{2}u^2\log\left(1+\frac{1}{n}\right)-u+n\log\left(1+\frac{u}{n}\right)\right]$$

$$\sum_{n=1}^{\infty}\left[\frac{1}{6}x^3\log\left(1+\frac{1}{n}\right)-\frac{1}{2}x^2-nx+n(n+x)\log\left(1+\frac{x}{n}\right)\right]$$

$$=\left(\frac{1}{4}-2\log A\right)x+\frac{1}{4}x^2+(x-1)\log G(1+x)-2\log\Gamma_3(1+x)$$

$$\sum_{k=1}^{\infty}\frac{\left(H_k^{(1)}\right)^2+H_k^{(2)}-\left(\log k+\gamma\right)^2+\varsigma(2)+2\gamma^2}{k}=2\varsigma(2)\gamma+\frac{4}{3}\gamma^3+\frac{2}{3}\varsigma(3)-2\gamma\gamma_1-\gamma_2$$

$$\frac{1}{2}[\gamma-\log\pi]=\int_0^1\left[\frac{1}{(1+t)\log t}+\frac{2t}{1-t^2}\right]dt$$

$$\sum_{n=1}^{\infty}\frac{u^n}{n^p}\log\left(1+\frac{1}{n}\right)=-\sum_{n=1}^{\infty}\frac{(-1)^n}{n}Li_{n+p}(u)$$

$$\int_0^1\log\Gamma(x)\log\cos\pi x\,dx=-\frac{1}{2}\log 2\log 2\pi+\frac{\pi^2}{48}$$

$$\int_0^1\log x\log\sin(\pi x)dx=\frac{1}{2\pi}\sum_{n=1}^{\infty}\frac{si(2n\pi)}{n^2}+\log 2+\frac{\pi^2}{24}$$

where $A$ is the Glaisher-Kinkelin constant and $\Gamma_3(x)$ is the Barnes triple gamma function.

## 1. INTRODUCTION

The Riemann zeta function $\varsigma(s)$ is defined for complex values $s$ by [126, p.96]

$$(1.1)\quad \varsigma(s)=\sum_{n=1}^{\infty}\frac{1}{n^s}=\frac{1}{1-2^{-s}}\sum_{n=1}^{\infty}\frac{1}{(2n-1)^s}=\frac{1}{1-2^{-s}}\sum_{n=0}^{\infty}\frac{1}{(2n+1)^s}\qquad ,(\text{Re}(s)>1)$$

$$=\frac{1}{1-2^{1-s}}\sum_{n=1}^{\infty}\frac{(-1)^{n+1}}{n^s}=\frac{1}{1-2^{1-s}}\varsigma_a(s)\qquad ,(\text{Re}(s)>0;\ s\neq 1)$$

where $\varsigma_a(s)$ is the alternating zeta function.

The zeta function $\varsigma(s)$ is named after George Friedrich Bernhard Riemann (1826-1866) who, rather surprisingly, wrote only one paper on number theory during his tragically short life [81, p.49]. This paper, "Üeber die Anzahl der Primzahlen unter einer gegebenen Grösse", which was published in 1859 [112], was only eight pages long (an English translation, On the Number of Primes less than a given Magnitude, is contained in the appendix to Edward's book [57] and an electronic version is also available on the Trinity College Dublin website [112]). It was described by E.T.Bell, in his book "Men of Mathematics" [18, p.487], as one of Riemann's most profound and suggestive works and it gave birth to the famous Riemann Hypothesis:

$$\varsigma(s)=0 \quad \Rightarrow \quad s=1/2+i\tau$$

(except for the so-called trivial zeros which are briefly referred to later in (3.11)).

A vast literature now exists on the Riemann zeta function and includes treatises by Edwards [57], Ivić [81], Srivastava and Choi [126] and Titchmarsh [129] to name but a few. The purpose of this series of papers is to consider some aspects of $\varsigma(s)$, particularly in the case where $s$ is a real number: very little reference is made to complex analysis (largely because unfortunately I have forgotten most of the complex analysis which took up temporary residence in my brain when I was at university more than 30 years ago!) . Wherever possible, I have tried to make these papers self-contained and, with this in mind, I have presented detailed proofs so that a reader new to the topic may follow the analysis more easily (without constantly needing to refer to a multitude of other sources). Four years ago, I too was a "reader new to the topic" and my learning curve is presented herein.

According to Ivić [81, p.1], it appears that the real-valued $\varsigma(s)$ was first considered by the great Swiss mathematician Leonard Euler (1707-1783) [135, p.265]. Euler's first exploration involved a long outstanding question, the so-called Basel problem, which was to find a closed form expression for the sum of the infinite series:

$$\varsigma(2)=\sum_{n=1}^{\infty}\frac{1}{n^2}$$

Prima facie, it is difficult to even obtain a numerical approximation for the sum of this infinite series because it converges very slowly; for example, its partial sum to 1,000 terms is correct to only two decimal places (as mentioned in Dunham's marvellous book "Euler, The Master of us All" [56, p.43]). Writing from Basel, in his Tractatus de seriebus infinitis (1689), Jakob Bernoulli (1654-1705) included the following plea for help in determining the sum of the series, "If anyone finds and communicates to us that which thus far has eluded our efforts, great will be our gratitude" [56, p.42].

The first breakthrough came when Euler [63] proved the following identity in his paper, "De summatione innumerabilium progressionum" which was presented to the St. Petersburg Academy on 5 March 1731

(1.2) $$\varsigma(2) = \log x \log(1-x) + \sum_{n=1}^{\infty} \frac{x^n}{n^2} + \sum_{n=1}^{\infty} \frac{(1-x)^n}{n^2}$$

and, upon putting $x = 1/2$, he concluded that

(1.3) $$\varsigma(2) = \log^2 2 + 2\sum_{n=1}^{\infty} \frac{1}{2^n n^2}$$

Using the familiar Maclaurin series for $\log(1-x)$, this is equivalent to

(1.3a) $$\varsigma(2) = \left( \sum_{n=1}^{\infty} \frac{1}{2^n n} \right)^2 + 2\sum_{n=1}^{\infty} \frac{1}{2^n n^2}$$

and, by only taking 14 terms of the second series, Euler was able to show that $\varsigma(2)$ was equal to approximately 1.644934... , an answer correct to six decimal places [56, p.45].

Using more modern notation, (1.3a) can be written in terms of the harmonic numbers $H_n$ and the dilogarithm $Li_2(x)$ (refer to (3.11) for a proof of (1.4))

(1.4) $$\varsigma(2) = \frac{1}{4}\left( \sum_{n=1}^{\infty} \frac{H_n}{2^n} \right)^2 + 2Li_2(1/2)$$

where

$$H_n = \sum_{k=1}^{n} \frac{1}{k}$$

and $Li_2(x)$, which is a Dirichlet power series, is a particular case of the polylogarithm function $Li_n(x)$ defined by

(1.5a) $$Li_n(x) = \sum_{k=1}^{\infty} \frac{x^k}{k^n} \qquad , (|x| \le 1)$$

(1.5b) $$Li_1(x) = \sum_{k=1}^{\infty} \frac{x^k}{k} = -\log(1-x) \quad , (-1 \le x < 1)$$

We also have the degenerate cases

$$Li_0(x) = \frac{x}{1-x} \qquad Li_{-1}(x) = \frac{x}{(1-x)^2}$$

Later in this paper we will also use the notation $H_n^{(r)} = \sum_{k=1}^{n} \frac{1}{k^r}$ so that $H_n^{(1)} = H_n$.
As mentioned by Lewin [101], the dilogarithm was first studied in the second half of the 18th century by Euler and also by the less well known English mathematician, John Landen (1719-1790). Landen [94a] also investigated the trilogarithm, $Li_3(x)$, and some of his results are employed in Section 3 of this paper. Kirillov [85] indicates that the dilogarithm was also studied by Leibniz in 1696 and, according to Maximon [101b], this study was initiated in one of the letters from Leibniz to Jakob Bernoulli.

The identity (1.2) can be easily verified by differentiation to show that the right hand side is a constant, and then evaluating the constant as $\varsigma(2)$ by letting $x \to 0$. However, this proof is not as illuminating as the original derivation and, in this regard, the reader is recommended to read the excellent article by Ayoub, "Euler and the Zeta Function", [15] (for which he won the Lester R Ford prize for expository writing in 1975). The formula (1.3) corrects the misprint in [15, p.1073]. Euler's proof of (1.2) is set out below.

From the Maclaurin expansion it is easily seen that

$$\frac{\log(1-x)}{x} = -\frac{1}{x}\sum_{n=1}^{\infty}\frac{x^n}{n} = -\sum_{n=1}^{\infty}\frac{x^{n-1}}{n}$$

and by integrating we obtain

$$-\varsigma(2) = \int_0^1 \frac{\log(1-x)\,dx}{x}$$

Using the substitution $u = 1-x$ this becomes

$$-\varsigma(2) = \int_0^1 \frac{\log u\,du}{1-u}$$

(1.5c) $$= \int_0^x \frac{\log u\,du}{1-u} + \int_x^1 \frac{\log u\,du}{1-u} = I_1 + I_2$$

By the binomial theorem, $I_1$ can be represented as

$$I_1 = \int_0^x \sum_{k=0}^{\infty} u^k \log u \, du$$

Using integration by parts we have

$$\int_o^x u^k \log u \, du = \frac{u^{k+1}}{k+1} \log u \Bigg|_0^x - \int_0^x \frac{u^k du}{k+1}$$

$$= \frac{x^{k+1}}{k+1} \log x - \frac{x^{k+1}}{(k+1)^2} \quad , k \geq 0$$

Therefore we have

$$I_1 = \sum_{k=0}^{\infty} \frac{x^{k+1}}{k+1} \log x - \sum_{k=0}^{\infty} \frac{x^{k+1}}{(k+1)^2}$$

(1.6a)
$$= -\log(1-x)\log x - \sum_{k=1}^{\infty} \frac{x^k}{k^2}$$

We now put $t = 1-u$ in $I_2$, and expand $\log(1-t)$ by the Maclaurin power series, to obtain

$$I_2 = \int_0^{1-x} \frac{\log(1-t)\, dt}{t}$$

$$= -\sum_{n=1}^{\infty} \int_0^{1-x} \frac{t^{n-1}}{n} dt$$

(1.6b)
$$= -\sum_{n=1}^{\infty} \frac{(1-x)^n}{n^2}$$

and combining (1.6a) and (1.6b) we obtain Euler's formula

(1.6c)
$$\varsigma(2) = \log x \log(1-x) + \sum_{n=1}^{\infty} \frac{x^n}{n^2} + \sum_{n=1}^{\infty} \frac{(1-x)^n}{n^2}$$

This may be written as

$$\varsigma(2) = \log x \log(1-x) + Li_2(x) + Li_2(1-x)$$

(a much more direct derivation of the above identity is contained in [126, p.107] and also at (3.110a) in this paper).

As is well known, Euler's first triumph came in 1734 when, by making the bold, and then unproven, assumption that the $\sin x$ function could be written as an infinite product of its factors, he declared that [133, p.66]

(1.6d) $$\sin x = x\prod_{n=1}^{\infty}\left(1-\frac{x^2}{n^2\pi^2}\right)$$

Then, taking the standard Maclaurin expansion

(1.6e) $$\sin x = \sum_{n=0}^{\infty}(-1)^n\frac{x^{2n+1}}{(2n+1)!}$$

and, comparing coefficients of $x^3$ in both expansions, Euler demonstrated that

(1.6f) $$\varsigma(2) = \frac{\pi^2}{6}$$

A truly wonderful result! However, some of Euler's contemporaries were concerned about the validity of his method and, in an attempt to silence his critics, Euler subsequently produced another quite different proof employing greater rigour with the use of the $\sin^{-1} x$ function (see Ayoub's paper [15, p.1079] for more details and also Kimble's one page proof [84]). Many different proofs of (1.6d) are now available: a particularly elegant proof was given by Kortram [92] in 1996 and an even more concise exposition was provided in an answer to a question posed by Caris [40] in 1914. A short proof, based on Euler's reflection formula for the gamma function, is contained in Appendix E of Volume VI.

The next major triumph occurred in 1750 (according to Kline [87, p.449]) when Euler showed that

(1.7) $$\varsigma(2n) = \sum_{k=1}^{\infty}\frac{1}{k^{2n}} = (-1)^{n+1}\frac{2^{2n-1}\pi^{2n}B_{2n}}{(2n)!}$$

In [49] Cvijović and Klinowski refer to the year 1740. It should be noted however that the Russian journal was subject to significant publication time lags: Euler's paper "De seriebus quibusdam considerations" was actually presented to the St. Petersburg Academy in October 1739 but was not published in Commentarii Academiae Scientiarum Imperialis Petropolitanae until 1750.

In (1.7) the $B_n$ are the (rational) Bernoulli numbers given by the generating function

(1.8) $$\frac{x}{e^x-1} = \sum_{n=0}^{\infty}B_n\frac{x^n}{n!} \qquad , (|x| < 2\pi)$$

The closely associated Bernoulli polynomials, $B_n(t)$, are defined by

(1.9) $$\frac{xe^{tx}}{e^x-1} = \sum_{n=0}^{\infty}B_n(t)\frac{x^n}{n!} \qquad , (|x| < 2\pi)$$

and it is easily seen from (1.9) that $B_n = B_n(0)$. The Bernoulli numbers and polynomials are used later in this paper and further important properties are outlined in Appendix A of Volume VI (where it is also explained why the radius of convergence of the above series is equal to $2\pi$ ).
The celebrated formula in (1.7) expresses $\varsigma(2n)$ as a rational multiple of $\pi^{2n}$ but neither Euler, nor anyone else since, has been able to determine any such formula for $\varsigma(2n+1)$. Euler did however conjecture that

(1.10) $$\varsigma(3) = \alpha \log^2 2 + \beta \frac{\pi^2}{6} \log 2 \ \ (?)$$

for rational numbers $\alpha$ and $\beta$ [56, p.60] (albeit "dimensional" analysis would suggest that we should have $\log^3 2$ instead of $\log^2 2$. For my part, I suspect that a formula such as (1.10) does not exist because, using integer relation detection techniques, such a simple relationship would undoubtedly have been experimentally discovered by now. My conjecture is that either $\alpha$ or $\beta$ contains a factor of $\sqrt{2}$, or another small surd or perhaps $\log(2\pi)$ makes an appearance). As regards dimensions, I believe that there is a misprint in Dunham's book [56, p.60] because Kimoto and Wakayama have stated in [84a] that Euler's conjecture was of the form

(1.10a) $$\varsigma(3) = \alpha \log^3 2 + \beta \varsigma(2) \log 2$$

To date, even $\varsigma(3)$ has defied the many efforts made to find a closed form expression for it. (Many years ago Knopp [90, p.230] provided an interesting discussion as to what constitutes a closed form expression and, more recently, Borwein et al.[29, p.284] have commented that one already exists for $\varsigma(3)$ on the basis that $\mathrm{Li}_3(\tau^{-2})$ should be regarded as a fundamental constant in its own right (where $\tau$ is the golden mean) ). If, for example, $\pi$, $e$, $\log 2$ and $\sqrt{2}$ can themselves only be represented by infinite series, what does a closed form expression really mean?).

There was however another spectacular breakthrough in 1978 (just over 200 years after Euler's efforts!) when Roger Apéry (1916-1994), at the ripe old age of 62, proved that $\varsigma(3)$ was irrational (see in particular the papers by Apéry [11], Beukers [23], van der Poorten [131a], [131b] and [132], Zeilberger [141] and [142] and Zudilin [144]): $\varsigma(3)$ is accordingly now referred to as Apéry's constant in honour of the French mathematician. An interesting biography of Apéry appears in [10] which, inter alia, recounts an alarming encounter which he had with the Gestapo in 1944 when he was working for the French Resistance.

More than 20 years later, Rivoal in 2000 [113] succeeded in proving that there are infinitely many integers $n$ such that $\varsigma(2n+1)$ is irrational and, a few months later, he refined his search [114] by showing that at least one of the nine numbers $\varsigma(5)$, $\varsigma(7)$,..., $\varsigma(21)$ was irrational. The search was further narrowed in 2001 by Zudilin [143] to the four numbers $\varsigma(5)$, $\varsigma(7)$, $\varsigma(9)$ and $\varsigma(11)$.

According to Ayoub [15, p.1084], Euler returned to the zeta function, for what appears to be the last time, in 1772 in a paper entitled "Exercitationes Analyticae". Notwithstanding that by this time Euler had been blind for six years, through what Ayoub describes as "a striking and elaborate scheme", he was able to prove that

$$\int_0^{\pi/2} x \log \sin x \, dx = \frac{7}{16}\varsigma(3) - \frac{\pi^2}{8}\log 2 \tag{1.11}$$

Equation (1.11) has been derived by many other mathematicians in many different ways since Euler's time and, more recently, these include Amigó [7], Choi, Srivastava and Adamchik [45], Crandall and Buhler [48], Chen [43c], Espinosa and Moll [59], Nash and O'Connor [104] and Srivastava [125a]. Several additional proofs, obtained by more elementary methods, are contained in this series of papers.

The reader interested in Euler's work should consult The Euler Archive [63], a website which is making Euler's original publications available in pdf format (together with English translations of some of his papers). And remember, it is still good advice to heed Laplace's famous exhortation "Read Euler, read Euler. He is the master of us all" [56, p. xiii] (especially in view of the fact that 2007 is the 300th anniversary of Euler's birth).

Using integration by parts, it is easily seen that

$$\int_0^{\pi/2} x^2 \cot x \, dx = x^2 \log \sin x \Big|_0^{\pi/2} - 2\int_0^{\pi/2} x \log \sin x \, dx = -2\int_0^{\pi/2} x \log \sin x \, dx \tag{1.12}$$

and hence the integral in (1.12) is connected with $\varsigma(3)$ as represented by Euler in (1.11). Indeed, in this paper we give an elementary proof that

$$\varsigma(2n+1) = (-1)^{n+1} \frac{(2\pi)^{2n+1}}{(2n+1)!} \int_0^{1/2} B_{2n+1}(x) \cot(\pi x) \, dx \tag{1.13}$$

where $B_{2n+1}(x)$ are the Bernoulli polynomials defined in (1.9). Equation (1.13) is not new: other less elementary proofs are shown in [48] and [49]. These other proofs however usually require knowledge of identities such as the decomposition of the cotangent function or knowledge of the theory of Fourier series.

In the next Section we consider integrals involving $\cot x$ (on the premise that it makes sense to follow on from where the Master left off).

## 2. AN INTEGRAL INVOLVING COT X

The following identity involving $\cot x$ is easily verified by multiplying the numerator and denominator by the complex conjugate $(1 - \cos x e^{-ix})$.

(2.1)
$$\frac{1}{1-\cos x e^{ix}} = 1 + i\cot x$$

I first encountered this identity in 1970 when I was studying A level mathematics in Northern Ireland; the source was Further Mathematics by R.I. Porter [108, p.191]. Now, in 2004, and three times older than I was then ($n = 3!$), I think that it is apt that I should have consulted one of my earliest textbooks for my first ever mathematics paper (a large part of this paper was written in 2004: the next year or so just crept by unexpectedly!).

The left hand side of (2.1) can be expanded as a finite geometric series

(2.2)
$$\frac{1}{1-\cos x e^{ix}} = \sum_{n=0}^{N} (\cos x e^{ix})^n + \frac{(\cos x e^{ix})^{N+1}}{1-\cos x e^{ix}}$$

where we are simply using the basic identity

(2.3)
$$\frac{1}{1-y} = \sum_{n=0}^{N} y^n + \frac{y^{N+1}}{1-y}$$

Combining (2.1) and (2.2) and equating real and imaginary parts we have

(2.4)
$$\sum_{n=0}^{N} \cos^n x \cos nx + Q_N(x) = 1$$

(2.5)
$$\sum_{n=0}^{N} \cos^n x \sin nx + R_N(x) = \cot x$$

where

(2.6)
$$Q_N(x) + iR_N(x) = \cos^{N+1} x(\cos(N+1)x + i\sin(N+1)x)(1 + i\cot x)$$

$$= \cos^{N+1} x(\cos(N+1)x - \sin(N+1)x\cot x) + i\cos^{N+1} x(\sin(N+1)x + \cos(N+1)x\cot x)$$

If we now multiply (2.4) by a Riemann integrable function $p(x)$ and integrate over the range $[a,b]$, we have the finite summation (having legitimately changed the order of integration and summation in the finite sum)

(2.7)
$$\sum_{n=0}^{N} \int_a^b p(x)\cos^n x \cos nx\, dx + \int_a^b p(x) Q_N(x) dx = \int_a^b p(x) dx$$

and similarly from (2.5) we obtain

(2.8) $$\sum_{n=0}^{N}\int_a^b p(x)\cos^n x\sin nx\,dx+\int_a^b p(x)R_N(x)dx=\int_a^b p(x)\cot x\,dx$$

(where we now have our first sight of an integral involving $\cot x$).

First, let us consider how to evaluate (2.7). An obvious method of attack would be to employ the well-known expansion of $\cos^n x$ in terms of multiple angles. Reference to Gradshteyn and Ryzhik [74, p.30] provides the formulae:

$$\cos^{2n} x=\frac{1}{2^{2n}}\left[\sum_{k=0}^{n-1}2\binom{2n}{k}\cos 2(n-k)x+\binom{2n}{n}\right]$$

$$\cos^{2n-1} x=\frac{1}{2^{2n-2}}\left[\sum_{k=0}^{n-1}2\binom{2n-1}{k}\cos(2n-2k-1)x\right]$$

(and elegant derivations of these formulae are given in a paper by B. Wiener and J. Wiener, "De Moivre's Formula to the Rescue" [136]).

However, these formulae appeared rather cumbersome at first sight and hence, as before, I resorted to my rather simpler textbook, namely Further Mathematics. This indicated that the analysis could be considerably simplified by using the De Moivre type formula

(2.9) $$2^n\cos^n x(\cos nx+i\sin nx)=\left(1+e^{2ix}\right)^n$$

whose proof depends on the half-angle formula

$$1+e^{2ix}=1+\cos 2x+i\sin 2x=2\cos^2 x+2i\sin^2 x$$

$$=2\cos x(\cos x+i\sin x)$$

The right-hand side of (2.9) can then be expanded by the binomial theorem to produce

(2.10) $$\left(1+e^{2ix}\right)^n=\sum_{k=0}^{n}\binom{n}{k}e^{2ikx}=\sum_{k=0}^{n}\binom{n}{k}(\cos 2kx+i\sin 2kx)$$

Equating the real and imaginary parts of (2.9) and (2.10) we obtain

(2.11) $$2^n\cos^n x\cos nx=\sum_{k=0}^{n}\binom{n}{k}\cos 2kx$$

(2.12) $$2^n\cos^n x\sin nx=\sum_{k=1}^{n}\binom{n}{k}\sin 2kx$$

About four years after I did this work, I discovered that the following formulae are contained in Ramanujan's Notebook [21, Part I, p.246] for $x \le \pi/2$

(2.12a) $$2^n \cos^n x \cos(a+n)x = \sum_{k=1}^{\infty} \binom{n}{k} \cos(a+2k)x$$

(2.12b) $$2^n \cos^n x \sin(a+n)x = \sum_{k=1}^{\infty} \binom{n}{k} \sin(a+2k)x$$

and it may be noted that (2.12b) may be obtained by differentiating (2.12a) with respect to $a$.

Substituting (2.11) and (2.12) in (2.7) and (2.8) respectively we have

(2.13) $$\sum_{n=0}^{N} \frac{1}{2^n} \sum_{k=0}^{n} \int_a^b \binom{n}{k} p(x) \cos 2kx \, dx + Q_N = \int_a^b p(x) dx$$

(2.14) $$\sum_{n=1}^{N} \frac{1}{2^n} \sum_{k=0}^{n} \int_a^b \binom{n}{k} p(x) \sin 2kx \, dx + R_N = \int_a^b p(x) \cot x \, dx$$

where

$$Q_N = \int_a^b p(x) Q_N(x) dx$$

$$R_N = \int_a^b p(x) R_N(x) dx$$

Let us now consider the form of $Q_N$ as $N \to \infty$. We have

(2.15) $$Q_N = \int_a^b p(x) \cos^{N+1} x (\cos(N+1)x - \sin(N+1)x \cot x) dx$$

which, at first glance, appears to be a fairly complex integral.

However, since $|cos^{N+1} x| \le 1$ we have

(2.16) $$|Q_N| \le \left| \int_a^b p(x) \left( \cos(N+1)x - \sin(N+1)x \cot x \right) dx \right|$$

$$\le \int_a^b \left| p(x) \left( \cos(N+1)x - \sin(N+1)x \cot x \right) \right| dx$$

(2.16a) $$\leq \int_a^b \left| \, p(x)\cos(N+1)x \, \right| dx + \int_a^b \left| \, \frac{p(x)\cos x \sin(N+1)x}{\sin x} \, \right| dx$$

Using integration by parts we have for the first integral in (2.16a)

$$\int_a^b p(x)\cos(N+1)x\,dx = p(x)\frac{\sin(N+1)x}{N+1}\Bigg|_a^b - \int_a^b p'(x)\frac{\sin(N+1)x\,dx}{N+1}$$

Therefore, assuming that $p'(x)$ is bounded on $[a,b]$, it is clear that

$$\lim_{N\to\infty}\int_a^b p(x)\cos(N+1)x\,dx = 0$$

We therefore need to investigate the behaviour of the second integral in (2.16a) as $N\to\infty$. To do this, we require a weak version of the Riemann-Lebesgue lemma: a proof of the strong version based on Lebesgue integration theory is contained in Apostol's book, "Mathematical Analysis" [13, p.313] but the following version, given by Berndt in 1975 in an article entitled "Elementary Evaluation of $\varsigma(2n)$" [19], is sufficient for our purposes (we shall refer to Berndt's paper again in Volume V). In 1997, J.B. Dence employed a similar method in [53] where he made use of the Euler polynomials (which are defined in Appendix A of Volume VI).

**The Riemann-Lebesgue Lemma:**

Let $f(x)$ be twice continuously differentiable on $[a,b]$ and suppose that either (i) $\sin x$ has no zero in $[a,b]$ or (ii) if $\sin a = 0$ then $f(a)=0$ also. Then the Riemann-Lebesgue lemma states that

(2.17) $$\lim_{M\to\infty}\int_a^b f(x)\frac{\sin Mx\,dx}{\sin x} = 0$$

**Proof:**

Put $g(x) = f(x)/\sin x$ and integrate by parts to obtain

(2.18) $$\int_a^b g(x)\sin Mx\,dx = -\frac{g(x)\cos Mx}{M}\Bigg|_a^b + \int_a^b \frac{g'(x)\cos Mx\,dx}{M}$$

We now need to consider the following possibilities:

(i) If $\sin x \neq 0$ at $x=a$, then $g(x)$ is finite there. By hypothesis $f(x)$ is finite at $x=b$ and so the integrated term in (2.18) tends to 0 as $M\to\infty$.

(ii) Alternatively, if $\sin a = 0$ then by L'Hôpital's rule (remembering that $f(a)=0$)

$$\lim_{x\to a} g(x) = \lim_{x\to a} \frac{f'(x)}{\cos x}$$

and this limit exists because $f(x)$ is stated to be continuously differentiable at $x = a$ and $\cos a \neq 0$ (because we have assumed here that $\sin a = 0$. Hence the integrated term in (2.18) again approaches 0 as $M \to \infty$.

(iii) If $\sin x \neq 0$ at $x = a$, then

$$g'(x) = \frac{\sin x\, f'(x) - f(x)\cos x}{\sin^2 x} \tag{2.19}$$

is finite at $x = a$ and the integrand on the right hand side of (2.18) is continuous and hence bounded on $[a,b]$. It follows that the integral on the right hand side of (2.18) also tends to 0 as $M \to \infty$.

(iv) Alternatively, if $\sin a = 0$ we have, using L'Hôpital's rule and (2.19),

$$\lim_{x\to a} g'(x) = \lim_{x\to a} \frac{\sin x\, f''(x) + \cos x\, f'(x) + f(x)\sin x - f'(x)\cos x}{2\sin x \cos x}$$

$$= \lim_{x\to a} \frac{f''(x) + f(x)}{2\cos x} = \frac{f''(a)}{2\cos a}$$

because $f(x)$ is stated to be twice continuously differentiable on $[a,b]$ and $\cos a \neq 0$. So the integrand in (2.18) is again bounded on $[a,b]$ and the integral approaches 0 as $M \to \infty$.

In a similar fashion, we could also prove that

$$\lim_{M\to\infty} \int_a^b f(x) \frac{\cos Mx\, dx}{\sin x} = 0 \tag{2.20}$$

However, a more direct proof of (2.20) may be obtained by noting that in the proof of the Riemann-Lebesgue lemma (2.17) we did not assume that $M$ was an integer (or even a rational number). We could therefore substitute $M = M' + \pi/2$ in (2.17), with the result that $\sin Mx = \sin(M' + \pi/2)x = \cos M'x$, and the proof of (2.20) follows immediately (the Riemann-Lebesgue lemma will be used again in Sections 5, 6 and 7 of Volume V).

In equations (2.13) and (2.14) let us assume that $p(x)$ is twice continuously differentiable on $[a,b]$ and that $p(a) = 0$. Therefore, $p(x)\cos x$ satisfies the conditions of the Riemann-Lebesgue lemma and hence the second integral in (2.16a) approaches 0 as $N \to \infty$. Therefore we have

$$\lim_{N\to\infty} R_N = 0 \text{ and similarly } \lim_{N\to\infty} Q_N = 0$$

and, in the limit as $N \to \infty$, the equations (2.13) and (2.14) become the fundamental identities:

(2.21) $$\sum_{n=0}^{\infty} \frac{1}{2^n} \sum_{k=0}^{n} \int_a^b \binom{n}{k} p(x)\cos 2kx\, dx = \int_a^b p(x)dx$$

(2.22) $$\sum_{n=1}^{\infty} \frac{1}{2^n} \sum_{k=1}^{n} \int_a^b \binom{n}{k} p(x)\sin 2kx\, dx = \int_a^b p(x)\cot x dx$$

The only restrictions placed on $p(x)$ are that it is twice continuously differentiable and that $p(a) = 0$ if $\sin a = 0$. The range of integration is such that $\sin x \neq 0$ for all $x$ in $(a, b]$ (these restrictions can be made less prescriptive using the strong version of the Riemann-Lebesgue lemma).

From inspection it is apparent that in (2.2) we could validly substitute $x \to \alpha x$ where $\alpha$ may be regarded as a constant (or an independent variable). In this case, (2.21) and (2.22) can be generalised to

(2.23) $$\sum_{n=0}^{\infty} \frac{1}{2^n} \sum_{k=0}^{n} \int_a^b \binom{n}{k} p(x)\cos 2\alpha kx\, dx = \int_a^b p(x)dx$$

(2.24) $$\sum_{n=1}^{\infty} \frac{1}{2^n} \sum_{k=1}^{n} \int_a^b \binom{n}{k} p(x)\sin 2\alpha kx\, dx = \int_a^b p(x)\cot \alpha x\, dx$$

provided (i) $\sin \alpha x \neq 0 \;\forall\; x \in [a, b]$ or, alternatively, (ii) if $\sin \alpha a = 0$ then $p(a) = 0$ also.

It should be noted that if $b$ is regarded as a variable, the above equations effectively represent Fourier series. Equations (2.21) and (2.22) are referred to in this paper as the basic identities.

We shall also show in Volume V that the following identities hold under similar conditions

$$\frac{1}{2}\int_a^b p(x)\, dx = \sum_{n=0}^{\infty} \int_a^b p(x)\cos \alpha nx\, dx$$

$$\frac{1}{2}\int_a^b p(x)\cot(\alpha x/2)\, dx = \sum_{n=0}^{\infty} \int_a^b p(x)\sin \alpha nx\, dx$$

Subject to the same conditions, reference to (2.7) and (2.8) shows that we also have the identities

(2.25) $$\sum_{n=0}^{\infty}\int_a^b p(x)\cos^n x\,\cos nx\,dx=\int_a^b p(x)dx$$

and similarly from (2.5) we obtain

(2.26) $$\sum_{n=0}^{\infty}\int_a^b p(x)\cos^n x\,\sin nx\,dx=\int_a^b p(x)\cot x\,dx$$

Identities (2.25) and (2.26) are considered further in (8.56a) in Volume V. The Wolfram Integrator evaluates $\int_a^b x^2\cos^n x\,\sin nx\,dx$ in terms of hypergeometric functions, but I have yet to explore this aspect further.

## 3. RESULTS OBTAINED FROM THE BASIC IDENTITIES

Let us consider a simple case where $p(x)=x$ and the range of integration is $[0,\pi/2]$: it is clear that the necessary conditions for the Riemann-Lebesgue lemma are satisfied. We therefore have

**Theorem 3.1:**

(3.1) $$\sum_{n=1}^{\infty}\frac{1}{2^n}\sum_{k=1}^{n}\binom{n}{k}\frac{(-1)^{k+1}}{k}=\sum_{n=1}^{\infty}\frac{H_n}{2^n}=2\log 2$$

**Proof:**

We have using integration by parts

$$\int_0^{\pi/2} x\cot x\,dx = x\log\sin x\Big|_0^{\pi/2}-\int_0^{\pi/2}\log\sin x\,dx=-\int_0^{\pi/2}\log\sin x\,dx$$

It is believed that Euler was the first person to show that

(3.2) $$\int_0^{\pi/2}\log\sin x\,dx=-\frac{\pi}{2}\log 2$$

A proof is outlined as a question in Further Mathematics [108, p.253] and another elementary proof is given as a by-product in Russell’s one page paper entitled, “Another Eulerian-Type Proof” [116]. For completeness, a proof is summarised below.

$$2I = 2\int_0^{\pi/2} \log\sin x\, dx = \int_0^{\pi} \log\sin x\, dx = \int_0^{\pi} \log 2\sin\left(x/2\right)\cos\left(x/2\right) dx$$

$$= \pi\log 2 + \int_0^{\pi} \log\sin(x/2)\, dx + \int_0^{\pi} \log\cos\left(x/2\right) dx$$

and we then use the substitutions $x/2 \to t$ in the first integral, and $x/2 \to \pi/2 - t$ in the second integral. This result is quite fascinating insofar as we can evaluate a non-trivial integral without actually doing any integration per se.

Wiener [138a] gave the following interesting proof in 2001. Define $F(p),\ p \geq 0$ by

$$F(p) = \int_0^{\pi/2} \frac{\tan^{-1}(p\tan x)}{\tan x} dx$$

Differentiating with respect to $p$ we get

$$F'(p) = \int_0^{\pi/2} \frac{dx}{p^2\tan^2 x + 1}$$

With the substitution $u = \tan x$ we obtain

$$F'(p) = \frac{\pi}{2(p+1)}$$

and, since $F(0) = 0$, we obtain

$$F(p) = \frac{\pi}{2}\log(1+p)$$

Letting $p = 1$ we get

(3.2a) $$\int_0^{\pi/2} x\cot x\, dx = \frac{\pi}{2}\log 2$$

From (2.22) we have

(3.3) $$S_1 = \sum_{n=1}^{\infty} \frac{1}{2^n} \sum_{k=1}^{n} \int_0^{\pi/2} \binom{n}{k} x\sin 2kx\, dx = \int_0^{\pi/2} x\cot x\, dx = \frac{\pi}{2}\log 2$$

and integration by parts shows that

$$\int_0^{\pi/2} x\sin 2kx\, dx = -\frac{x\cos 2kx}{2k} + \frac{\sin 2kx}{4k^2}\Big|_0^{\frac{\pi}{2}} = \frac{\pi(-1)^{k+1}}{4k} \qquad ,(k \geq 1)$$

Therefore we have

(3.4) $$S_1 = \frac{\pi}{4}\sum_{n=1}^{\infty}\frac{1}{2^n}\sum_{k=1}^{n}\binom{n}{k}\frac{(-1)^{k+1}}{k} = \frac{\pi}{2}\log 2$$

In Section 4 it is demonstrated that

(3.5) $$\sum_{k=1}^{n}\binom{n}{k}\frac{(-1)^{k+1}}{k} = \sum_{k=1}^{n}\frac{1}{k} = H_n$$

where $H_n$ is the harmonic number. Rather unsurprisingly, we learn from M.E. Hoffman [80] that the identity (3.5) dates back to Euler [61] no less. Accordingly, we have

(3.6) $$S_1 = \frac{\pi}{4}\sum_{n=1}^{\infty}\frac{H_n}{2^n}$$

or

(3.7) $$\sum_{n=1}^{\infty}\frac{H_n}{2^n} = 2\log 2$$

The above result is well known and an alternative (and more direct) proof is outlined below.

As will be demonstrated in (3.27), we have the Cauchy product of two series

(3.8) $$-\frac{\log(1-x)}{1-x} = \sum_{n=1}^{\infty} H_n x^n$$

and letting $x = 1/2$ results in (3.7). Our use of the basic identity has therefore simply reproduced a known result in this case.

The following integral is contained in G&R [74, 3.832 1] and this provides yet another proof of (3.7).

(3.8a) $$\int_0^{\pi/2} x\cos^{p-1} x\sin ax\, dx = \frac{\pi}{2^{p+1}}\Gamma(p)\frac{\psi\left(\frac{p+a+1}{2}\right) - \psi\left(\frac{p-a+1}{2}\right)}{\Gamma\left(\frac{p+a+1}{2}\right)\Gamma\left(\frac{p-a+1}{2}\right)}$$

where $p > 0$ and $-(p+1) < a < p+1$ and $\Gamma(x)$ and $\psi(x)$ are the gamma and digamma functions defined in Volume II(a). Letting $p = n+1$ and $a = n$ this becomes

(3.8b) $$\int_0^{\pi/2} x\cos^n x\sin nx\,dx = \frac{\pi}{2^{n+2}}\Gamma(n+1)\frac{\psi(n+1)-\psi(1)}{\Gamma(n+1)\Gamma(1)}$$

$$= \frac{\pi}{2^{n+2}} H_n$$

because $\psi(n+1)-\psi(1) = H_n$ using (4.1.7a) in Volume II(a). The above formula is also contained in Ramanujan's Notebook [21, Part I, p.290].

The integral (3.8a) is considered further in (8.48) in Volume V.

Now using (2.8) with $p(x) = x$ we have

$$\sum_{n=0}^{\infty}\int_0^{\pi/2} x\cos^n x\sin nx\,dx = \frac{\pi}{4}\sum_{n=1}^{\infty}\frac{H_n}{2^n}$$

$$= \int_0^{\pi/2} x\cot x\,dx$$

(3.8c) $$= \frac{\pi}{2}\log 2$$

where we have used the identity in (3.2).

**Theorem 3.2:**

$$S_2 = \sum_{n=1}^{\infty}\frac{1}{2^n}\sum_{k=1}^{n}\int_0^{\pi/2}\binom{n}{k} x^2\sin 2kx\,dx = \int_0^{\pi/2} x^2\cot x\,dx = -\frac{7}{8}\varsigma(3) + \frac{\pi^2}{4}\log 2$$

**Proof:**

Using (2.22) and putting $p(x) = x^2$ we have

(3.9) $$S_2 = \sum_{n=1}^{\infty}\frac{1}{2^n}\sum_{k=1}^{n}\int_0^{\pi/2}\binom{n}{k} x^2\sin 2kx\,dx = \int_0^{\pi/2} x^2\cot x\,dx$$

$$= -2\int_0^{\pi/2} x\log\sin x\,dx$$

$$= -\frac{7}{8}\varsigma(3) + \frac{\pi^2}{4}\log 2$$

where we have used Euler's 1772 integral (1.11) in evaluating the final integral.

Integration by parts gives

(3.9a) $$\int_0^{\pi/2} x^2 \sin 2kx\, dx = \frac{\cos 2kx}{4k^3} - \frac{x^2\cos 2kx}{2k} + \frac{x\sin 2kx}{2k^2}\Bigg|_0^{\frac{\pi}{2}}$$

$$= \frac{(-1)^k}{4k^3} - \frac{1}{4k^3} - \frac{\pi^2(-1)^k}{8k} \qquad , k \ge 1$$

Therefore we have

(3.10) $$S_2 = \frac{1}{4}A - \frac{1}{4}B - \frac{\pi^2}{8}C$$

where

(3.10a) $$A = \sum_{n=1}^{\infty} \frac{1}{2^n} \sum_{k=1}^{n} \binom{n}{k} \frac{(-1)^k}{k^3}$$

(3.10b) $$B = \sum_{n=1}^{\infty} \frac{1}{2^n} \sum_{k=1}^{n} \binom{n}{k} \frac{1}{k^3}$$

(3.10c) $$C = \sum_{n=1}^{\infty} \frac{1}{2^n} \sum_{k=1}^{n} \binom{n}{k} \frac{(-1)^k}{k}$$

From (3.1) we know that $C = -2\log 2$, therefore we only need to determine the values of $A$ and $B$. At the end of this Section we will show how $A, B$ and $C$ may be more easily evaluated using $P_s(x)$ which is defined in (3.56).

The expression for $A$ has a structural resemblance to the formula found by Sondow [121] by applying the Euler series transformation method (which is covered in Knopp's excellent book [90, p.240]) to the alternating Riemann zeta function. Sondow's result was

(3.11) $$\varsigma_a(s) = \sum_{n=0}^{\infty} \frac{1}{2^{n+1}} \sum_{k=0}^{n} \binom{n}{k} \frac{(-1)^k}{(k+1)^s}$$

where the alternating Riemann zeta function is defined by

$$\varsigma_a(s) = \sum_{n=1}^{\infty} \frac{(-1)^{n+1}}{n^s} = \sum_{n=0}^{\infty} \frac{(-1)^n}{(n+1)^s}$$

and is sometimes called the Dirichlet eta function and often designated by $\eta(s)$. It is known that $\varsigma_a(s)$ is an analytic function for $\mathrm{Re}(s) > 0$.

We shall see in (4.4.79) in Volume III that

$$\varsigma_a(s,u) = \sum_{n=0}^{\infty} \frac{1}{2^{n+1}} \sum_{k=0}^{n} \binom{n}{k} \frac{(-1)^k}{(k+u)^s}$$

where $\varsigma_a(s,u)$ may be regarded as an alternating Hurwitz zeta function and, using (4.4.24a), this may be written as

$$\varsigma_a(s,u) = \sum_{n=0}^{\infty} \frac{(-1)^n}{(n+u)^s}$$

We have for $\text{Re}(s) > 1$

$$\varsigma_a(s) = \sum_{n=1}^{\infty} \frac{1}{n^s} - \sum_{n=1}^{\infty} \frac{2}{(2n)^s}$$

$$= \varsigma(s) - \frac{2}{2^s} \varsigma(s)$$

and hence

$$\varsigma_a(s) = (1-2^{1-s})\varsigma(s)$$

The above formula therefore enables us to define $\varsigma(s)$ for all $\text{Re}(s) > 0$ except for $s=1$.

The formula (3.11) is also reported in Havil's delightful book, "Gamma: Exploring Euler's Constant" [78, p.206]. Equation (3.11) is a globally convergent series for $\varsigma(s)$ and, except for $s=1$, provides an analytic continuation of $\varsigma(s)$ to the entire complex plane. As shown below it can be used, for example, to derive values for $\varsigma(0), \varsigma(-1)$ and $\varsigma(-2)$:

$$\varsigma_a(0) = \sum_{n=0}^{\infty} \frac{1}{2^{n+1}} \sum_{k=0}^{n} \binom{n}{k} (-1)^k$$

$$= \sum_{n=0}^{\infty} \frac{\delta_{0,n}}{2^{n+1}} = \frac{1}{2}$$

where $\delta_{i,j}$ is the Kronecker delta (and we have used the binomial theorem (4.1.2) from Volume II(a) with $x=1$). Therefore we have

(3.11a) $$\varsigma(0) = -\frac{1}{2}$$

It should be noted that we cannot automatically substitute $s=0$ in the formula $\varsigma_a(s) = (1-2^{1-s})\varsigma(s)$ because that equation is only valid for $\text{Re}(s) > 0$ (excluding $s=1$). Fortunately, Hardy [129, p.16] gave the following functional equation for the alternating zeta function

$$\varsigma_a(-s)=\left(1-\left[2^{-s}-1\right]^{-1}\right)\pi^{-s-1}s\Gamma(s)\sin(\pi s/2)\varsigma_a(1+s)$$

$$=2\frac{\left[2^{-s-1}-1\right]}{\left[2^{-s}-1\right]}\pi^{-s-1}s\Gamma(s)\sin(\pi s/2)\varsigma_a(1+s)$$

and it is this equation that enables us to equate $\varsigma_a(0)=-\varsigma(0)$. As can be seen from Ayoub's paper [15], this is precisely the functional equation for the zeta function which was first postulated by Euler many years before Riemann.

Similarly, using (3.11) with $s=-1$ we obtain

$$\varsigma_a(-1)=\sum_{n=0}^{\infty}\frac{1}{2^{n+1}}\sum_{k=0}^{n}\binom{n}{k}(-1)^k(k+1)$$

Using the binomial theorem we have

$$x(1-x)^n=\sum_{k=0}^{n}\binom{n}{k}(-1)^k x^{k+1}$$

and differentiation produces

(3.11aa) $$(1-x)^n-nx(1-x)^{n-1}=\sum_{k=0}^{n}\binom{n}{k}(-1)^k(k+1)x^k$$

Letting $x=1$ we obtain

$$\sum_{k=0}^{n}\binom{n}{k}(-1)^k(k+1)=\delta_{0,n}-n\delta_{1,n}$$

Therefore we have

$$\varsigma_a(-1)=\sum_{n=0}^{\infty}\frac{\delta_{0,n}-n\delta_{1,n}}{2^{n+1}}=\frac{1}{2^{0+1}}-\frac{1}{2^{1+1}}=\frac{1}{4}$$

and, using Hardy's functional equation we see that $\varsigma_a(-1)=-3\varsigma(-1)$, and hence we have

(3.11b) $$\varsigma(-1)=-\frac{1}{12}$$

We now multiply equation (3.11aa) by $x$ and differentiate to obtain

$$(1-x)^n - 3nx(1-x)^{n-1} + n(n-1)x^2(1-x)^{n-2} = \sum_{k=0}^{n} \binom{n}{k} (-1)^k (k+1)^2 x^k$$

and with $x = 1$ we have

$$\varsigma_a(-2) = \sum_{n=0}^{\infty} \frac{\delta_{0,n} - 3n\delta_{1,n} + n(n-1)\delta_{2,n}}{2^{n+1}} = \frac{1}{2^{0+1}} - \frac{3}{2^{1+1}} + \frac{2}{2^{2+1}} = 0$$

and hence we have found the first trivial zero of the Riemann zeta function. More generally we have [126, p.97]

(3.11c) $$\varsigma(-n) = -\frac{B_{n+1}}{n+1}$$

and since $B_{2n+1} = 0$ we obtain $\varsigma(-2n) = 0$. Some particular values of the zeta function are considered further in Appendix F of Volume VI.

The identity (3.11) has some history: it was conjectured by Knopp (1882-1957) around 1930, then proved by Hasse [77] in 1930 and subsequently rediscovered by Sondow in 1994. Hasse (1898-1979) also showed that

(3.12) $$\varsigma(s) = \frac{1}{s-1} \sum_{n=0}^{\infty} \frac{1}{n+1} \sum_{k=0}^{n} \binom{n}{k} \frac{(-1)^k}{(k+1)^{s-1}}$$

(3.12a) $$\varsigma(s,a) = \frac{1}{s-1} \sum_{n=0}^{\infty} \frac{1}{n+1} \sum_{k=0}^{n} \binom{n}{k} \frac{(-1)^k}{(k+a)^{s-1}}$$

where $\varsigma(s,a) = \sum_{n=0}^{\infty} \frac{1}{(n+a)^s}$ for $\mathrm{Re}(s) > 1$ is the Hurwitz zeta function. The above two formulae are valid for all $s$ except $s = 1$. It may be noted that $\varsigma(s,1) = \varsigma(s)$.

Proofs of (3.11) and (3.12) are shown (4.4.79) and (4.4.85) respectively of Volume III. A further proof of (3.11) has recently been given by Amore [8]-see (3.83) of this paper. A proof of (3.12a) is given in (4.4.24pi). Curiously, there are no references whatsoever to either Hasse or Sondow in the very extensive bibliography compiled by Srivastava and Choi [126], a surprising omission by these renowned mathematicians.

At first sight, the two Hasse identities look rather different. However, if we consider the function defined by

$$f(t,s) = \sum_{n=0}^{\infty} t^n \sum_{k=0}^{n} \binom{n}{k} \frac{(-1)^k}{(k+1)^s}$$

we see that they are in fact intimately related. Indeed we have

$$(s-1)\varsigma(s)=\sum_{n=0}^{\infty}\frac{1}{n+1}\sum_{k=0}^{n}\binom{n}{k}\frac{(-1)^k}{(k+1)^{s-1}}=\int_0^1 f(t,s-1)dx$$

$$\varsigma_a(s)=\sum_{n=0}^{\infty}\frac{1}{2^{n+1}}\sum_{k=0}^{n}\binom{n}{k}\frac{(-1)^k}{(k+1)^{s}}=\frac{1}{2}f(1/2,s)$$

In (4.4.99aiv) in Volume III we shall see that

$$\sum_{n=1}^{\infty}t^n\sum_{k=0}^{n}\binom{n}{k}\frac{x^k}{(k+y)^s}=\frac{1}{\Gamma(s)}\int_0^{\infty}\frac{u^{s-1}e^{-yu}\left(1+xe^{-u}\right)t}{\left[1-\left(1+xe^{-u}\right)t\right]}du$$

A quick internet search on the topic of Euler series transformations led me to a 1994 paper by Flajolet and Sedgewick [68] which, in turn, considered more esoteric subjects such as Mellin transforms, Rice/Nörlund integrals and Bell polynomials (see for example [86] and [105]] ). Their paper [68] also proved the following identities:

If we define $S_n(m)$ by

(3.13) $$S_n(m)=\sum_{k=1}^{n}\binom{n}{k}\frac{(-1)^k}{k^m}$$

for $n$ an integer, then $S_n(m)$ can be expressed in terms of the generalised harmonic numbers as

(3.14) $$-S_n(m)=\sum_{1m_1+2m_2+3m_3...=m}\frac{1}{m_1!\,m_2!\,m_3!...}\left(\frac{H_n^{(1)}}{1}\right)^{m_1}\left(\frac{H_n^{(2)}}{2}\right)^{m_2}\left(\frac{H_n^{(3)}}{3}\right)^{m_3}...$$

where $H_n^{(k)}$ are the generalised harmonic numbers defined by

(3.15) $$H_n^{(k)}=\sum_{j=1}^{n}\frac{1}{j^k}$$

The modus operandi of the summation in (3.14) is easily illustrated by the following example: if $m=4$, then $m_1+2m_2+3m_3+4m_4=4$ is satisfied by the integers in the following array

$$\begin{bmatrix} m_1 & m_2 & m_3 & m_4 \\ 4 & 0 & 0 & 0 \\ 0 & 2 & 0 & 0 \\ 1 & 0 & 1 & 0 \\ 0 & 0 & 0 & 1 \\ 2 & 2 & 0 & 0 \end{bmatrix}$$

and these powers are shown below (see also (4.3.29) in Volume II(a)).

$$\varsigma_a(5) = \sum_{n=1}^{\infty} \frac{1}{n2^n} \left\{ \left(H_n^{(1)}\right)^4 + 6\left(H_n^{(1)}\right)^2 H_n^{(2)} + 8H_n^{(1)} H_n^{(3)} + 3\left(H_n^{(2)}\right)^2 + 6H_n^{(4)} \right\}$$

It was at this very early stage of my mathematical adventure that I knew that I should begin to give up: I was only a dabbling amateur…and I really should leave this serious stuff to serious mathematicians. However, I did persevere some more!

The first few values of $S_n(m)$ given by Flajolet and Sedgewick are:

(3.16a) $$-S_n(1) = -\sum_{k=1}^{n} \binom{n}{k} \frac{(-1)^k}{k} = H_n^{(1)}$$

(3.16b) $$-S_n(2) = -\sum_{k=1}^{n} \binom{n}{k} \frac{(-1)^k}{k^2} = \frac{1}{2}\left(H_n^{(1)}\right)^2 + \frac{1}{2} H_n^{(2)}$$

(3.16c) $$-S_n(3) = -\sum_{k=1}^{n} \binom{n}{k} \frac{(-1)^k}{k^3} = \frac{1}{6}\left(H_n^{(1)}\right)^3 + \frac{1}{2} H_n^{(1)} H_n^{(2)} + \frac{1}{3} H_n^{(3)}$$

Well, at least with $S_n(1)$ we are back in familiar territory, since this corresponds with the well-known Euler formula in (3.5). However the other formulae, especially $S_n(3)$, seemed to be adding complexity rather than clarity to the problem. After much cogitation, I finally remembered that I had seen something similar in Adamchik's 1996 paper, "On Stirling Numbers and Euler Sums" [2]. As if by magic, Adamchik's paper contained the following identities:

(3.17) $$\sum_{k=1}^{n} \frac{H_k^{(1)}}{k} = \frac{1}{2}\left(H_n^{(1)}\right)^2 + \frac{1}{2} H_n^{(2)}$$

(3.18) $$\sum_{k=1}^{n} \frac{H_k^{(2)}}{k} + \sum_{k=1}^{n} \frac{H_k^{(1)}}{k^2} = H_n^{(3)} + H_n^{(1)} H_n^{(2)}$$

(3.19) $$\sum_{k=1}^{n} \frac{\left(H_k^{(1)}\right)^2}{k} + \sum_{k=1}^{n} \frac{H_k^{(2)}}{k} = \frac{1}{3}\left(H_n^{(1)}\right)^3 + H_n^{(1)} H_n^{(2)} + \frac{2}{3} H_n^{(3)}$$

$$= 2\sum_{k=1}^{n} \frac{1}{k} \sum_{j=1}^{k} \frac{H_j^{(1)}}{j}$$

(for completeness, elementary proofs of the Adamchik identities are contained in (4.4.169) et seq of Volume IV). Identity (3.19) is also employed in (E.59) of Appendix E in Volume VI.

The third identity (3.19) is equal to $-2S_n(3)$ and hence we have

(3.20a) $$S_n(3) = \sum_{k=1}^{n} \binom{n}{k} \frac{(-1)^k}{k^3}$$

(3.20b) $$= -\left\{ \frac{1}{6}\left(H_n^{(1)}\right)^3 + \frac{1}{2} H_n^{(1)} H_n^{(2)} + \frac{1}{3} H_n^{(3)} \right\}$$

(3.20c) $$= -\frac{1}{2} \left\{ \sum_{k=1}^{n} \frac{\left(H_k^{(1)}\right)^2}{k} + \sum_{k=1}^{n} \frac{H_k^{(2)}}{k} \right\}$$

We therefore have

(3.21a) $$A = \sum_{n=1}^{\infty} \frac{1}{2^n} \sum_{k=1}^{n} \binom{n}{k} \frac{(-1)^k}{k^3} = -\frac{1}{2} \sum_{n=1}^{\infty} \frac{1}{2^n} \left\{ \sum_{k=1}^{n} \frac{\left(H_k^{(1)}\right)^2}{k} + \sum_{k=1}^{n} \frac{H_k^{(2)}}{k} \right\}$$

(3.21b) $$= -\frac{1}{2} \sum_{n=1}^{\infty} \frac{1}{2^n} \sum_{k=1}^{n} \frac{\left(H_k^{(1)}\right)^2}{k} - \frac{1}{2} \sum_{n=1}^{\infty} \frac{1}{2^n} \sum_{k=1}^{n} \frac{H_k^{(2)}}{k}$$

(3.21c) $$= -\frac{1}{2} D - \frac{1}{2} E$$

As this paper neared completion, I found another way to derive (3.21a) without the need for the Flajolet and Sedgewick analysis: this was fortuitous since I had prefaced this paper with the comment that it primarily employed elementary methods of classical analysis. Details of this alternative proof are contained in Section 4.

Using (3.18) and (3.19) we can eliminate the term $H_n^{(1)} H_n^{(2)}$ to obtain

(3.21d) $$\sum_{k=1}^{n} \frac{\left(H_k^{(1)}\right)^2}{k} - \sum_{k=1}^{n} \frac{H_k^{(1)}}{k^2} = \frac{2}{3}\left[ H_n^{(3)} - \left(H_n^{(1)}\right)^3 \right]$$

We will in fact show in (4.4.155fi), (4.4.155h) and (4.4.155q) in Volume IV that

(3.21ei) $$-n \int_0^1 (1-t)^{n-1} \log t \, dt = H_n$$

(3.21eii) $$n \int_0^1 (1-t)^{n-1} \log^2 t \, dt = H_n^{(2)} + \left(H_n^{(1)}\right)^2$$

(3.21eiii) $$-n \int_0^1 (1-t)^{n-1} \log^3 t \, dt = 6\left( \frac{1}{6}\left[H_n^{(1)}\right]^3 + \frac{1}{2} H_n^{(1)} H_n^{(2)} + \frac{1}{3} H_n^{(3)} \right)$$

and we have also proved in (4.4.155zi) that

(3.21eiv) $$(-1)^{p+1} n\int_0^1 (1-t)^{n-1} \log^p t\, dt = p!\sum_{k=1}^{n} \binom{n}{k} \frac{(-1)^k}{k^p}$$

Indeed, since $\int_0^1 (1-t)^{n-1} \log^k t\, dt = \int_0^1 t^{n-1} \log^k (1-t) t\, dt$ one would automatically expect a connection with the Stirling numbers of the first kind defined in (3.105).

Let us first of all consider the infinite series for $D$ in (3.21c).

**Lemma 3.1:**

(3.22) $$D = \sum_{n=1}^{\infty} \frac{1}{2^n} \sum_{k=1}^{n} \frac{\left(H_k^{(1)}\right)^2}{k} = 2\sum_{n=1}^{\infty} \frac{\left(H_n^{(1)}\right)^2}{n2^n} = \frac{7}{4}\varsigma(3)$$

**Proof:**

The following identity holds whenever the series is absolutely convergent [90, p.138]

(3.23) $$\sum_{n=1}^{\infty} a_n \sum_{k=1}^{n} b_k = \sum_{n=1}^{\infty} b_n \sum_{k=n}^{\infty} a_k$$

The proof is very easy to visualise. Let $B_n = \sum_{k=1}^{n} b_k$ , then

$$\sum_{n=1}^{\infty} a_n B_n = a_1 B_1 + a_2 B_2 + a_3 B_3 + ...$$

$$= a_1 b_1 + a_2 (b_1 + b_2) + a_3 (b_1 + b_2 + b_3) + ...$$

$$= \quad b_1 \{a_1 + a_2 + a_3 + ... \}$$

$$+ b_2 \{0 \ + \ a_2 + a_3 + ...\}$$

$$+ b_3 \{0 \ \ + 0 \ + a_3 + ...\}$$

$$+ ...........................$$

$$= \sum_{n=1}^{\infty} b_n \left( \sum_{k=n}^{\infty} a_k \right)$$

Therefore we have

(3.24) $$D = \sum_{n=1}^{\infty} \frac{1}{2^n} \sum_{k=1}^{n} \frac{\left(H_k^{(1)}\right)^2}{k} = \sum_{n=1}^{\infty} \frac{\left(H_n^{(1)}\right)^2}{n} \sum_{k=n}^{\infty} \frac{1}{2^k}$$

and, from the elementary geometric series, we have

$$\sum_{k=n}^{\infty}\frac{1}{2^k}=\sum_{k=1}^{\infty}\frac{1}{2^k}-\sum_{k=1}^{n-1}\frac{1}{2^k}=1-\left(1-\frac{1}{2^{n-1}}\right)=\frac{1}{2^{n-1}}$$

Hence we get

(3.25)
$$D=2\sum_{n=1}^{\infty}\frac{\left(H_n^{(1)}\right)^2}{n2^n}$$

To proceed further, we require the following lemma.

**Lemma 3.2:**

(3.26)
$$\frac{1}{2}\log^2(1-x)+Li_2(x)=\sum_{n=1}^{\infty}\frac{H_n}{n}x^n$$

**Proof:**

From Knopp's excellent book, Theory and Application of Infinite Series [90, p.179], we have the familiar Cauchy product

(3.27)
$$\frac{1}{1-x}\sum_{n=0}^{\infty}a_n x^n=\sum_{n=0}^{\infty}s_n x^n$$

where $s_n=\sum_{k=1}^{n}a_k$, provided $|x|<1$ and $x$ is also less than the radius of convergence of the series $\sum_{n=0}^{\infty}a_n x^n$.

Therefore, using the Maclaurin series for $\log(1-x)$

$$-\log(1-x)=\sum_{n=1}^{\infty}\frac{x^n}{n}\qquad,\ (-1\le x<1)$$

we have

(3.28)
$$-\frac{\log(1-x)}{1-x}=\sum_{n=1}^{\infty}x^n\sum_{k=1}^{n}\frac{1}{k}=\sum_{n=1}^{\infty}H_n\,x^n$$

From (3.28) we have

(3.29)
$$-\frac{\log(1-t)}{1-t}=\sum_{n=1}^{\infty}H_n\,t^n$$

and integrating from 0 to $x$ we obtain

(3.30) $$\frac{1}{2}\log^2(1-x) = \sum_{n=1}^{\infty}\frac{H_n}{n+1}x^{n+1}$$

with both series convergent for $|x| < 1$. Now we can write the right-hand side of (3.30) as

$$\sum_{n=1}^{\infty}\frac{H_n}{n+1}x^{n+1} = \sum_{n=1}^{\infty}\frac{H_{n+1}}{n+1}x^{n+1} - \sum_{n=1}^{\infty}\frac{x^{n+1}/(n+1)}{n+1}$$

$$= \left(\frac{H_2x^2}{2} + \frac{H_3x^3}{3} + \ldots\right) - \left(\frac{x^2}{2^2} + \frac{x^3}{3^2} + \ldots\right)$$

(3.30a)
$$= \left(\frac{H_1x}{1} + \frac{H_2x^2}{2} + \frac{H_3x^3}{3} + \ldots\right) - \left(\frac{x}{1^2} + \frac{x^2}{2^2} + \frac{x^3}{3^2} + \ldots\right)$$

$$= \sum_{n=1}^{\infty}\frac{H_n}{n}x^n - Li_2(x)$$

Combining (3.30) and (3.30a) proves the lemma

(3.31) $$\frac{1}{2}\log^2(1-x) + Li_2(x) = \sum_{n=1}^{\infty}\frac{H_n}{n}x^n \qquad , |x|<1$$

A further proof of this is given in (3.105a).

Dividing (3.31) by $(1-x)$ and employing the identity in (3.27), we obtain for the right-hand side

(3.32) $$\sum_{n=1}^{\infty}\frac{H_n}{n}x^n\frac{1}{1-x} = \sum_{n=1}^{\infty}\frac{H_n}{n}x^n\sum_{n=0}^{\infty}x^n$$

$$= \sum_{n=0}^{\infty}x^n\sum_{k=1}^{n}\frac{H_k}{k}$$

$$= \frac{1}{2}\left[\sum_{n=1}^{\infty}\left(H_n^{(1)}\right)^2 x^n + \sum_{n=1}^{\infty}H_n^{(2)}x^n\right]$$

where we used Adamchik's formula (3.17). Combining (3.31) and (3.32), and employing the following identity,

(3.33) $$\frac{1}{1-x}\sum_{n=1}^{\infty}\frac{x^n}{n^r}=\sum_{n=1}^{\infty}H_n^{(r)}x^n=\frac{Li_r(x)}{1-x} \qquad ,x\in[0,1)$$

The above identity is easily proved by substituting $s_n = H_n^{(r)}$ in (3.27) and specifically for $r=2$ we have

(3.34) $$\frac{Li_2(x)}{1-x}=\sum_{n=1}^{\infty}H_n^{(2)}x^n \qquad ,x\in[0,1)$$

We then have

(3.35) $$\frac{\log^2(1-x)+Li_2(x)}{1-x}=\sum_{n=1}^{\infty}\left(H_n^{(1)}\right)^2x^n \qquad ,x\in[0,1)$$

In their paper, "Explicit Evaluation of Euler Sums" [28], Borwein et al. describe (3.35) as "an easily verified generating function": I suspected that there was a much more direct proof, but it eluded me for a long time (some two years after I wrote that, I discovered a shorter approach set out in (3.221b)).

Dividing (3.35) by $x$ and integrating over the interval [0, 1/2] we have (L'Hôpital's rule tells us that $\frac{\log^2(1-x)}{x}$ remains finite at $x=0$)

(3.36) $$\sum_{n=1}^{\infty}\frac{\left(H_n^{(1)}\right)^2}{n2^n}=\int_0^{1/2}\frac{\log^2(1-x)}{x(1-x)}dx+\int_0^{1/2}\frac{Li_2(x)}{x(1-x)}dx$$

(3.37) $$=\int_0^{1/2}\frac{\log^2(1-x)}{x}dx+\int_0^{1/2}\frac{\log^2(1-x)}{1-x}dx+\int_0^{1/2}\frac{Li_2(x)}{x}dx+\int_0^{1/2}\frac{Li_2(x)}{1-x}dx$$

The second integral in (3.37) is seen to be equal to $\frac{1}{3}\log^3 2$ (using the obvious substitution $u=\log(1-x)$) and the third integral is equal to $Li_3(1/2)$ (using the series definition (1.6) of the polylogarithm function).

Using integration by parts, it is easily seen that

(3.38) $$\int\frac{\log^2(1-x)}{x}dx=\log^2(1-x)\log x+2\int\frac{\log(1-x)\log x}{1-x}dx$$

Integration by parts again, and using the identities (obtained from the series definition of the relevant polylogarithm)

(3.39) $$\frac{d}{dx}Li_2(1-x)=-\frac{Li_1(1-x)}{1-x}=\frac{\log x}{1-x} \qquad \text{and} \qquad \frac{d}{dx}Li_3(1-x)=-\frac{Li_2(1-x)}{1-x}$$

we obtain

(3.40) $$\int \frac{\log(1-x)\log x}{1-x}dx = \log(1-x)Li_2(1-x) + \int \frac{Li_2(1-x)}{1-x}dx$$

(3.41) $$= \log(1-x)Li_2(1-x) - Li_3(1-x)$$

Therefore we have

(3.42) $$\int \frac{\log^2(1-x)}{x}dx = \log^2(1-x)\log x + 2\log(1-x)Li_2(1-x) - 2Li_3(1-x)$$

Therefore, using the Euler and Landen identities for the low dimension polylogarithms [126, pp.107, 114] (note that (3.43a) may be obtained from (1.2) by letting $x = 1/2$ : a proof of (3.43b) is contained at (3.115a)).

(3.43a) $$Li_2(1/2) = \frac{\pi^2}{12} - \frac{1}{2}\log^2 2$$

(3.43b) $$Li_3(1/2) = \frac{7}{8}\varsigma(3) - \frac{\pi^2}{12}\log 2 + \frac{1}{6}\log^3 2$$

(3.43c) $$Li_3(1) = \varsigma(3)$$

and we have

(3.44) $$\int_0^{1/2} \frac{\log^2(1-x)}{x}dx = \frac{1}{4}\varsigma(3) - \frac{1}{3}\log^3 2$$

The fourth integral in (3.37) $\int \frac{Li_2(x)}{1-x}dx$ can also be evaluated in a similar way using integration by parts

(3.45) $$\int \frac{Li_2(x)}{1-x}dx = -Li_2(x)\log(1-x) - \int \frac{\log^2(1-x)}{x}dx$$

because $Li_2'(x) = -\frac{\log(1-x)}{x}$ (using the series definition of $Li_2(x)$).
Therefore, using (3.44) we obtain

(3.46) $$\int_0^{1/2} \frac{Li_2(x)}{1-x}dx = Li_2(1/2)\log 2 - \left\{\frac{1}{4}\varsigma(3) - \frac{1}{3}\log^3 2\right\}$$

$$= \frac{\pi^2}{12}\log 2 - \frac{1}{4}\varsigma(3) - \frac{1}{6}\log^3 2$$

Collecting the terms together for (3.36) we obtain

$$\sum_{n=1}^{\infty}\frac{\left(H_n^{(1)}\right)^2}{n2^n}=\left(\frac{1}{4}\varsigma(3)-\frac{1}{3}\log^3 2\right)+\frac{1}{3}\log^3 2+Li_3(1/2)+\left\{\frac{\pi^2}{12}\log 2-\frac{1}{4}\varsigma(3)-\frac{1}{6}\log^3 2\right\}$$

More generally, as set out below, we can easily show that (see also (3.106))

(3.46a) $$\sum_{n=1}^{\infty}\frac{\left(H_n^{(1)}\right)^2}{n}x^n=-\frac{1}{3}\log^3(1-x)+Li_3(x)-Li_2(x)\log(1-x)$$

We note from (3.35) that

$$\frac{\log^2(1-x)+Li_2(x)}{1-x}=\sum_{n=1}^{\infty}\left(H_n^{(1)}\right)^2 x^n \qquad ,x\in[0,1)$$

Dividing this by $x$ and integrating we get

$$\int_0^x\frac{\log^2(1-x)}{x}dx+\int_0^x\frac{\log^2(1-x)}{1-x}dx+\int_0^x\frac{Li_2(x)}{x(1-x)}dx=\sum_{n=1}^{\infty}\frac{\left(H_n^{(1)}\right)^2}{n}x^n$$

We have

$$\int_0^x\frac{\log^2(1-x)}{1-x}dx=-\frac{1}{3}\log^3(1-x)$$

$$\int\frac{\log^2(1-x)}{x}dx=\log^2(1-x)\log x+2\int\frac{\log(1-x)\log x}{1-x}dx$$

$$\int\frac{\log(1-x)\log x}{1-x}dx=\log(1-x)Li_2(1-x)+\int\frac{Li_2(1-x)}{1-x}dx$$

$$=\log(1-x)Li_2(1-x)-Li_3(1-x)$$

Therefore we have

$$\int\frac{\log^2(1-x)}{x}dx=\log^2(1-x)\log x+2\log(1-x)Li_2(1-x)-2Li_3(1-x)$$

$$\int_0^x\frac{\log^2(1-x)}{x}dx=\log^2(1-x)\log x+2\log(1-x)Li_2(1-x)-2Li_3(1-x)+2\varsigma(3)$$

$$\int_0^x \frac{Li_2(x)}{x(1-x)}dx = \int_0^x \frac{Li_2(x)}{x}dx + \int_0^x \frac{Li_2(x)}{1-x}dx = Li_3(x) + \int_0^x \frac{Li_2(x)}{1-x}dx$$

The integral $\int \frac{Li_2(x)}{1-x}dx$ can also be evaluated using integration by parts

$$\int \frac{Li_2(x)}{1-x}dx = -Li_2(x)\log(1-x) - \int \frac{\log^2(1-x)}{x}dx$$

We then see that

$$\int_0^x \frac{\log^2(1-x)}{x}dx + \int_0^x \frac{Li_2(x)}{1-x}dx = -Li_2(x)\log(1-x)$$

We have seen above that

$$\int \frac{\log^2(1-x)}{x}dx = \log^2(1-x)\log x + 2\log(1-x)Li_2(1-x) - 2Li_3(1-x)$$

and hence we get

$$\int \frac{Li_2(x)}{1-x}dx = -Li_2(x)\log(1-x) - \log^2(1-x)\log x - 2\log(1-x)Li_2(1-x) + 2Li_3(1-x)$$

This gives us

$$\int_0^x \frac{Li_2(x)}{1-x}dx =$$

$$-Li_2(x)\log(1-x) - \log^2(1-x)\log x - 2\log(1-x)Li_2(1-x) + 2Li_3(1-x) - 2\varsigma(3)$$

and therefore we have

$$\int_0^x \frac{Li_2(x)}{x(1-x)}dx =$$

$$Li_3(x) - Li_2(x)\log(1-x) - \log^2(1-x)\log x - 2\log(1-x)Li_2(1-x) + 2Li_3(1-x) - 2\varsigma(3)$$

This then gives us

$$\sum_{n=1}^{\infty} \frac{\left(H_n^{(1)}\right)^2}{n} x^n = -\frac{1}{3}\log^3(1-x) + Li_3(x) - Li_2(x)\log(1-x)$$

Hence we have

(3.47) $$\sum_{n=1}^{\infty}\frac{\left(H_n^{(1)}\right)^2}{n2^n}=\frac{7}{8}\varsigma(3)$$

Therefore we obtain

(3.48) $$D=\sum_{n=1}^{\infty}\frac{1}{2^n}\sum_{k=1}^{n}\frac{\left(H_k^{(1)}\right)^2}{k}=\frac{7}{4}\varsigma(3)$$

We also note that

$$\frac{1}{1-x}\sum_{n=1}^{\infty}\frac{x^n}{n^2}=\sum_{n=1}^{\infty}H_n^{(2)}x^n=\frac{Li_2(x)}{1-x}$$

and the Wolfram Integrator gives the integral for $|x| < 1$

$$\int\frac{x^{n-1}}{1-x}dx=\frac{1}{n}x^n\,{}_2F_1(n,1;n+1;x)$$

Therefore we have

$$\int_0^u\frac{Li_2(x)}{x(1-x)}dx=\sum_{n=1}^{\infty}\frac{u^n}{n^3}\,{}_2F_1(n,1;n+1;u)$$

$$=Li_3(u)-Li_2(u)\log(1-u)-\log^2(1-u)\log u-2\log(1-u)Li_2(1-u)+2Li_3(1-u)-2\varsigma(3)$$

We now need to determine the sum of the second infinite series in (3.21b)

**Lemma 3.3:**

(3.49) $$E=\sum_{n=1}^{\infty}\frac{1}{2^n}\sum_{k=1}^{n}\frac{H_k^{(2)}}{k}=2\sum_{n=1}^{\infty}\frac{H_n^{(2)}}{n2^n}=\frac{5}{4}\varsigma(3)$$

**Proof:**

Since the above series is absolutely convergent, it can be rearranged in the same manner as (3.23) to give

(3.50) $$E=2\sum_{n=1}^{\infty}\frac{H_n^{(2)}}{n2^n}$$

We now use (3.34) above

(3.51) $$\sum_{n=1}^{\infty} H_n^{(2)} x^n = \frac{Li_2(x)}{1-x}$$

Again, using the same technique of dividing by $x$ and integrating (3.51) over the interval [0,1/2], we obtain

(3.52) $$\sum_{n=1}^{\infty} \frac{H_n^{(2)}}{n2^n} = \int_0^{1/2} \frac{Li_2(x)}{x(1-x)} dx = \int_0^{1/2} \frac{Li_2(x)}{x} dx + \int_0^{1/2} \frac{Li_2(x)}{(1-x)} dx$$

$$= Li_3(1/2) + \frac{\pi^2}{12}\log 2 - \frac{1}{4}\varsigma(3) - \frac{1}{6}\log^3 2$$

(3.53) $$= \frac{5}{8}\varsigma(3)$$

where we have used (3.46).Therefore

(3.54) $$E = \sum_{n=1}^{\infty} \frac{1}{2^n} \sum_{k=1}^{n} \frac{H_k^{(2)}}{k} = \frac{5}{4}\varsigma(3)$$

It now remains to evaluate $B$ contained in equation (3.10b). We can in fact significantly generalise the analysis at this stage by reference to the following lemma.

**Lemma 3.4:**

(3.55) $$Li_s(x) = \sum_{n=1}^{\infty} \frac{1}{2^{n+1}} \sum_{k=1}^{n} \binom{n}{k} \frac{x^k}{k^s}$$

**Proof:**

Let us define $P_s(x)$ by

(3.56) $$P_s(x) = \sum_{n=1}^{\infty} \frac{1}{2^n} \sum_{k=1}^{n} \binom{n}{k} \frac{x^k}{k^s}$$

Writing the terms out explicitly gives

$$P_s(x) = \frac{1}{2^1}\left\{\binom{1}{1}\frac{x}{1^s}\right\}$$

$$+\frac{1}{2^2}\left\{\binom{2}{1}\frac{x}{1^s} + \binom{2}{2}\frac{x^2}{2^s}\right\}$$

$$+\frac{1}{2^3}\left\{\binom{3}{1}\frac{x}{1^s}+\binom{3}{2}\frac{x^2}{2^s}+\binom{3}{3}\frac{x^3}{3^s}\right\}$$

$$+\frac{1}{2^4}\left\{\binom{4}{1}\frac{x}{1^s}+\binom{4}{2}\frac{x^2}{2^s}+\binom{4}{3}\frac{x^3}{3^s}+\binom{4}{2}\frac{x^4}{4^s}\right\}$$

$$+\ldots\ldots\ldots\ldots\ldots\ldots\ldots\ldots\ldots\ldots$$

Assuming that the series is absolutely convergent for $\mathrm{Re}(s)>1$, it may be rearranged by adding the terms vertically in columns as follows (the reason why $s$ must be greater than 1 in order to secure convergence for $P_s(1)$ is explained in Volume III).

(3.57) $$P_s(x)=\frac{x}{1^s}\left\{\frac{1}{2^1}\binom{1}{1}+\frac{1}{2^2}\binom{2}{1}+\frac{1}{2^3}\binom{3}{1}+\frac{1}{2^4}\binom{4}{1}+\ldots\right\}$$

$$+\frac{x^2}{2^s}\left\{\quad 0\quad+\frac{1}{2^2}\binom{2}{2}+\frac{1}{2^3}\binom{3}{2}+\frac{1}{2^4}\binom{4}{2}+\ldots\right\}$$

$$+\frac{x^3}{3^s}\left\{\quad 0\quad+\quad 0\quad+\frac{1}{2^3}\binom{3}{3}+\frac{1}{2^4}\binom{4}{3}+\ldots\right\}$$

$$+\ldots\ldots\ldots\ldots\ldots\ldots\ldots\ldots\ldots\ldots$$

The series in parentheses in the j[th] row of (3.57) may be expressed as

(3.58) $$\sum_{i=0}^{\infty}\frac{1}{2^i}\binom{i}{j}$$

where, in this case, the summation is over the **upper** index. The above expression correctly gives the zero terms in the j[th] row because

$$\binom{i}{j}=0\ \text{ if } j>i\ge 0$$

by virtue of the definition of the binomial coefficients [75, p.155]. Hence, (3.56) can be replaced by

(3.59) $$P_s(x)=\sum_{j=1}^{\infty}\frac{x^j}{j^s}\left\{\sum_{i=0}^{\infty}\frac{1}{2^i}\binom{i}{j}\right\}$$

We now employ a neat trick revealed in [75, p.199]. If we expand $(1-y)^{-j-1}$ by the binomial theorem, the coefficient of $y^i$ is $\binom{-j-1}{i}(-1)^i$ and this can be written as $\binom{j+i}{i}$ by negating the upper index. This then gives

(3.60) $$\frac{1}{(1-y)^{j+1}}=\sum_{i=0}^{\infty}\binom{j+i}{j}y^i$$

Now multiplying this identity by $y^j$ gives

(3.61) $$\frac{y^j}{(1-y)^{j+1}}=\sum_{i=0}^{\infty}\binom{i}{j}y^i$$

Letting $y=1/2$ in (3.61) produces

(3.62) $$\sum_{i=0}^{\infty}\frac{1}{2^i}\binom{i}{j}=2$$

and, curiously, the result does not depend on the value of $j$ (and therefore the sum of each row in parentheses in (3.57) is equal to 2). Hence we have the result

(3.64) $$P_s(x)=\sum_{n=1}^{\infty}\frac{1}{2^n}\sum_{k=1}^{n}\binom{n}{k}\frac{x^k}{k^s}=2\sum_{k=1}^{\infty}\frac{x^k}{k^n}=2Li_s(x)$$

and specifically we have $B=P_3(1)=2Li_3(1)=2\varsigma(3)$.

An alternative proof of (3.61) is set out below. We have

$$\frac{d^j}{dy^j}y^i=i(i-1)(i-2)...(i-j+1)y^{i-j}$$

and therefore we directly obtain

$$\sum_{i=0}^{\infty}\binom{i}{j}y^i=\frac{y^j}{j!}\sum_{i=0}^{\infty}\frac{d^j}{dy^j}y^i=\frac{y^j}{j!}\frac{d^j}{dy^j}\left(\frac{1}{1-y}\right)=\frac{y^j}{(1-y)^{j+1}}$$

We may also note the connection with the Stirling numbers $s(n,k)$ of the first kind (3.104) in the above derivation.

$$\frac{d^n}{dy^n}y^x=x(x-1)(x-2)...(x-n+1)y^{x-n}=y^{x-n}\sum_{k=0}^{n}s(n,k)x^k$$

If $x=e^{i\theta}$ then we have

(3.64a) $$\sum_{k=1}^{\infty}\frac{\cos k\theta}{k^n}=\sum_{n=1}^{\infty}\frac{1}{2^{n+1}}\sum_{k=1}^{n}\binom{n}{k}\frac{\cos k\theta}{k^s}$$

(3.64b) $$\sum_{k=1}^{\infty}\frac{\sin k\theta}{k^n}=\sum_{n=1}^{\infty}\frac{1}{2^{n+1}}\sum_{k=1}^{n}\binom{n}{k}\frac{\sin k\theta}{k^s}$$

It should be noted that s may be a complex number provided $\operatorname{Re}(s) > 1$. The structural similarity of $P_s(-1) = 2\,Li_s(-1)$ with the Hasse/Sondow formula (3.11) should be noted. This is discussed in more detail in Part IV of Section 4. The identity (3.64) is significantly generalised in (3.67a).

When I first derived this lemma, I did not include $x$ as a variable and hence I originally obtained the result

(3.65) $$P_s(1)=\sum_{n=1}^{\infty}\frac{1}{2^n}\sum_{k=1}^{n}\binom{n}{k}\frac{1}{k^s}=2\varsigma(s)$$

It was therefore not initially apparent to me that the series $P_s(x)$ was so intimately connected with the polylogarithm function. Furthermore, it was not until the final stages of the preparation of this paper, that I even realised that

(3.66) $$P_s(-1)=\sum_{n=1}^{\infty}\frac{1}{2^n}\sum_{k=1}^{n}\binom{n}{k}\frac{(-1)^k}{k^s}=2Li_s(-1)$$

and this explains why I independently proved the following Theorem 3.3 (and indeed the analysis in Section 5 in connection with the Bernoulli polynomials). A further proof of this lemma, in the case where $s$ is a positive integer, is shown below.

But, before that, let's just finish off the proof of Theorem 3.2 which has already occupied some 19 pages of this paper!

We have from (3.10)

$$S_2=\frac{1}{4}A-\frac{1}{4}B-\frac{\pi^2}{8}C$$

where

$$A=\sum_{n=1}^{\infty}\frac{1}{2^n}\sum_{k=1}^{n}\binom{n}{k}\frac{(-1)^k}{k^3}=-\frac{1}{2}D-\frac{1}{2}E$$

$$=-\frac{3}{2}\varsigma(3)=2Li_3(-1) \qquad \text{[using (3.22), (3.49) and (3.55)]}$$

$$B=\sum_{n=1}^{\infty}\frac{1}{2^n}\sum_{k=1}^{n}\binom{n}{k}\frac{1}{k^3}=2\varsigma(3)=2Li_3(1) \qquad \text{[using (3.65)]}$$

$$C=\sum_{n=1}^{\infty}\frac{1}{2^n}\sum_{k=1}^{n}\binom{n}{k}\frac{(-1)^k}{k}=-2\log 2=2Li_1(-1) \qquad \text{[using (3.1)]}$$

Therefore we have (at long last!)

(3.67) $$\sum_{n=1}^{\infty}\frac{1}{2^n}\sum_{k=1}^{n}\int_0^{\pi/2}\binom{n}{k}x^2\sin 2kx\,dx=-\frac{7}{8}\varsigma(3)+\frac{\pi^2}{4}\log 2$$

An interesting generalisation of Lemma 3.4 is set out below.

**Lemma 3.4(a):**

(3.67a) $$\sum_{n=1}^{\infty}t^n\sum_{k=1}^{n}\binom{n}{k}\frac{x^k}{k^s}=\frac{1}{1-t}Li_s\left(\frac{xt}{1-t}\right)$$

$$\sum_{n=1}^{\infty}t^n\sum_{k=1}^{n}\binom{n}{k}\frac{x^k}{(k+y)^s}=\frac{1}{1-t}\sum_{n=1}^{\infty}\frac{1}{(n+y)^s}\left[\frac{xt}{(1-t)}\right]^n$$

See also (4.4.24a) and (4.4.43m) in Volume II(b).

**Proof:**

Following on from Lemma 3.4, let us define

$$P_s(t,x)=\sum_{n=1}^{\infty}t^n\sum_{k=1}^{n}\binom{n}{k}\frac{x^k}{k^s}$$

Writing the terms out explicitly gives us

(3.67b) $$P_s(t,x)=t\left\{\binom{1}{1}\frac{x}{1^s}\right\}$$

$$+t^2\left\{\binom{2}{1}\frac{x}{1^s}+\binom{2}{2}\frac{x^2}{2^s}\right\}$$

$$+t^3\left\{\binom{3}{1}\frac{x}{1^s}+\binom{3}{2}\frac{x^2}{2^s}+\binom{3}{3}\frac{x^3}{3^s}\right\}$$

$$+t^4\left\{\binom{4}{1}\frac{x}{1^s}+\binom{4}{2}\frac{x^2}{2^s}+\binom{4}{3}\frac{x^3}{3^s}+\binom{4}{2}\frac{x^4}{4^s}\right\}$$

$$+ \ldots\ldots\ldots\ldots\ldots\ldots\ldots\ldots\ldots\ldots$$

Assuming that the series is absolutely convergent for $\operatorname{Re}(s) > 1$, it may be rearranged by adding the terms vertically in columns as follows

$$P_s(t,x) = \frac{x}{1^s}\left\{t\binom{1}{1} + t^2\binom{2}{1} + t^2\binom{3}{1} + t^4\binom{4}{1} + \ldots\right\}$$

$$+\frac{x^2}{2^s}\left\{\quad 0 \; + t^2\binom{2}{2} + t^3\binom{3}{2} + t^4\binom{4}{2} + \ldots.\right\}$$

$$+\frac{x^3}{3^s}\left\{\quad 0 + \quad 0 \; + t^3\binom{3}{2} + t^4\binom{4}{3} + \ldots\right\}$$

$$+ \ldots\ldots\ldots\ldots\ldots\ldots\ldots\ldots\ldots\ldots$$

The series in parentheses in the j[th] row of (3.67b) may be expressed as $\sum_{i=0}^{\infty} t^i \binom{i}{j}$ where, in this case, the summation is again over the upper index. The above expression correctly gives the zero terms in the j[th] row because

$$\binom{i}{j} = 0 \;\text{ if } j > i \geq 0$$

by virtue of the definition of the binomial coefficients [75, p.155]. Hence, (3.67b) can be replaced by

(3.67c) $$P_s(t,x) = \sum_{j=1}^{\infty} \frac{x^j}{j^s}\left\{\sum_{i=0}^{\infty} t^i \binom{i}{j}\right\}$$

As shown previously in (3.61)

$$\frac{t^j}{(1-t)^{j+1}} = \sum_{i=0}^{\infty} \binom{i}{j} t^i$$

and therefore we get

$$P_s(t,x) = \sum_{j=1}^{\infty} \frac{x^j}{j^s} \frac{t^j}{(1-t)^{j+1}} = \frac{1}{1-t}\sum_{j=1}^{\infty} \frac{1}{j^s}\left[\frac{xt}{(1-t)}\right]^j = \frac{1}{1-t} Li_s\left(\frac{xt}{1-t}\right)$$

Hence we have the result

(3.67d) $$P_s(t,x)=\sum_{n=1}^{\infty}t^n\sum_{k=1}^{n}\binom{n}{k}\frac{x^k}{k^s}=\frac{1}{1-t}Li_s\left(\frac{xt}{1-t}\right)$$

and specifically we have as before $P_s(1/2,x)=P_s(x)=2Li_s(x)$.

By exactly the same method we easily find that

$$\sum_{n=1}^{\infty}t^n\sum_{k=1}^{n}\binom{n}{k}\frac{x^k}{(k+y)^s}=\frac{1}{1-t}\sum_{n=1}^{\infty}\frac{1}{(n+y)^s}\left[\frac{xt}{(1-t)}\right]^n$$

$$=\frac{1}{1-t}\sum_{n=0}^{\infty}\frac{1}{(n+y)^s}\left[\frac{xt}{(1-t)}\right]^n-\frac{1}{(1-t)y^s}$$

and therefore we have

(3.67d) $$\sum_{n=1}^{\infty}t^n\sum_{k=1}^{n}\binom{n}{k}\frac{x^k}{(k+y)^s}=\frac{1}{1-t}\Phi\left(\frac{xt}{(1-t)},s,y\right)-\frac{1}{(1-t)y^s}$$

where $\Phi(z,s,u)$ is the Hurwitz-Lerch zeta function defined in (4.4.82). Reference should also be made to (4.4.43).

The following identity [75aa] is easily derived

$$z\Phi\big(z,s,y+1\big)=\Phi\big(z,s,y\big)-\frac{1}{y^s}$$

and we therefore have

$$\frac{1}{1-t}\left[\Phi\left(\frac{xt}{(1-t)},s,y\right)-\frac{1}{y^s}\right]=\frac{xt}{(1-t)^2}\Phi\left(\frac{xt}{(1-t)},s,y+1\right)$$

Hence we obtain a more "balanced" identity

(3.67e) $$\sum_{n=1}^{\infty}t^n\sum_{k=1}^{n}\binom{n}{k}\frac{x^k}{(k+y)^s}=\frac{xt}{(1-t)^2}\Phi\left(\frac{xt}{(1-t)},s,y+1\right)$$

When $x=-1$ we may write (3.67a) as

$$(1-t)\sum_{n=1}^{\infty}t^n\sum_{k=1}^{n}\binom{n}{k}\frac{(-1)^k}{k^s}=Li_s\left(\frac{-t}{1-t}\right)$$

and from (3.11b) we have

$$Li_2\left(\frac{-t}{1-t}\right) = -\sum_{n=1}^{\infty} \frac{H_n^{(1)}}{n} t^n$$

Letting $s = 2$ gives us

(3.67f) $$(1-t)\sum_{n=1}^{\infty} t^n \sum_{k=1}^{n} \binom{n}{k} \frac{(-1)^k}{k^2} = Li_2\left(\frac{-t}{1-t}\right)$$

Equating coefficients of $t^n$ results in

$$\sum_{k=1}^{n} \binom{n}{k} \frac{(-1)^k}{k^2} - \sum_{k=1}^{n-1} \binom{n-1}{k} \frac{(-1)^k}{k^2} = -\frac{H_n^{(1)}}{n}$$

Since $\binom{n-1}{n} = 0$ we may write this as

$$\sum_{k=1}^{n} \left[\binom{n}{k} - \binom{n-1}{k}\right] \frac{(-1)^k}{k^2} = -\frac{H_n^{(1)}}{n}$$

and using $\frac{n-k}{n}\binom{n}{k} = \binom{n-1}{k}$ we obtain

$$\sum_{k=1}^{n} \left[\frac{k}{n}\binom{n}{k}\right] \frac{(-1)^k}{k^2} = -\frac{H_n^{(1)}}{n}$$

This may obviously be expressed as

$$\sum_{k=1}^{n} \binom{n}{k} \frac{(-1)^{k+1}}{k} = H_n^{(1)}$$

and hence we have an alternative derivation of Euler's identity (3.5).

Letting $s = 3$ in (3.67a) gives us

(3.67g) $$(1-t)\sum_{n=1}^{\infty} t^n \sum_{k=1}^{n} \binom{n}{k} \frac{(-1)^k}{k^3} = Li_3\left(\frac{-t}{1-t}\right)$$

and we will see in (4.4.155l) that

$$Li_3\left(\frac{-t}{1-t}\right) = -\frac{1}{2}\left(\sum_{n=1}^{\infty} \frac{H_n^{(2)}}{n} t^n + \sum_{n=1}^{\infty} \frac{\left(H_n^{(1)}\right)^2}{n} t^n\right)$$

Equating coefficients of $t^n$ results in the Flajolet and Sedgewick identity (3.16b)

$$\sum_{k=1}^{n} \binom{n}{k} \frac{(-1)^{k+1}}{k^3} = \frac{1}{2}\left(H_n^{(1)}\right)^2 + \frac{1}{2} H_n^{(2)}$$

Similarly, using (4.4.155u)

$$-Li_4\left(\frac{-t}{1-t}\right) = \frac{1}{6}\sum_{n=1}^{\infty} \frac{\left[H_n^{(1)}\right]^3}{n} t^n + \frac{1}{2}\sum_{n=1}^{\infty} \frac{H_n^{(1)} H_n^{(2)}}{n} t^n + \frac{1}{3}\sum_{n=1}^{\infty} \frac{H_n^{(3)}}{n} t^n$$

gives us the Flajolet and Sedgewick identity (3.16c)

$$-\sum_{k=1}^{n} \binom{n}{k} \frac{(-1)^k}{k^3} = \frac{1}{6}\left(H_n^{(1)}\right)^3 + \frac{1}{2} H_n^{(1)} H_n^{(2)} + \frac{1}{3} H_n^{(3)}$$

This very neatly unifies some of the foregoing analysis.

Dividing (3.67g) by $t$ and integrating results in

$$\sum_{n=1}^{\infty} \left[\frac{x^n}{n} - \frac{x^{n+1}}{n+1}\right] \sum_{k=1}^{n} \binom{n}{k} \frac{(-1)^k}{k^3} = \int_0^x Li_3\left(\frac{-t}{1-t}\right) \frac{dt}{t}$$

From (4.4.155lii) we will see that

$$\int_0^x Li_3\left(\frac{-t}{1-t}\right) \frac{dt}{t} = -\frac{1}{2}\left(\sum_{n=1}^{\infty} \frac{H_n^{(2)}}{n^2} x^n + \sum_{n=1}^{\infty} \frac{\left(H_n^{(1)}\right)^2}{n^2} x^n\right)$$

Equating coefficients of $x^n$ results in

$$\frac{1}{n}\sum_{k=1}^{n} \binom{n}{k} \frac{(-1)^k}{k^3} - \frac{1}{n}\sum_{k=1}^{n-1} \binom{n-1}{k} \frac{(-1)^k}{k^3} = -\frac{1}{2} \frac{H_n^{(2)} + \left(H_n^{(1)}\right)^2}{n^2}$$

and this again easily results in the Flajolet and Sedgewick identity (3.16b).

**Theorem 3.3:**

(3.68) $$\varsigma(2) = \sum_{n=1}^{\infty} \frac{1}{2^n} \sum_{k=1}^{n} \binom{n}{k} \frac{(-1)^{k+1}}{k^2}$$

**Proof:**

We now know that this is a direct consequence of Lemma 3.4 with $s = 2$ and $x = 1$

(and employing $Li_2(-1) = -\pi^2/12$) but, for completeness, my original proof is also shown below.

Let us consider the real part of the basic identity contained in (2.21) in the case where $p(x) = x$. We have

(3.69) $$S = \sum_{n=0}^{\infty} \frac{1}{2^n} \sum_{k=0}^{n} \int_0^{\pi/2} \binom{n}{k} x \cos 2kx \, dx = \int_0^{\pi/2} x dx = \frac{\pi^2}{8}$$

$$\int_0^{\pi/2} x \cos 2kx \, dx = \frac{\cos 2kx}{4k^2} + \frac{x \sin 2kx}{2k} \bigg|_0^{\pi/2} \qquad ,(k \ge 1)$$

$$= \frac{(-1)^k}{4k^2} - \frac{1}{4k^2} \qquad ,(k \ge 1)$$

$$= \frac{\pi^2}{8} \qquad ,(k = 0)$$

Segregating the term in $S$ for $n = 0$ we have

(3.70) $$S = \frac{\pi^2}{8} + \frac{\pi^2}{8} \sum_{n=1}^{\infty} \frac{1}{2^n} + \frac{1}{4} \sum_{n=1}^{\infty} \frac{1}{2^n} \sum_{k=1}^{n} \binom{n}{k} \frac{(-1)^k}{k^2} - \sum_{n=1}^{\infty} \frac{1}{2^n} \sum_{k=1}^{n} \binom{n}{k} \frac{1}{k^2}$$

As a result of (3.65) we have already shown that

(3.71) $$2\varsigma(2) = \sum_{n=1}^{\infty} \frac{1}{2^n} \sum_{k=1}^{n} \binom{n}{k} \frac{1}{k^2}$$

and accordingly we have

(3.72) $$-\varsigma(2) = \sum_{n=1}^{\infty} \frac{1}{2^n} \sum_{k=1}^{n} \binom{n}{k} \frac{(-1)^k}{k^2}$$

Adding (3.71) and (3.72) together we obtain

(3.73) $$\varsigma(2) = \sum_{n=1}^{\infty} \frac{1}{2^n} \sum_{k=1}^{n} \binom{n}{k} \frac{\left[1-(-1)^k\right]}{k^2}$$

There is an alternative derivation of (3.68). From the Flajolet and Sedgewick formula (3.16b) we already know that (see also Section 4)

(3.74) $$\sum_{k=1}^{n} \binom{n}{k} \frac{(-1)^k}{k^2} = S_n(2) = -\frac{1}{2}\left(H_n^{(1)}\right)^2 - \frac{1}{2} H_n^{(2)}$$

which we can use to generate the following infinite series

(3.75) $$S=\sum_{n=1}^{\infty}\frac{1}{2^n}\sum_{k=1}^{n}\binom{n}{k}\frac{(-1)^k}{k^2}=-\frac{1}{2}\sum_{n=1}^{\infty}\frac{1}{2^n}\left(H_n^{(1)}\right)^2-\frac{1}{2}\sum_{n=1}^{\infty}\frac{1}{2^n}H_n^{(2)}$$

From (3.35) we know that

(3.76) $$\sum_{n=1}^{\infty}\frac{\left(H_n^{(1)}\right)^2}{2^n}=2\log^2 2+2Li_2(1/2)=\log^2 2+\varsigma(2)$$

and from (3.34) we have

(3.77) $$\sum_{n=1}^{\infty}\frac{H_n^{(2)}}{2^n}=2Li_2(1/2)=\varsigma(2)-\log^2 2$$

Therefore we have

(3.78) $$S=-\log^2 2-Li_2(1/2)=-\frac{\pi^2}{6}=-\varsigma(2)$$

Simple algebra also highlights an interesting connection between $\varsigma(2)$ and $\log^2 2$

(3.79) $$2\varsigma(2)=\sum_{n=1}^{\infty}\frac{\left(H_n^{(1)}\right)^2}{2^n}+\sum_{n=1}^{\infty}\frac{H_n^{(2)}}{2^n}$$

(3.80) $$2\log^2 2=\sum_{n=1}^{\infty}\frac{\left(H_n^{(1)}\right)^2}{2^n}-\sum_{n=1}^{\infty}\frac{H_n^{(2)}}{2^n}$$

Referring to (3.21eii)

$$n\int_0^1(1-t)^{n-1}\log^2 t\,dt=H_n^{(2)}+\left(H_n^{(1)}\right)^2$$

we see from (3.79) that

$$2\varsigma(2)=\sum_{n=1}^{\infty}\frac{\left(H_n^{(1)}\right)^2+H_n^{(2)}}{2^n}=\sum_{n=1}^{\infty}\frac{n}{2^n}\int_0^1(1-t)^{n-1}\log^2 t\,dt$$

Since $\sum_{n=1}^{\infty}nx^{n-1}=\frac{1}{(1-x)^2}$ we have

$$\varsigma(2) = \int_0^1 \frac{\log^2 t}{(1+t)^2} dt$$

This is easily verified since

$$\int \frac{\log^2 t}{(1+t)^2} dt = \log t \left[ \frac{t \log t}{1+t} - 2\log(1+t) \right] - Li_2(-t)$$

**Theorem 3.4:**

(3.81) $$\sum_{n=1}^{\infty} \frac{H_n}{n^2 2^n} = \varsigma(3) - \frac{\pi^2}{12} \log 2$$

**Proof:**

From (3.26) we have

(3.82) $$\frac{1}{2}\log^2(1-x) + Li_2(x) = \sum_{n=1}^{\infty} \frac{H_n}{n} x^n$$

Now divide (3.82) by $x$ and integrate to obtain

$$\frac{1}{2}\int_0^{1/2} \frac{\log^2(1-x)}{x} dx + \int_0^{1/2} \frac{Li_2(x)}{x} dx = \sum_{n=1}^{\infty} \frac{H_n}{n^2 2^n}$$

Using (3.37) and (3.44) the two integrals are evaluated as

$$= \frac{1}{8}\varsigma(3) - \frac{1}{6}\log^3 2 + Li_2(1/2)$$

which simplifies to $\varsigma(3) - \frac{\pi^2}{12}\log 2$ and this completes the proof. This identity is contained, inter alia, in [30].

With hindsight, the previous sections of the paper could be written in a more logical order: notwithstanding that, I have endeavoured to adhere to the order in which my thought processes developed in real time (it makes the story more interesting, especially since I was originally convinced that a closed form formula for $\varsigma(3)$ would result!). Some of the referencing will however give readers an inkling of those areas where time travel was prevalent!

**Theorem 3.5:**

(3.83) $$\varsigma_a(s) = \sum_{n=0}^{\infty} \frac{1}{(1+\lambda)^{n+1}} \sum_{k=0}^{n} \binom{n}{k} \lambda^{n-k} \frac{(-1)^k}{(1+k)^s} \quad \text{for } \lambda > 0$$

where $s$ is a positive integer.

**Proof:**

In 2004, Amore [8] published an interesting paper on the acceleration of the convergence of series through a variational approach. He starts with the following integral representation of $\varsigma_a(s)$ (which may easily be derived as in (4.4.42))

(3.84) $$\varsigma_a(s) = \frac{(-1)^{s-1}}{\Gamma(s)} \int_0^1 \frac{\log^{s-1} x}{1+x} dx$$

and then introduces $\lambda$ as a variation parameter

(3.85) $$\varsigma_a(s) = \frac{(-1)^{s-1}}{\Gamma(s)} \int_0^1 \frac{1}{1+\lambda} \frac{\log^{s-1} x}{1-\left(\frac{\lambda-x}{1+\lambda}\right)} dx$$

The inequality $\left|\frac{\lambda - x}{1+\lambda}\right| < 1$ is satisfied in the interval $[0,1]$ provided that $\lambda > 0$ and we can then employ the binomial theorem to expand the denominator as follows

(3.85a) $$\varsigma_a(s) = \frac{(-1)^{s-1}}{\Gamma(s)} \sum_{n=0}^{\infty} \frac{1}{(1+\lambda)^{n+1}} \int_0^1 (\lambda - x)^n \log^{s-1} x \, dx$$

(3.85b) $$= \frac{(-1)^{s-1}}{\Gamma(s)} \sum_{n=0}^{\infty} \frac{1}{(1+\lambda)^{n+1}} \sum_{k=0}^{n} \binom{n}{k} (-1)^k \lambda^{n-k} \int_0^1 x^k \log^{s-1} x \, dx$$

(3.85c) $$= \frac{(-1)^{s-1}}{\Gamma(s)} \sum_{n=0}^{\infty} \frac{1}{(1+\lambda)^{n+1}} \sum_{k=0}^{n} \binom{n}{k} (-1)^{s+k+1} \lambda^{n-k} \frac{\Gamma(s)}{(1+k)^s}$$

(3.85d) $$= \sum_{n=0}^{\infty} \frac{1}{(1+\lambda)^{n+1}} \sum_{k=0}^{n} \binom{n}{k} \lambda^{n-k} \frac{(-1)^k}{(1+k)^s}$$

The transition from (3.85b) to (3.85c) is facilitated by parametric differentiation: we have

$$\int_0^1 x^k dx = \frac{1}{1+k}$$

and differentiating $s-1$ times we obtain

(3.86) $$\int_0^1 x^k \log^{s-1} x \, dx = \frac{(-1)^{s+1}(s-1)!}{(1+k)^s}$$

Whilst $\lambda$ appears explicitly in (3.85d), the series itself is independent of $\lambda$. With $\lambda = 1$ we obtain the Hasse/Sondow identity (3.11)

$$\varsigma_a(s) = \sum_{n=0}^{\infty} \frac{1}{2^{n+1}} \sum_{k=0}^{n} \binom{n}{k} \frac{(-1)^k}{(k+1)^s}$$

Using the principle of minimum sensitivity (this concept is apparently considered in [127]), Amore has numerically computed that the convergence of (3.85d) is optimised when

$$\lambda_m = 0.4494081497877167793073271775714 09...$$

It would be interesting to know if $\lambda_m$ has a meaningful closed form expression (unfortunately my untutored attempts with Plouffe's Inverter [107] and the Inverse Symbolic Calculator, http://oldweb.cecm.sfu.ca/projects/ISC/ISCmain.html , did not come up with any matches).

We may also write (3.85d) as

$$\varsigma_a(s) = \sum_{n=0}^{\infty} \frac{1}{(1+\lambda)^{n+1}} \sum_{k=0}^{n} \binom{n}{k} \lambda^{n-k} \frac{(-1)^k}{(1+k)^s} = \sum_{n=0}^{\infty} \frac{\lambda^n}{(1+\lambda)^{n+1}} \sum_{k=0}^{n} \binom{n}{k} \frac{(-1)^k}{\lambda^k (1+k)^s}$$

and with $\lambda = s$ we have provided that $s > 0$

$$\varsigma_a(s) = \sum_{n=0}^{\infty} \frac{1}{(1+s)^{n+1}} \sum_{k=0}^{n} \binom{n}{k} s^{n-k} \frac{(-1)^k}{(1+k)^s} = \frac{1}{1+s} \sum_{n=0}^{\infty} \left(\frac{s}{1+s}\right)^n \sum_{k=0}^{n} \binom{n}{k} \frac{(-1)^k}{s^k (1+k)^s}$$

Therefore we have

(3.86a)

$$\varsigma_a'(s) = -\sum_{n=0}^{\infty} \frac{n+1}{(1+s)^{n+2}} \sum_{k=0}^{n} \binom{n}{k} s^{n-k} \frac{(-1)^k}{(1+k)^s} + \sum_{n=0}^{\infty} \frac{1}{(1+s)^{n+1}} \sum_{k=0}^{n} \binom{n}{k} (n-k) s^{n-k-1} \frac{(-1)^k}{(1+k)^s}$$

$$-\sum_{n=0}^{\infty} \frac{1}{(1+s)^{n+1}} \sum_{k=0}^{n} \binom{n}{k} s^{n-k} \frac{(-1)^k \log(1+k)}{(1+k)^s}$$

$$= -\frac{1}{1+s} \sum_{n=0}^{\infty} \frac{n+1}{(1+s)^{n+1}} \sum_{k=0}^{n} \binom{n}{k} s^{n-k} \frac{(-1)^k}{(1+k)^s} + \frac{1}{s} \sum_{n=0}^{\infty} \frac{n}{(1+s)^{n+1}} \sum_{k=0}^{n} \binom{n}{k} s^{n-k} \frac{(-1)^k}{(1+k)^s}$$

$$-\frac{1}{s} \sum_{n=0}^{\infty} \frac{1}{(1+s)^{n+1}} \sum_{k=0}^{n} \binom{n}{k} k\, s^{n-k} \frac{(-1)^k}{(1+k)^s} - \sum_{n=0}^{\infty} \frac{1}{(1+s)^{n+1}} \sum_{k=0}^{n} \binom{n}{k} s^{n-k} \frac{(-1)^k \log(1+k)}{(1+k)^s}$$

We therefore obtain

(3.86b)

$$\varsigma_a'(s) = \frac{1}{s(1+s)}\sum_{n=0}^{\infty}\frac{n-s}{(1+s)^{n+1}}\sum_{k=0}^{n}\binom{n}{k}s^{n-k}\frac{(-1)^k}{(1+k)^s} - \frac{1}{s}\sum_{n=0}^{\infty}\frac{1}{(1+s)^{n+1}}\sum_{k=0}^{n}\binom{n}{k}k\,s^{n-k}\frac{(-1)^k}{(1+k)^s}$$

$$-\sum_{n=0}^{\infty}\frac{1}{(1+s)^{n+1}}\sum_{k=0}^{n}\binom{n}{k}s^{n-k}\frac{(-1)^k\log(1+k)}{(1+k)^s}$$

Now returning to the $\lambda$ formulation we have

$$\varsigma_a(s) = \sum_{n=0}^{\infty}\frac{\lambda^n}{(1+\lambda)^{n+1}}\sum_{k=0}^{n}\binom{n}{k}\frac{(-1)^k}{\lambda^k(1+k)^s}$$

$$\frac{d}{d\lambda}\varsigma_a(s) = \sum_{n=0}^{\infty}\frac{(n-\lambda)\lambda^{n-1}}{(1+\lambda)^{n+2}}\sum_{k=0}^{n}\binom{n}{k}\frac{(-1)^k}{\lambda^k(1+k)^s} - \sum_{n=0}^{\infty}\frac{\lambda^n}{(1+\lambda)^{n+1}}\sum_{k=0}^{n}\binom{n}{k}\frac{k(-1)^k}{\lambda^{k+1}(1+k)^s}$$

$$= \frac{1}{\lambda(1+\lambda)}\sum_{n=0}^{\infty}\frac{(n-\lambda)}{(1+\lambda)^{n+1}}\sum_{k=0}^{n}\binom{n}{k}\frac{\lambda^{n-k}(-1)^k}{(1+k)^s} - \frac{1}{\lambda}\sum_{n=0}^{\infty}\frac{1}{(1+\lambda)^{n+1}}\sum_{k=0}^{n}\binom{n}{k}\frac{k\lambda^{n-k}(-1)^k}{(1+k)^s}$$

Since $\frac{d}{d\lambda}\varsigma_a(s) = 0$ we obtain

(3.86c)

$$\frac{1}{(1+\lambda)}\sum_{n=0}^{\infty}\frac{(n-\lambda)}{(1+\lambda)^{n+1}}\sum_{k=0}^{n}\binom{n}{k}\frac{\lambda^{n-k}(-1)^k}{(1+k)^s} = \sum_{n=0}^{\infty}\frac{1}{(1+\lambda)^{n+1}}\sum_{k=0}^{n}\binom{n}{k}\frac{k\lambda^{n-k}(-1)^k}{(1+k)^s}$$

Letting $\lambda = 1$ in (3.86c) we get

$$\frac{1}{2}\sum_{n=0}^{\infty}\frac{(n-1)}{2^{n+1}}\sum_{k=0}^{n}\binom{n}{k}\frac{(-1)^k}{(1+k)^s} = \sum_{n=0}^{\infty}\frac{1}{2^{n+1}}\sum_{k=0}^{n}\binom{n}{k}\frac{k(-1)^k}{(1+k)^s}$$

Therefore we have

$$\frac{1}{2}\sum_{n=0}^{\infty}\frac{n}{2^{n+1}}\sum_{k=0}^{n}\binom{n}{k}\frac{(-1)^k}{(1+k)^s} - \frac{1}{2}\sum_{n=0}^{\infty}\frac{1}{2^{n+1}}\sum_{k=0}^{n}\binom{n}{k}\frac{(-1)^k}{(1+k)^s} = \sum_{n=0}^{\infty}\frac{1}{2^{n+1}}\sum_{k=0}^{n}\binom{n}{k}\frac{k(-1)^k}{(1+k)^s}$$

and hence by adding $\sum_{n=0}^{\infty}\frac{1}{2^{n+1}}\sum_{k=0}^{n}\binom{n}{k}\frac{(-1)^k}{(1+k)^s}$ to both sides we get

$$\frac{1}{2}\sum_{n=0}^{\infty}\frac{n}{2^{n+1}}\sum_{k=0}^{n}\binom{n}{k}\frac{(-1)^k}{(1+k)^s} + \sum_{n=0}^{\infty}\frac{1}{2^{n+1}}\sum_{k=0}^{n}\binom{n}{k}\frac{(-1)^k}{(1+k)^s} = \sum_{n=0}^{\infty}\frac{1}{2^{n+1}}\sum_{k=0}^{n}\binom{n}{k}\frac{(-1)^k}{(1+k)^{s-1}}$$

This is equivalent to

(3.86ci) $$\frac{1}{2}\sum_{n=0}^{\infty}\frac{n}{2^n}\sum_{k=0}^{n}\binom{n}{k}\frac{(-1)^k}{(1+k)^s}=\varsigma_a(s-1)-\varsigma_a(s)$$

Since $\lambda$ is a variable we may again replace it by a fixed $\lambda = s$ in (3.86c) and, after dividing by $s$ we get

(3.86cii) $$\frac{1}{s(s+1)}\sum_{n=0}^{\infty}\frac{(n-s)}{(1+s)^{n+1}}\sum_{k=0}^{n}\binom{n}{k}\frac{s^{n-k}(-1)^k}{(1+k)^s}=\frac{1}{s}\sum_{n=0}^{\infty}\frac{1}{(1+s)^{n+1}}\sum_{k=0}^{n}\binom{n}{k}\frac{ks^{n-k}(-1)^k}{(1+k)^s}$$

Hence we obtain from (3.86c) and (3.86cii) the result

(3.86ciii) $$\varsigma_a'(s)=-\sum_{n=0}^{\infty}\frac{s^n}{(1+s)^{n+1}}\sum_{k=0}^{n}\binom{n}{k}\frac{(-1)^k\log(1+k)}{s^k(1+k)^s}$$

and with $s=1$ this concurs with the derivative of the Hasse/Sondow formula (3.11). With, for example, $s=2$ we obtain

(3.86civ) $$\varsigma_a'(2)=-\frac{1}{3}\sum_{n=0}^{\infty}\left(\frac{2}{3}\right)^n\sum_{k=0}^{n}\binom{n}{k}\frac{(-1)^k\log(1+k)}{2^k(1+k)^2}$$

which may have more rapid convergence properties.

In a later version of his paper Amore [8] applied the same technique to the integral representation

$$\varsigma(s)=\frac{1}{1-2^{1-s}}\frac{(-1)^{s-1}}{\Gamma(s)}\int_0^1\frac{\log^{s-1}x}{1+x}dx$$

to obtain

$$\varsigma(s)=\frac{1}{1-2^{1-s}}\sum_{n=0}^{\infty}\frac{\lambda^n}{(1+\lambda)^{n+1}}\sum_{k=0}^{n}\binom{n}{k}\frac{1}{\lambda^j}\frac{(-1)^k}{(1+k)^s}$$

and in 2007 Coffey [45i] noted that this is in fact valid for all complex $s\neq 1$ and that $\lambda$ may be complex with $\mathrm{Re}(\lambda) > 0$. Therefore with $s=0$ we have

$$\varsigma(0)=-\frac{1}{1+\lambda}\sum_{n=0}^{\infty}\frac{\lambda^n}{(1+\lambda)^n}\sum_{k=0}^{n}\binom{n}{k}(-1)^k\frac{1}{\lambda^j}$$

$$=-\frac{1}{1+\lambda}\sum_{n=0}^{\infty}\frac{\lambda^n}{(1+\lambda)^n}\left(1-\frac{1}{\lambda}\right)^n$$

$$=-\frac{1}{1+\lambda}\sum_{n=0}^{\infty}\frac{(\lambda-1)^n}{(1+\lambda)^n}=-\frac{1}{2}$$

and we therefore get $\varsigma(0) = -1/2$.

In addition, as noted by Coffey [45i] we have

$$-\left[1-2^{1-s}\right]\varsigma'(s) = 2^{1-s}\varsigma(s)\log 2 + \frac{1}{1+\lambda}\sum_{n=0}^{\infty}\frac{\lambda^n}{(1+\lambda)^n}\sum_{k=0}^{n}\binom{n}{k}(-1)^k\frac{1}{\lambda^j}\frac{\log(1+k)}{(1+k)^s}$$

and this gives us

$$\varsigma'(0) = -\log 2 + \frac{1}{1+\lambda}\sum_{n=0}^{\infty}\frac{\lambda^n}{(1+\lambda)^n}\sum_{k=0}^{n}\binom{n}{k}(-1)^k\frac{1}{\lambda^j}\log(1+k)$$

With $\lambda = 1$ this becomes

$$\varsigma'(0) = -\log 2 + \sum_{n=0}^{\infty}\frac{1}{2^{n+1}}\sum_{k=0}^{n}\binom{n}{k}(-1)^k\log(1+k)$$

We will see in (4.4.113) in Volume III that

$$\sum_{n=0}^{\infty}\frac{1}{2^{n+1}}\sum_{k=0}^{n}\binom{n}{k}(-1)^k\log(1+k) = \frac{1}{2}\log\frac{\pi}{2}$$

and hence we obtain $\varsigma'(0) = -\frac{1}{2}\log(2\pi)$.
Similarly, we obtain

$$\varsigma'(-2m) = \frac{1}{2^{2m+1}-1}\frac{1}{1+\lambda}\sum_{n=0}^{\infty}\frac{\lambda^n}{(1+\lambda)^n}\sum_{k=0}^{n}\binom{n}{k}(-1)^k\frac{1}{\lambda^j}(1+k)^{2m}\log(1+k)$$

From (1.1) we may obtain another formula for $\varsigma_a'(s)$

$$\varsigma_a(s) = \varsigma(s)\left(1-2^{1-s}\right) = \sum_{k=1}^{\infty}\frac{(-1)^{k+1}}{k^s}$$

and upon differentiation we obtain

$$\varsigma_a'(s) = \varsigma'(s)\left(1-2^{1-s}\right) + 2^{1-s}\varsigma(s)\log 2 = \sum_{k=1}^{\infty}\frac{(-1)^k\log k}{k^s}$$

We previously demonstrated that for $s > 0$

(3.86d) $$\varsigma_a(s) = \frac{1}{1+s}\sum_{n=0}^{\infty}\left(\frac{s}{1+s}\right)^n\sum_{k=0}^{n}\binom{n}{k}\frac{(-1)^k}{s^k(1+k)^s}$$

This may also be proved as follows. We first of all note from (4.4.11) that

$$I(x)=\int_0^1 t^{x-1}\left(1-\frac{t}{s}\right)^n dt=\int_0^1 t^{x-1}\sum_{k=0}^{n}\binom{n}{k}(-1)^k\frac{t^k}{s^k}dt$$

$$=\sum_{k=0}^{n}\binom{n}{k}(-1)^k\frac{1}{s^k(x+k)}$$

Differentiating $s-1$ times we obtain

$$\frac{d^{s-1}}{dx^{s-1}}I(x)=\int_0^1 t^{x-1}\left(1-\frac{t}{s}\right)^n \log^{s-1}t\,dt=(-1)^{s-1}(s-1)!\sum_{k=0}^{n}\binom{n}{k}(-1)^k\frac{1}{s^k(x+k)^s}$$

and we then have

$$\frac{1}{1+s}\sum_{n=0}^{\infty}\left(\frac{s}{1+s}\right)^n\sum_{k=0}^{n}\binom{n}{k}\frac{(-1)^k}{s^k(x+k)^s}=\frac{(-1)^{s-1}}{(1+s)(s-1)!}\sum_{n=0}^{\infty}\left(\frac{s}{1+s}\right)^n\int_0^1 t^{x-1}\left(1-\frac{t}{s}\right)^n\log^{s-1}t\,dt$$

$$=\frac{(-1)^{s-1}}{(1+s)(s-1)!}\int_0^1\sum_{n=0}^{\infty}\left(\frac{s}{1+s}\right)^n\left(1-\frac{t}{s}\right)^n t^{x-1}\log^{s-1}t\,dt$$

Employing the geometric series this becomes

$$=\frac{(-1)^{s-1}}{(s+1)(s-1)!}\int_0^1\frac{t^{x-1}\log^{s-1}t}{1-\frac{s-t}{s+1}}dt$$

We therefore obtain

(3.86e)

$$\frac{1}{1+s}\sum_{n=0}^{\infty}\left(\frac{s}{1+s}\right)^n\sum_{k=0}^{n}\binom{n}{k}\frac{(-1)^k}{s^k(x+k)^s}=\frac{(-1)^{s-1}}{(s-1)!}\int_0^1\frac{t^{x-1}\log^{s-1}t}{1+t}dt=\frac{(-1)^{s-1}}{(s-1)!}\frac{d^{s-1}}{dx^{s-1}}\int_0^1\frac{t^{x-1}}{1+t}dt$$

and upon letting $x=1$ this becomes

$$\frac{1}{1+s}\sum_{n=0}^{\infty}\left(\frac{s}{1+s}\right)^n\sum_{k=0}^{n}\binom{n}{k}\frac{(-1)^k}{s^k(1+k)^s}=\frac{(-1)^{s-1}}{(s-1)!}\int_0^1\frac{\log^{s-1}t}{1+t}dt=\varsigma_a(s)$$

As noted in [75aa] we may reverse the order of summation to obtain yet another derivation of the above identity

$$\frac{1}{1+s}\sum_{n=0}^{\infty}\left(\frac{s}{1+s}\right)^n\sum_{k=0}^{n}\binom{n}{k}\frac{(-1)^k}{s^k(1+k)^s}=\frac{1}{1+s}\sum_{k=0}^{\infty}\frac{(-1)^k}{s^k(1+k)^s}\sum_{n=k}^{\infty}\binom{n}{k}\left(\frac{s}{1+s}\right)^n$$

$$=\frac{1}{1+s}\sum_{k=0}^{\infty}\frac{(-1)^k}{s^k(1+k)^s}(1+s)s^k$$

$$=\sum_{k=0}^{\infty}\frac{(-1)^k}{(1+k)^s}=\sum_{k=1}^{\infty}\frac{(-1)^{k+1}}{k^s}=\varsigma_a(s)$$

Similarly we have

$$\frac{1}{1+s}\sum_{n=0}^{\infty}\left(\frac{s}{1+s}\right)^n\sum_{k=0}^{n}\binom{n}{k}\frac{(-1)^k\log(1+k)}{s^k(1+k)^s}=\frac{1}{1+s}\sum_{k=0}^{\infty}\frac{(-1)^k\log(1+k)}{s^k(1+k)^s}\sum_{n=k}^{\infty}\binom{n}{k}\left(\frac{s}{1+s}\right)^n$$

$$=\frac{1}{1+s}\sum_{k=0}^{\infty}\frac{(-1)^k\log(1+k)}{s^k(1+k)^s}(1+s)s^k$$

$$=\sum_{k=0}^{\infty}\frac{(-1)^k\log(1+k)}{(1+k)^s}=\sum_{k=1}^{\infty}\frac{(-1)^{k+1}\log k}{k^s}=\varsigma_a'(s)$$

Letting $s=1$ in (3.86e) and integrating with respect to $x$ we get for $a,b>0$

$$\sum_{n=0}^{\infty}\frac{1}{2^{n+1}}\sum_{k=0}^{n}\binom{n}{k}(-1)^k\log\frac{b+k}{a+k}=\int_a^b dx\int_0^1\frac{t^{x-1}}{1+t}dt$$

$$=\int_0^1\frac{1}{1+t}dt\int_a^b t^{x-1}dx$$

$$=\int_0^1\frac{t^{b-1}-t^{a-1}1}{(1+t)\log t}dt$$

We have therefore shown that

(3.86h) $$\sum_{n=0}^{\infty}\frac{1}{2^{n+1}}\sum_{k=0}^{n}\binom{n}{k}(-1)^k\log\frac{b+k}{a+k}=\int_0^1\frac{t^{b-1}-t^{a-1}}{(1+t)\log t}dt$$

and we will see this integral again in (4.4.112g).

Letting $s>1$ in (3.86e), and integrating with respect to $x$, we get for $a,b>0$

(3.86i)

$$\frac{1}{1-s^2}\sum_{n=0}^{\infty}\left(\frac{s}{1+s}\right)^n\sum_{k=0}^{n}\binom{n}{k}\frac{(-1)^k}{s^k}\left[\frac{1}{(b+k)^{s-1}}-\frac{1}{(a+k)^{s-1}}\right]=\frac{(-1)^{s-1}}{(s-1)!}\int_0^1\frac{(t^{b-1}-t^{a-1})\log^{s-2}t}{1+t}dt$$

Reversing the order of summation, we get as before

$$\frac{1}{1-s^2}\sum_{n=0}^{\infty}\left(\frac{s}{1+s}\right)^n\sum_{k=0}^{n}\binom{n}{k}\frac{(-1)^k}{s^k}\left[\frac{1}{(b+k)^{s-1}}-\frac{1}{(a+k)^{s-1}}\right]=$$

$$\frac{1}{1-s^2}\sum_{k=0}^{\infty}\frac{(-1)^k}{s^k}\left[\frac{1}{(b+k)^{s-1}}-\frac{1}{(a+k)^{s-1}}\right]\sum_{n=k}^{\infty}\binom{n}{k}\left(\frac{s}{1+s}\right)^n$$

$$=\frac{1}{1-s^2}\sum_{k=0}^{\infty}\frac{(-1)^k}{s^k}\left[\frac{1}{(b+k)^{s-1}}-\frac{1}{(a+k)^{s-1}}\right](1+s)s^k$$

$$=\frac{1}{1-s}\sum_{k=0}^{\infty}(-1)^k\left[\frac{1}{(b+k)^{s-1}}-\frac{1}{(a+k)^{s-1}}\right]$$

$$=\frac{1}{1-s}\sum_{k=0}^{\infty}\frac{(-1)^k}{(b+k)^{s-1}}-\frac{1}{1-s}\sum_{k=0}^{\infty}\frac{(-1)^k}{(a+k)^{s-1}}$$

Therefore we have

$$\frac{(-1)^s}{(s-2)!}\int_0^1\frac{(t^{b-1}-t^{a-1})\log^{s-2}t}{1+t}dt=\sum_{k=0}^{\infty}\frac{(-1)^k}{(b+k)^{s-1}}-\sum_{k=0}^{\infty}\frac{(-1)^k}{(a+k)^{s-1}}$$

and this clearly indicates that

$$\frac{(-1)^s}{(s-2)!}\int_0^1\frac{t^{x-1}\log^{s-2}t}{1+t}dt=\sum_{k=0}^{\infty}\frac{(-1)^k}{(x+k)^{s-1}}+c$$

where $c$ is a constant. The constant $c$ may be determined by letting $x=1$ because we then have

$$\frac{(-1)^s}{(s-2)!}\int_0^1\frac{\log^{s-2}t}{1+t}dt=\sum_{k=0}^{\infty}\frac{(-1)^k}{(1+k)^{s-1}}+c=\varsigma_a(s-1)+c$$

Since we have from [25, p.240] for $s\ge 3$

$$\int_0^1\frac{\log^{s-2}t}{1+t}dt=\frac{(-1)^s(s-2)!(2^{s-2}-1)}{2^{s-2}}\varsigma(s-1)$$

$$=(-1)^s(s-2)!(1-2^{1-(s-1)})\varsigma(s-1)$$

(3.86ii) $$\int_0^1 \frac{\log^{s-2} t}{1+t} dt = (-1)^s (s-2)! \varsigma_a(s-1)$$

it is seen that $c = 0$. We then obtain with $s \to s+1$

(3.86j) $$\frac{(-1)^{s-1}}{(s-1)!} \int_0^1 \frac{t^{x-1} \log^{s-1} t}{1+t} dt = \sum_{k=0}^{\infty} \frac{(-1)^k}{(x+k)^s}$$

and we shall meet this integral again in (4.4.44) in the case where $x = 1$.

We will see in (C.48) of Volume VI that

(3.86ji) $$\int_0^1 \frac{t^{x-1}}{1+t} dt = \sum_{k=0}^{\infty} \frac{(-1)^k}{k+x}$$

and hence the following connection with the zeta function is easily seen

(3.86jii) $$\frac{d^{s-1}}{dx^{s-1}} \int_0^1 \frac{t^{x-1}}{1+t} dt = \int_0^1 \frac{t^{x-1} \log^{s-1} t}{1+t} dt = (-1)^{s-1} (s-1)! \sum_{k=0}^{\infty} \frac{(-1)^k}{(k+x)^s}$$

The Wolfram Integrator easily evaluates the related integrals

$$\int \frac{\log^2 t}{1+t} dt = \log^2 t \log(1+t) + 2\log t \, Li_2(-t) - 2Li_3(-t)$$

$$\int \frac{\log^3 t}{1+t} dt = 6\left[\frac{1}{6}\log^3 t \log(1+t) + \frac{1}{2}\log^2 t \, Li_2(-t) - \log t \, Li_3(-t) + Li_4(-t)\right]$$

$$\int \frac{\log^4 t}{1+t} dt = 24\left[\frac{1}{24}\log^4 t \log(1+t) + \frac{1}{6}\log^3 t \, Li_2(-t) - \frac{1}{2}\log^2 t \, Li_3(-t) + \log t \, Li_4(-t) - Li_5(-t)\right]$$

$$\int \frac{t \log^2 t}{1+t} dt = \log^2 t \log(1+t) + 2\log t \, Li_2(-t) - 2Li_3(-t) + 2t\left[1 - \log t + \frac{1}{2}\log^2 t\right]$$

Integrating (3.86ji) we get for $b > a > 0$

$$\int_a^b dx \int_0^1 \frac{t^{x-1}}{1+t} dt = \int_a^b \sum_{k=0}^{\infty} \frac{(-1)^k}{k+x} dx = \sum_{k=0}^{\infty} (-1)^k \log \frac{k+b}{k+a}$$

$$\int_a^b dx \int_0^1 \frac{t^{x-1}}{1+t} dt = \int_0^1 \frac{t^{b-1} - t^{a-1}}{(1+t)\log t} dt$$

Hence we get (see also (4.4.112g))

(3.86k) $$\int_0^1 \frac{t^{b-1} - t^{a-1}}{(1+t)\log t} dt = \sum_{k=0}^{\infty} (-1)^k \log \frac{k+b}{k+a}$$

As mentioned previously, by using elementary calculus we can also prove (3.64) in the case where $s$ is a positive integer. First of all, consider the case where $s = 1$

(3.87) $$P_1(x) = \sum_{n=1}^{\infty} \frac{1}{2^n} \sum_{k=1}^{n} \binom{n}{k} \frac{x^k}{k}$$

Assuming that term by term differentiation is valid, we obtain

(3.88) $$xP_1'(x) = \sum_{n=1}^{\infty} \frac{1}{2^n} \sum_{k=1}^{n} \binom{n}{k} x^k$$

$$= \sum_{n=1}^{\infty} \frac{1}{2^n} \left( (1+x)^n - 1 \right)$$

$$= \frac{2x}{1-x} \qquad \text{, provided } \left| \frac{1+x}{2} \right| < 1$$

By integrating (3.88) we therefore prove the formula for $s = 1$

(3.89) $$P_1(t) = 2\int_0^t \frac{dx}{1-x} = -2\log(1-t) = 2Li_1(t)$$

It is easily seen that

(3.90) $$x\frac{d}{dx}\{xP_2'(x)\} = \sum_{n=1}^{\infty} \frac{1}{2^n} \sum_{k=1}^{n} \binom{n}{k} x^k = \frac{2x}{1-x}$$

and therefore

(3.91) $$tP_2'(t) = 2\int_0^t \frac{dx}{1-x} = 2Li_1(t)$$

Hence we obtain

(3.92) $$P_2(x) = 2\int_0^x \frac{Li_1(t)}{t} dt = 2Li_2(x)$$

The case where $s = 3$ is easily seen to be

(3.93) $$x\frac{d}{dx}\left[x\frac{d}{dx}\{\,xP_3'(x)\}\right] = \frac{2x}{1-x}$$

from which we can derive

(3.94) $$P_3(x) = 2Li_3(x)$$

Using mathematical induction it can be proved that

(3.95) $$\left(x\frac{d}{dx}\right)^n P_n(x) = \frac{2x}{1-x}$$

where the symbol $\left(x\frac{d}{dx}\right)^n$ means that we apply the operator $\left(x\frac{d}{dx}\right)$ a total of $n$ times in succession, and hence we have

(3.96) $$P_n(x) = 2Li_n(x)$$

From Knopf's recent paper [88] we have

(3.97) $$\left(x\frac{d}{dx}\right)^n f(x) = \sum_{k=1}^{n} S(n,k)x^k f^{(k)}(x)$$

where $S(n,k)$ are the rather inelegantly named Stirling numbers of the second kind [126, p.58] (the brace notation $\left\{\begin{matrix} n \\ k \end{matrix}\right\}$ is sometimes employed instead of $S(n,k)$): perhaps one day we shall encounter Stirling numbers of the third kind!

These Stirling numbers are defined by various generating functions including [126, p.58]:

(3.98) $$x^n = \sum_{k=0}^{n} S(n,k)x(x-1)...(x-k+1)$$

and an equation due to Cauchy [70aa]

(3.99) $$(e^x - 1)^k = k!\sum_{n=k}^{\infty} S(n,k)\frac{x^n}{n!}$$

A proof of (3.99) was given by Póyla and Szegö in [108a, p.225].

We also have [73b]

(3.99a) $$\prod_{k=0}^{n}\frac{1}{1-kx}=\sum_{k=0}^{\infty}S(n,k)x^k$$

(3.99b) $$\prod_{j=0}^{k}\frac{x^k}{1-jx}=\sum_{n=0}^{\infty}S(n,k)x^k$$

The Scottish mathematician James Stirling (1692-1770) published his most important work, "Methodus differentalis", in 1730 and it was in this book that the theory of Stirling numbers was first developed (they were named in honour of Stirling by Nielsen); Stirling's famous asymptotic approximation of $n!$ was also contained in this book (the intimate connection between $S(n,k)$ and the Bernoulli numbers is shown in Appendix A).

The numbers $S(n,k)$ can also be expressed explicitly as [126, p.58]

(3.100) $$S(n,k)=\frac{1}{k!}\sum_{j=0}^{k}(-1)^{k-j}\binom{k}{j}j^n$$

Using (3.95) and (3.97) we see that

(3.101) $$\sum_{k=1}^{n}S(n,k)x^k Li_n^{(k)}(x)=\frac{x}{1-x}=Li_0(x)$$

Indeed, in [101b] Maximon reports that the dilogarithm satisfies the second-order inhomogeneous differential equation

$$x(1-x)Li_2''(x)+(1-x)Li_2'(x)=1$$

Letting $y=e^x$ in (3.99) we get

$$(1-y)^k=(-1)^k k!\sum_{n=k}^{\infty}S(n,k)\frac{\log^n y}{n!}$$

Dividing by $\log y$ and integrating we get

(3.101a)
$$\int_0^1\frac{(1-y)^k}{\log y}dy=(-1)^k k!\int_0^1\sum_{n=k}^{\infty}S(n,k)\frac{\log^{n-1}y}{n!}dy=(-1)^k k!\sum_{n=k}^{\infty}\frac{S(n,k)}{n!}\int_0^1\log^{n-1}y\,dy$$

Since [25, p.239] $\int_0^1\log^{n-1}y\,dy=(-1)^{n-1}(n-1)!$ we obtain

$$\int_0^1\frac{(1-y)^k}{\log y}dy=(-1)^k k!\sum_{n=k}^{\infty}(-1)^{n-1}\frac{S(n,k)}{n}$$

Making reference to (4.4.94a) in Volume III

$$\int_0^1 \frac{(1-y)^k}{\log y} dy = \sum_{j=0}^{k} \binom{k}{j} (-1)^j \log(1+j)$$

then gives us an unusual identity

$$(-1)^k k! \sum_{n=k}^{\infty} (-1)^{n-1} \frac{S(n,k)}{n} = \sum_{j=0}^{k} \binom{k}{j} (-1)^j \log(1+j)$$

For example, letting $k=1$ and noting that $S(n,1)=1$ [75, p.264] we get the familiar logarithmic series

$$\sum_{n=1}^{\infty} \frac{(-1)^{n-1}}{n} = \log 2$$

Since $S(n,2) = 2^{n-1} - 1$, for $n \geq 2$ the validity of the operation in (3.101a) must be highly questionable, leading as it does to an apparently divergent series.

An alternative proof of (3.89) is shown below.

In (4.1.6) we will show that

$$\sum_{k=1}^{n} \binom{n}{k} (-1)^{k+1} \frac{t^k}{k} = \sum_{k=1}^{n} \frac{1-(1-t)^k}{k}$$

and letting $t \to -x$ we obtain

(3.102) $$\sum_{k=1}^{n} \binom{n}{k} \frac{x^k}{k} = \sum_{k=1}^{n} \frac{(1+x)^k - 1}{k}$$

We have the following integral identity by elementary integration (I also obtained this idea from Further Mathematics [108, p.206] as I approached what I initially thought were the very closing stages of the preparation of this paper. However, that was more than two years ago and for quite a long time I must admit that I never knew when this series of papers was going to end: like Pinocchio's nose, it just grew and grew!)

(3.103) $$\int_0^x (1+y)^{k-1} dy = \frac{(1+x)^k - 1}{k}$$

and therefore we have

$$\sum_{k=1}^{n} \binom{n}{k} \frac{x^k}{k} = \sum_{k=1}^{n} \frac{(1+x)^k - 1}{k} = \sum_{k=1}^{n} \int_0^x (1+y)^{k-1} dy$$

From the definition of $P_1(x)$ we have

$$P_1(x)=\sum_{n=1}^{\infty}\frac{1}{2^n}\sum_{k=1}^{n}\binom{n}{k}\frac{x^k}{k}=\sum_{n=1}^{\infty}\frac{1}{2^n}\int_0^x\frac{(1+y)^n-1}{y}dy$$

$$=\int_0^x\frac{1}{y}\left\{\sum_{n=1}^{\infty}\frac{(1+y)^n}{2^n}-\sum_{n=1}^{\infty}\frac{1}{2^n}\right\}dy$$

$$=\int_0^x\frac{1}{y}\left\{\frac{1+y}{1-y}-1\right\}dy$$

$$=\int_0^x\frac{2}{1-y}dy$$

$$=-2\log(1-x)=2Li_1(x)$$

Yet another derivation is shown below.

We have from the binomial theorem

$$(1+xe^{-t})^n-1=\sum_{k=1}^{n}\binom{n}{k}x^k e^{-kt}$$

and therefore

$$\int_0^u\frac{(1+xe^{-t})^n-1}{x}dx=\sum_{k=1}^{n}\binom{n}{k}\frac{u^k}{k}e^{-kt}$$

Integrating with respect to $t$ we obtain

$$\int_0^{\infty}dt\int_0^u\frac{(1+xe^{t})^n-1}{x}dx=\sum_{k=1}^{n}\binom{n}{k}\frac{u^k}{k^2}$$

We now consider the summation

$$\sum_{n=1}^{\infty}\frac{1}{2^{n+1}}\sum_{k=1}^{n}\binom{n}{k}\frac{u^k}{k^2}=\int_0^{\infty}dt\int_0^u\sum_{n=1}^{\infty}\frac{(1+xe^{-t})^n-1}{x2^{n+1}}dx$$

and we have

$$\sum_{n=1}^{\infty}\frac{(1+xe^{-t})^n-1}{2^{n+1}}=\frac{1}{2}\left[\frac{1+xe^{-t}}{1-xe^{-t}}-1\right]=\frac{xe^{-t}}{1-xe^{-t}}$$

Therefore we get

$$\int_0^\infty dt \int_0^u \sum_{n=1}^\infty \frac{(1+xe^{-t})^n - 1}{x2^{n+1}} dx = \int_0^u \frac{dx}{x} \int_0^\infty \frac{xe^{-t}}{1-xe^{-t}} dt$$

$$\int_0^\infty \frac{xe^{-t}}{1-xe^{-t}} dt = \log\left(1-xe^{-t}\right)\Big|_0^\infty = -\log\left(1-x\right)$$

$$\int_0^u \frac{dx}{x} \int_0^\infty \frac{xe^{-t}}{1-xe^{-t}} dt = -\int_0^u \frac{\log(1-x)}{x} = Li_2(u)$$

Therefore we have shown that

$$Li_2(u) = \sum_{n=1}^\infty \frac{1}{2^{n+1}} \sum_{k=1}^n \binom{n}{k} \frac{u^k}{k^2}$$

Upon integration we find

$$Li_3(x) = \int_0^x \frac{Li_2(u)}{u} = \sum_{n=1}^\infty \frac{1}{2^{n+1}} \sum_{k=1}^n \binom{n}{k} \frac{x^k}{k^3}$$

and so on (which goes to show that there's more than one way to skin a cat!).

## SOME STUFF ON STIRLING NUMBERS OF THE FIRST KIND

The Stirling numbers $s(n,k)$ of the first kind [126, p.56] are defined by the following generating function (the bracket symbol $\begin{bmatrix} n \\ k \end{bmatrix}$ is also employed)

(3.104) $$x(x-1)...(x-n+1) = \sum_{k=0}^n s(n,k)x^k$$

and the Maclaurin expansion

(3.105) $$\log^k(1+x) = k!\sum_{n=k}^\infty s(n,k)\frac{x^n}{n!} \qquad |x| < 1$$

Since $s(n,k) = \frac{1}{k!}\frac{d^n}{dx^n}\log^k(1+x)\Big|_{x=0}$ it is clear that $s(n,k) = 0 \;\; \forall \; n \le k-1$ (as is also evident from the polynomial expression in (3.104)). A proof of (3.105) is given below.

The first few Stirling numbers $s(n,k)$ of the first kind are given in [120] and also the book by Srivastava and Choi [126, p.57]

(3.105i) $$s(n,0)=\delta_{n,0}$$

$$s(n,1)=(-1)^{n+1}(n-1)!$$

$$s(n,2)=(-1)^{n}(n-1)!H_{n-1}$$

$$s(n,3)=(-1)^{n+1}\frac{(n-1)!}{2}\left\{\left(H_{n-1}\right)^{2}-H_{n-1}^{(2)}\right\}$$

$$s(n,4)=(-1)^{n}\frac{(n-1)!}{6}\left\{\left(H_{n-1}\right)^{3}-3H_{n-1}H_{n-1}^{(2)}+2H_{n-1}^{(3)}\right\}$$

The above representations should be compared with the identities found by Larcombe et al. as set out, for example, in (4.4.135) ff. in Volume IV.

Kölbig [91ac] has noted that

$$\sum_{k=1}^{n}(-1)^{k}s(n,k)\binom{k}{2}2^{2-k}=(-1)^{n}\frac{(2n-1)!!}{2^{n-1}}\left[O_n^{(1)}-O_n^{(2)}\right]$$

where $O_n^{(k)}=\sum_{j=1}^{n}\frac{1}{(2j-1)^{k}}$ .

We note from (3.26) and (3.105) that

$$\frac{1}{2}\log^{2}(1-x)+Li_{2}(x)=\sum_{n=1}^{\infty}\frac{H_{n}}{n}x^{n}$$

$$\log^{k}(1-x)=k!\sum_{n=k}^{\infty}(-1)^{n}s(n,k)\frac{x^{n}}{n!}\qquad |x|<1$$

Therefore we have

$$\sum_{n=k}^{\infty}(-1)^{n}s(n,2)\frac{x^{n}}{n!}+Li_{2}(x)=\sum_{n=1}^{\infty}\frac{H_{n}}{n}x^{n}$$

and equating coefficients of $x^n$ gives us

$$(-1)^{n}\frac{s(n,2)}{n!}+\frac{1}{n^{2}}=\frac{H_{n}}{n}$$

which gives us the above formula for $s(n,2)$.

With logarithmic differentiation we obtain

$$\frac{s'(x)}{s(x)} = \frac{1}{x} + \frac{1}{x-1} ... + \frac{1}{x-n+1}$$

where $s(x) = x(x-1)...(x-n+1)$.

Multiplying across by $s(x)$ we have

$$s'(x) = (x-1)(x-2)...(x-n+1) + x(x-2)...(x-n+1) + ... + x(x-2)...(x-n)$$

and therefore we get

$$s'(0) = (-1)(-2)....(-n+1) = (-1)^{n+1}(n-1)!$$

Hence, by reference to the Maclaurin expansion for $s(x)$ we see that

(3.105ii) $\quad s(n,1) = (-1)^{n+1}(n-1)!$.

We also may write $s'(x)$ as

$$s'(x) = s_0'(x) + s_1'(x) + ... + s_n'(x)$$

and note that

$$s_0''(x) = s_0'(x)\left[\frac{1}{x-1} ... + \frac{1}{x-n+1}\right]$$

Since $s_j'(x)$ contains a factor of $x$ $\forall j \ge 1$ it is clear that $s_j''(0) = 0$ $\forall j \ge 1$. Hence we have

(3.105iii) $\quad s''(0) = s_0''(0) = (-1)^n (n-1)! H_{n-1}$

We have therefore shown that $s(n,2) = (-1)^n (n-1)! H_{n-1}$. This may also be obtained by integrating (3.29). Further Stirling numbers may be computed in the same manner, albeit requiring more algebraic dexterity (see for example Shen's paper [120] which employs the unsigned Stirling numbers $s_k^n = (-1)^{n+k} s(n,k)$). Reference should also be made to (4.3.66g) in Volume II(a) which shows how the Stirling numbers may be determined in a more systematic manner.

Letting $x = \tan^2 u$ in (3.105) we get

(3.105iv) $\quad (-1)^k 2^k \log^k \cos u = k! \sum_{n=k}^{\infty} s(n,k) \frac{\tan^{2n} u}{n!} \quad , \tan^2 u < 1$

and with $k = 1$ we see that

$$2\log\cos u = \sum_{n=k}^{\infty}(-1)^n \frac{\tan^{2n} u}{n!}$$

The Stirling numbers satisfy the recurrence relation

(3.105v) $\qquad s(n+1,m) = s(n,m-1) - ns(n,m)$ for $1 \le m \le n$

Adamchik [2] has shown that

(3.105vi) $\qquad s(n,k) = (-1)^{n-k}\dfrac{(n-1)!}{(k-1)!}w(n,k-1)$

where $w(n,0)=1$ and

$$w(n,j) = \sum_{i=0}^{j-1}\frac{\Gamma(1-j+i)}{\Gamma(1-j)}H_{n-1}^{(i+1)}w(n,j-1-i)$$

The following proof of (3.105) was given by Póyla and Szegö in [108a, p.227]. From the binomial theorem we have

$$(1-t)^{-x} = \sum_{n=0}^{\infty}\frac{x(x+1)(x+2)...(x+n-1)}{n!}t^n$$

Using (3.104) this becomes

$$= \sum_{n=0}^{\infty}\frac{(-1)^n t^n}{n!}\sum_{k=0}^{n}s(n,k)(-1)^k x^k$$

$$= 1 + \sum_{n=1}^{\infty}\frac{(-1)^n t^n}{n!}\sum_{k=0}^{n}s(n,k)(-1)^k x^k$$

$$= 1 + \sum_{n=1}^{\infty}\frac{(-1)^n t^n}{n!}\sum_{k=1}^{n}s(n,k)(-1)^k x^k$$

$$= 1 + \sum_{k=1}^{\infty}(-1)^k x^k \sum_{n=k}^{\infty}\frac{s(n,k)}{n!}(-1)^n t^n$$

We see that

$$(1-t)^{-x} = \exp[-x\log(1-t)]$$

$$= \sum_{k=1}^{\infty}\frac{x^k}{k!}\left(\log\frac{1}{1-t}\right)^k$$

and upon comparing coefficients of $x^k$ we obtain (3.105).

We may also note that

$$\frac{d^n}{dz^n} z^x = x(x-1)(x-2)...(x-n+1)z^{x-n}$$

$$\left.\frac{d^n}{dz^n} z^x\right|_{z=1} = x(x-1)(x-2)...(x-n+1)$$

We obtain from (3.105i)

$$\log^2(1+x) = 2\sum_{n=2}^{\infty}(-1)^n (n-1)!H_{n-1}\frac{x^n}{n!}$$

$$= 2\sum_{n=2}^{\infty}(-1)^n \frac{H_{n-1}}{n} x^n$$

$$= 2\sum_{n=2}^{\infty}(-1)^n \left(\frac{H_n}{n} - \frac{1}{n^2}\right) x^n$$

and, since the term corresponding to $n = 1$ is zero, we have

$$= 2\sum_{n=1}^{\infty}(-1)^n \frac{H_n}{n} x^n - 2\sum_{n=1}^{\infty}\frac{(-x)^n}{n^2}$$

With $x \to -x$ we therefore have identity (3.26)

(3.105a) $$\frac{1}{2}\log^2(1-x) + Li_2(x) = \sum_{n=1}^{\infty}\frac{H_n}{n} x^n$$

Dividing (3.105a) by $x$ and integrating we get

(3.105b) $$\frac{1}{2}\int_0^t \frac{\log^2(1-x)}{x} dx + Li_3(t) = \sum_{n=1}^{\infty}\frac{H_n}{n^2} t^n$$

The above integral was computed in (3.42) where we obtained

(3.105c) $$\int_0^t \frac{\log^2(1-x)}{x} dx = \log^2(1-t)\log t + 2\log(1-t)Li_2(1-t) - 2Li_3(1-t) + 2\varsigma(3)$$

From (3.105b) and (3.105c) we deduce

(3.105d) $$\sum_{n=1}^{\infty}\frac{H_n}{n^2} t^n = \frac{1}{2}\log^2(1-t)\log t + \log(1-t)Li_2(1-t) + Li_3(t) - Li_3(1-t) + \varsigma(3)$$

Equating the coefficients of $t^n$ in (3.105d) could be fruitful; in particular the coefficient of $t^n$ could result in a (new?) identity involving Stirling numbers of the first kind together with binomial coefficients. Unfortunately, no new information would be obtained for $\varsigma(3)$.

Combining (3.105d) with (3.21ei) we obtain

$$-\sum_{n=1}^{\infty} t^n \int_0^1 \frac{(1-x)^{n-1}\log x\, dx}{n} = \frac{1}{2}\log^2(1-t)\log t + \log(1-t)Li_2(1-t) + Li_3(t) - Li_3(1-t) + \varsigma(3)$$

We see that

$$-\sum_{n=1}^{\infty} t^n \int_0^1 \frac{(1-x)^{n-1}\log x\, dx}{n} = -\sum_{n=1}^{\infty}\int_0^1 \frac{[(1-x)t]^n \log x\, dx}{(1-x)n}$$

$$= -\int_0^1 \sum_{n=1}^{\infty} \frac{[(1-x)t]^n}{n}\frac{\log x\, dx}{1-x}$$

$$= \int_0^1 \frac{\log[(1-x)t]\log x\, dx}{1-x}$$

and hence we have

$$\int_0^1 \frac{\log[(1-x)t]\log x\, dx}{1-x} = \frac{1}{2}\log^2(1-t)\log t + \log(1-t)Li_2(1-t) + Li_3(t) - Li_3(1-t) + \varsigma(3)$$

On the other hand, a straightforward integration by parts gives us

$$\int \frac{\log[(1-x)t]\log x\, dx}{1-x} = \log[(1-x)t]Li_2(1-x) + \int \frac{Li_2(1-x)}{1-x}dx$$

$$= \log[(1-x)t]Li_2(1-x) - Li_3(1-x)$$

and the definite integral becomes

$$\int_0^1 \frac{\log[(1-x)t]\log x\, dx}{1-x} = -\varsigma(2)\log t + \varsigma(3)$$

This seems to imply that

$$\frac{1}{2}\log^2(1-t)\log t + \log(1-t)Li_2(1-t) + Li_3(t) - Li_3(1-t) = -\varsigma(2)\log t$$

but this is not correct.

I discovered in 2007 that (3.105d) is not new. As reported by Berndt [21, Part I, p.251], Ramanujan proved an equivalent form

$$\sum_{n=1}^{\infty}\frac{H_n^{(1)}}{(n+1)^2}t^{n+1}=\frac{1}{2}\log^2(1-t)\log t+Li_2(1-t)\log(1-t)-Li_3(1-t)+\varsigma(3)$$

We see that $\displaystyle\sum_{n=1}^{\infty}\frac{H_n^{(1)}}{(n+1)^2}t^{n+1}=\sum_{n=1}^{\infty}\frac{H_n^{(1)}+\frac{1}{n+1}}{(n+1)^2}t^{n+1}-\sum_{n=1}^{\infty}\frac{t^{n+1}}{(n+1)^3}$

$$=\sum_{n=1}^{\infty}\frac{H_{n+1}^{(1)}}{(n+1)^2}t^{n+1}-\sum_{n=1}^{\infty}\frac{t^{n+1}}{(n+1)^3}=\sum_{n=1}^{\infty}\frac{H_n^{(1)}}{n^2}t^n-Li_3(t)$$

and (3.105d) immediately follows.

With $t=1$ we deduce the well-known identity (see also formula (4.2.33))

(3.105e) $$\sum_{n=1}^{\infty}\frac{H_n}{n^2}=2\varsigma(3)$$

With $t=1/2$ we get

(3.105f) $$2\sum_{n=1}^{\infty}\frac{H_n}{n^2 2^n}=-\log^3 2-2\log 2Li_2(1/2)+2\varsigma(3)$$

or, alternatively, using (3.43a) we have

(3.105g) $$\sum_{n=1}^{\infty}\frac{H_n}{n^2 2^n}=\varsigma(3)-\frac{1}{2}\varsigma(2)\log 2$$

Similarly we find from (3.105) for $-1\le x<1$

(3.106) $$\frac{1}{3}\log^3(1-x)+2Li_3(x)=2\sum_{n=1}^{\infty}\frac{H_n}{n^2}x^n+\sum_{n=1}^{\infty}\frac{H_n^{(2)}}{n}x^n-\sum_{n=1}^{\infty}\frac{\left(H_n\right)^2}{n}x^n$$

In passing, we note from [25, p.76] that

$$\log^3(1-x)=\sum_{n=1}^{\infty}\frac{x^{n+2}}{n+2}\sum_{k=1}^{n}\frac{H_k}{k+1}$$

It should be noted that the above identities can also be obtained in a more laborious manner by using equations (3.35) and (3.42).

Letting $x=1/2$ in (3.106) we get

(3.106a) $$-\frac{1}{3}\log^3 2+2Li_3(1/2)=2\sum_{n=1}^{\infty}\frac{H_n}{n^2 2^n}+\sum_{n=1}^{\infty}\frac{H_n^{(2)}}{n2^n}-\sum_{n=1}^{\infty}\frac{\left(H_n\right)^2}{n2^n}$$

Reference to (3.105g) gives us an expression for $\sum_{n=1}^{\infty}\frac{H_n}{n^2 2^n}$ and looking ahead to (4.2.30b) we have

$$\varsigma_a(3)=\frac{1}{2}\sum_{n=1}^{\infty}\frac{1}{n2^n}\left\{(H_n)^2+H_n^{(2)}\right\}$$

This then provides us with two simultaneous equations for the two unknown quantities $\sum_{n=1}^{\infty}\frac{H_n^{(2)}}{n2^n}$ and $\sum_{n=1}^{\infty}\frac{\left(H_n\right)^2}{n2^n}$ which in turn give us the following identities

(3.106b) $$\sum_{n=1}^{\infty}\frac{H_n^{(2)}}{n2^n}=2Li_3(1/2)-\frac{5}{4}\varsigma(3)+\varsigma(2)\log 2-\frac{1}{3}\log^3 2$$

(3.106c) $$\sum_{n=1}^{\infty}\frac{\left(H_n\right)^2}{n2^n}=-2Li_3(1/2)+\frac{11}{4}\varsigma(3)-\varsigma(2)\log 2+\frac{1}{3}\log^3 2$$

Using (3.21eii) we see that

$$\varsigma_a(3)=\frac{1}{2}\sum_{n=1}^{\infty}\frac{1}{n2^n}\left\{(H_n)^2+H_n^{(2)}\right\}=\frac{1}{4}\sum_{n=1}^{\infty}\frac{1}{2^{n-1}}\int_0^1(1-t)^{n-1}\log^2 t\,dt$$

$$=\frac{1}{2}\int_0^1\frac{\log^2 t}{1+t}dt$$

This may be easily verified by noting that

$$\int\frac{\log^2 t}{1+t}dt=2\log(1+t)\log^2 t+2Li_2(-t)\log t-2Li_3(-t)$$

From (3.106) we see that

$$\frac{1}{3}\log^3(1-x)+2Li_3(x)=2\sum_{n=1}^{\infty}\frac{H_n}{n^2}x^n+\sum_{n=1}^{\infty}\frac{H_n^{(2)}}{n}x^n+\sum_{n=1}^{\infty}\frac{\left(H_n\right)^2}{n}x^n-2\sum_{n=1}^{\infty}\frac{\left(H_n\right)^2}{n}x^n$$

and, referring to (3.21eii), we have

$$n\int_0^1(1-t)^{n-1}\log^2 t\,dt=H_n^{(2)}+\left(H_n^{(1)}\right)^2$$

Therefore we obtain

(3.106ci) $$\sum_{n=1}^{\infty}\frac{H_n^{(2)}}{n}x^n+\sum_{n=1}^{\infty}\frac{\left(H_n\right)^2}{n}x^n=\sum_{n=1}^{\infty}\int_0^1(1-t)^{n-1}x^n\log^2 t\,dt$$

$$=\int_0^1\frac{x\log^2 t}{1-x(1-t)}dt$$

The Wolfram Integrator gives us

$$\int\frac{x\log^2 t}{1-x(1-t)}dt=\log\left(\frac{xt}{1-x}+1\right)\log^2 t+2Li_2\left(-\frac{xt}{1-x}\right)\log t-2Li_3\left(-\frac{xt}{1-x}\right)$$

and hence we get

(3.106d) $$\int_0^1\frac{x\log^2 t}{1-x(1-t)}dt=-2Li_3\left(-\frac{x}{1-x}\right)$$

Therefore we have

$$\frac{1}{3}\log^3(1-x)+2Li_3(x)=2\sum_{n=1}^{\infty}\frac{H_n}{n^2}x^n-2Li_3\left(-\frac{x}{1-x}\right)-2\sum_{n=1}^{\infty}\frac{\left(H_n\right)^2}{n}x^n$$

Using (3.105d) we then get

(3.106e)

$$2\sum_{n=1}^{\infty}\frac{\left(H_n\right)^2}{n}x^n=\log^2(1-x)\log x-\frac{1}{3}\log^3(1-x)+2\log(1-x)Li_2(1-x)-2Li_3(1-x)$$

$$+2\varsigma(3)-2Li_3\left(-\frac{x}{1-x}\right)$$

In addition from (3.106) and (3.106e) we may also deduce that

(3.106f) $$\sum_{n=1}^{\infty}\frac{H_n^{(2)}}{n}x^n=-2Li_3\left(-\frac{x}{1-x}\right)$$

Dividing (3.106ci) by $x$ and integrating results in

(3.106g) $$\sum_{n=1}^{\infty}\frac{H_n^{(2)}}{n^2}t^n+\sum_{n=1}^{\infty}\frac{\left(H_n\right)^2}{n^2}t^n=-\int_0^1\frac{\log^2 x\log[1-t(1-x)]}{1-x}dx$$

and we may compare this with (3.109b). Another similar operation gives us

(3.106h) $$\sum_{n=1}^{\infty}\frac{H_n^{(2)}}{n^3}t^n+\sum_{n=1}^{\infty}\frac{(H_n)^2}{n^3}t^n=-\int_0^1\frac{\log^2 x\,Li_2[t(1-x)]}{1-x}dx$$

Integrating (3.106) we find

(3.107) $$\frac{1}{3}\int_0^t\frac{\log^3(1-x)}{x}dx+2Li_4(t)=2\sum_{n=1}^{\infty}\frac{H_n^{(1)}}{n^3}t^n+\sum_{n=1}^{\infty}\frac{H_n^{(2)}}{n^2}t^n-\sum_{n=1}^{\infty}\frac{\left(H_n^{(1)}\right)^2}{n^2}t^n$$

With $t=1/2$ we have

(3.108) $$\frac{1}{3}\int_0^{1/2}\frac{\log^3(1-x)}{x}dx+2Li_4(1/2)=2\sum_{n=1}^{\infty}\frac{H_n^{(1)}}{2^n n^3}+\sum_{n=1}^{\infty}\frac{H_n^{(2)}}{2^n n^2}-\sum_{n=1}^{\infty}\frac{\left(H_n^{(1)}\right)^2}{2^n n^2}$$

Using the Wolfram Integrator we obtain (in a scintilla temporis)

(3.108a)

$$\int_0^t\frac{\log^3(1-x)}{x}dx=\log^3(1-t)\log t+3\log^2(1-t)Li_2(1-t)-6\log(1-t)Li_3(1-t)+6Li_4(1-t)-6\varsigma(4)$$

(a human proof is given in (3.122)) and accordingly we have

(3.108b) $$2\sum_{n=1}^{\infty}\frac{H_n^{(1)}}{n^3}t^n+\sum_{n=1}^{\infty}\frac{H_n^{(2)}}{n^2}t^n-\sum_{n=1}^{\infty}\frac{\left(H_n^{(1)}\right)^2}{n^2}t^n=$$

$$\frac{1}{3}\log^3(1-t)\log t+\log^2(1-t)Li_2(1-t)-2\log(1-t)Li_3(1-t)+2Li_4(1-t)+2Li_4(t)-2\varsigma(4)$$

With $t=1$ we have

(3.108c) $$2\sum_{n=1}^{\infty}\frac{H_n^{(1)}}{n^3}+\sum_{n=1}^{\infty}\frac{H_n^{(2)}}{n^2}-\sum_{n=1}^{\infty}\frac{\left(H_n^{(1)}\right)^2}{n^2}=0$$

From (4.4.167s) we have

$$\sum_{n=1}^{\infty}\frac{H_n^{(2)}}{n^2}=\frac{7}{4}\varsigma(4)$$

In (4.4.168) we will show that

$$\sum_{n=1}^{\infty}\frac{\left(H_n^{(1)}\right)^2}{n^2}=\frac{17}{360}\pi^4=\frac{17}{4}\varsigma(4)$$

and we therefore obtain the known result

(3.108d) $$\sum_{n=1}^{\infty}\frac{H_n^{(1)}}{n^3}=\frac{5}{4}\varsigma(4)$$

Using (3.17) we have

(3.109) $$2\sum_{n=1}^{\infty}\frac{x^n}{n}\sum_{k=1}^{n}\frac{H_k^{(1)}}{k}=\sum_{n=1}^{\infty}\frac{\left(H_n^{(1)}\right)^2}{n}x^n+\sum_{n=1}^{\infty}\frac{H_n^{(2)}}{n}x^n$$

and adding (3.106) and (3.109) together we obtain

(3.109a) $$\frac{1}{3}\log^3(1-x)+2Li_3(x)=2\left[\sum_{n=1}^{\infty}\frac{H_n^{(1)}}{n^2}x^n+\sum_{n=1}^{\infty}\frac{H_n^{(2)}}{n}x^n-\sum_{n=1}^{\infty}\frac{x^n}{n}\sum_{k=1}^{n}\frac{H_k^{(1)}}{k}\right]$$

We also note that Spieß [123bi] has derived the following identities

$$(-1)^2\log^2(1-x)=2\sum_{n=0}^{\infty}\frac{1}{n+1}H_n^{(1)}x^{n+1}$$

$$(-1)^3\log^3(1-x)=3\sum_{n=0}^{\infty}\frac{1}{n+1}\left[\left(H_n^{(1)}\right)^2-H_n^{(2)}\right]x^{n+1}$$

$$=3\sum_{n=1}^{\infty}\frac{1}{n}\left[\left(H_{n-1}^{(1)}\right)^2-H_{n-1}^{(2)}\right]x^n$$

as well as a general identity for $\log^p(1-x)$. These are in fact equivalent to a combination of (3.105) and (3.105i). We may compare the latter identity with (3.106)

$$\frac{1}{3}\log^3(1-x)+2Li_3(x)=2\sum_{n=1}^{\infty}\frac{H_n}{n^2}x^n+\sum_{n=1}^{\infty}\frac{H_n^{(2)}}{n}x^n-\sum_{n=1}^{\infty}\frac{\left(H_n\right)^2}{n}x^n$$

and obtain the trivial identity

$$2Li_3(x)=2\sum_{n=1}^{\infty}\frac{H_n^{(1)}}{n^2}x^n+\sum_{n=1}^{\infty}\frac{H_n^{(2)}}{n}x^n-\sum_{n=1}^{\infty}\frac{\left(H_n^{(1)}\right)^2}{n}x^n+\sum_{n=1}^{\infty}\frac{\left(H_{n-1}^{(1)}\right)^2}{n}x^n-\sum_{n=1}^{\infty}\frac{H_{n-1}^{(2)}}{n}x^n$$

Equating coefficients of $x^n$ in (3.110) gives us

$$\frac{1}{3}3!s(n,3)\frac{(-1)^n}{n!}+\frac{2}{n^3}=2\frac{H_n^{(1)}}{n^2}+2\frac{H_n^{(2)}}{n}-\frac{1}{n}\sum_{k=1}^{n}\frac{H_k^{(1)}}{k}$$

Using (3.105i) $s(n,3) = (-1)^{n+1}\frac{(n-1)!}{2}\left[\left(H_{n-1}\right)^2 - H_{n-1}^{(2)}\right]$ this becomes

$$-\frac{1}{n}\left[\left(H_{n-1}\right)^2 - H_{n-1}^{(2)}\right] + \frac{2}{n^3} = 2\frac{H_n^{(1)}}{n^2} + 2\frac{H_n^{(2)}}{n} - \frac{1}{n}\sum_{k=1}^{n}\frac{H_k^{(1)}}{k}$$

and a little algebra will produce Adamchik's identity (3.17).

From (3.108b) we see that

$$2\sum_{n=1}^{\infty}\frac{H_n^{(1)}}{n^3}t^n + \left[\sum_{n=1}^{\infty}\frac{H_n^{(2)}}{n^2}t^n + \sum_{n=1}^{\infty}\frac{\left(H_n^{(1)}\right)^2}{n^2}t^n\right] - 2\sum_{n=1}^{\infty}\frac{\left(H_n^{(1)}\right)^2}{n^2}t^n =$$

$$\frac{1}{3}\log^3(1-t)\log t + \log^2(1-t)Li_2(1-t) - 2\log(1-t)Li_3(1-t) + 2Li_4(1-t) + 2Li_4(t) - 2\varsigma(4)$$

Using (3.21eii) again, we see that

$$\sum_{n=1}^{\infty}\frac{H_n^{(2)}}{n^2}t^n + \sum_{n=1}^{\infty}\frac{\left(H_n^{(1)}\right)^2}{n^2}t^n = \sum_{n=1}^{\infty}\frac{1}{n}\int_0^1 (1-x)^{n-1}t^n\log^2 x\,dx$$

Using $\frac{1}{n} = \int_0^1 u^{n-1}du$ we obtain

$$\sum_{n=1}^{\infty}\frac{1}{n}\int_0^1 (1-x)^{n-1}t^n\log^2 x\,dx = \sum_{n=1}^{\infty}\int_0^1 u^{n-1}du\int_0^1 (1-x)^{n-1}t^n\log^2 x\,dx$$

$$= \int_0^1 du\int_0^1 \frac{t\log^2 x}{1-ut(1-x)}dx$$

$$= -\int_0^1 \frac{\log^2 x\log[1-t(1-x)]}{1-x}dx$$

Unfortunately the Wolfram Integrator is unable to evaluate the above integral. Our analysis however shows that

(3.109b) $$2\sum_{n=1}^{\infty}\frac{H_n^{(1)}}{n^3}t^n - \int_0^1 \frac{\log^2 x\log[1-t(1-x)]}{1-x}dx - 2\sum_{n=1}^{\infty}\frac{\left(H_n^{(1)}\right)^2}{n^2}t^n =$$

$$\frac{1}{3}\log^3(1-t)\log t + \log^2(1-t)Li_2(1-t) - 2\log(1-t)Li_3(1-t) + 2Li_4(1-t) + 2Li_4(t) - 2\varsigma(4)$$

and we may simplify this slightly by using (3.108b).

An identity may be obtained for $\int_0^1 \frac{\log^2 x\, Li_2[t(1-x)]}{1-x}dx$ by dividing (3.109b) by $t$ and integrating.

With $t=1$ we see that

(3.109c) $$2\sum_{n=1}^\infty \frac{H_n^{(1)}}{n^3} - \int_0^1 \frac{\log^3 x}{1-x}dx - 2\sum_{n=1}^\infty \frac{\left(H_n^{(1)}\right)^2}{n^2} = 0$$

With $t=1/2$ we have

(3.109d) $$2\sum_{n=1}^\infty \frac{H_n^{(1)}}{2^n n^3} - \int_0^1 \frac{\log^2 x \log[(1+x)/2]}{1-x}dx - 2\sum_{n=1}^\infty \frac{\left(H_n^{(1)}\right)^2}{2^n n^2} =$$

$$\frac{1}{3}\log^4 2 + \log^2 2 Li_2(1/2) + 2\log 2 Li_3(1/2) + 2Li_4(1/2) + 2Li_4(1/2) - 2\varsigma(4)$$

In (3.46a) we proved that

(3.110a) $$\sum_{n=1}^\infty \frac{\left(H_n^{(1)}\right)^2}{n} x^n = -\frac{1}{3}\log^3(1-x) + Li_3(x) - Li_2(x)\log(1-x)$$

We now divide (3.110a) by $x$ and integrate to obtain

(3.110b) $$\sum_{n=1}^\infty \frac{\left(H_n^{(1)}\right)^2}{n^2} t^n =$$

$$-\frac{1}{3}\left[\log^3(1-t)\log t + 3\log^2(1-t)Li_2(1-t) - 6\log(1-t)Li_3(1-t) + 6Li_4(1-t) - 6\varsigma(4)\right]$$

$$+Li_4(t) - \int_0^t \frac{Li_2(x)\log(1-x)}{x}dx$$

With $t=1$ we have

(3.110c) $$\sum_{n=1}^\infty \frac{\left(H_n^{(1)}\right)^2}{n^2} = 3\varsigma(4) - \int_0^1 \frac{Li_2(x)\log(1-x)}{x}dx$$

In (4.4.168) we will show that

(3.110d) $$\sum_{n=1}^{\infty} \frac{\left(H_n^{(1)}\right)^2}{n^2} = \frac{17}{360}\pi^4 = \frac{17}{4}\varsigma(4)$$

and hence we have

(3.110d) $$\int_0^1 \frac{Li_2(x)\log(1-x)}{x}dx = -\frac{5}{4}\varsigma(4)$$

and this should be compared with (4.4.229).

Making the substitution $y = 1 - x$ we get

$$\int_0^1 \frac{Li_2(1-y)\log y}{1-y}dy = -\frac{5}{4}\varsigma(4)$$

An alternative proof is set out below

(3.110e) $$\int_0^1 \frac{Li_2(1-y)\log y}{1-y}dy = -\frac{1}{2}\int_0^1 \frac{d}{dy}\left[Li_2(y)\right]^2 dy = -\frac{1}{2}\left[Li_2(1)\right]^2$$

$$= -\frac{1}{2}\varsigma^2(2) = -\frac{\pi^2}{72} = -\frac{5}{4}\varsigma(4)$$

Hence we have also derived (3.110d) without needing a reference to (4.4.168).

In [69a] Freitas derived the following result (the Wolfram Integrator failed to evaluate this integral)

$$\int_0^1 \frac{Li_2(y)\log y}{1-y}dy = -\frac{3}{4}\varsigma(4)$$

In the same manner we obtain the more general identity

(3.110ea) $$\sum_{n=1}^{\infty} \frac{\left(H_n^{(1)}\right)^2}{n^2}t^n = Li_4(t) + \frac{1}{2}\left[Li_2(t)\right]^2$$

$$-\frac{1}{3}\left[\log^3(1-t)\log t + 3\log^2(1-t)Li_2(1-t) - 6\log(1-t)Li_3(1-t) + 6Li_4(1-t) - 6\varsigma(4)\right]$$

Dividing (3.110ea) by $t$ and completing the customary integration would provide us with an identity for $\sum_{n=1}^{\infty} \frac{\left(H_n^{(1)}\right)^2}{n^3} t^n$, but some of the required integrals are rather forbidding.

We also obtain from (3.105)

$$\log^4(1+x) = \sum_{n=2}^{\infty} (-1)^n (n-1)! \left\{ \left(H_{n-1}\right)^3 - 3H_{n-1}H_{n-1}^{(2)} + 2H_{n-1}^{(3)} \right\} \frac{x^n}{n!}$$

$$= \sum_{n=2}^{\infty} (-1)^n \left\{ \left(H_n - \frac{1}{n}\right)^3 - 3\left(H_n - \frac{1}{n}\right)\left(H_n^{(2)} - \frac{1}{n^2}\right) + 2H_n^{(3)} - 2\frac{1}{n^3} \right\} \frac{x^n}{n}$$

$$= \sum_{n=1}^{\infty} (-1)^n \left\{ \frac{\left(H_n^{(1)}\right)^3}{n} - 3\frac{\left(H_n^{(1)}\right)^2}{n^2} + 6\frac{H_n^{(1)}}{n^3} - 6\frac{1}{n^4} - 3\frac{H_n^{(1)}H_n^{(2)}}{n} + 3\frac{H_n^{(2)}}{n^2} + 2\frac{H_n^{(3)}}{n} \right\} x^n$$

Hence we have

(3.110f)

$$\frac{1}{3}\log^4(1-x) + 2Li_4(x) = \sum_{n=1}^{\infty} \left\{ \frac{1}{3}\frac{\left(H_n^{(1)}\right)^3}{n} - \frac{\left(H_n^{(1)}\right)^2}{n^2} + 2\frac{H_n^{(1)}}{n^3} - \frac{H_n^{(1)}H_n^{(2)}}{n} + \frac{H_n^{(2)}}{n^2} + \frac{2}{3}\frac{H_n^{(3)}}{n} \right\} x^n$$

Putting $x = 1/2$ and with simple algebra we get

$$\frac{1}{3}\log^4 2 + 2Li_4(1/2) = \sum_{n=1}^{\infty} \frac{1}{2^n} \left\{ \frac{1}{3}\frac{\left(H_n^{(1)}\right)^3}{n} - \frac{H_n^{(1)}H_n^{(2)}}{n} + \frac{2}{3}\frac{H_n^{(3)}}{n} \right\} + \sum_{n=1}^{\infty} \frac{1}{2^n} \left\{ -\frac{\left(H_n^{(1)}\right)^2}{n^2} + 2\frac{H_n^{(1)}}{n^3} + \frac{H_n^{(2)}}{n^2} \right\}$$

Looking ahead to (4.2.37) we see that

$$2\varsigma_a(4) = \sum_{n=1}^{\infty} \frac{1}{n2^n} \left\{ \frac{1}{3}\left(H_n^{(1)}\right)^3 + H_n^{(1)}H_n^{(2)} + \frac{2}{3}H_n^{(3)} \right\}$$

Adding these two equations results in

$$\frac{1}{3}\log^4 2 + 2Li_4(1/2) + 2\varsigma_a(4) = 2\sum_{n=1}^{\infty} \frac{1}{2^n} \left\{ \frac{1}{3}\frac{\left(H_n^{(1)}\right)^3}{n} + \frac{2}{3}\frac{H_n^{(3)}}{n} \right\} + \sum_{n=1}^{\infty} \frac{1}{2^n} \left\{ -\frac{\left(H_n^{(1)}\right)^2}{n^2} + 2\frac{H_n^{(1)}}{n^3} + \frac{H_n^{(2)}}{n^2} \right\}$$

which may be slightly simplified to

$$\log^4 2 + 6Li_4(1/2) + 6\varsigma_a(4) = 2\sum_{n=1}^{\infty}\frac{1}{2^n}\left\{\frac{\left(H_n^{(1)}\right)^3}{n} + 2\frac{H_n^{(3)}}{n}\right\} + \sum_{n=1}^{\infty}\frac{1}{2^n}\left\{6\frac{H_n^{(1)}}{n^3} + 6\frac{H_n^{(2)}}{n^2} - 3\frac{\left(H_n^{(1)}\right)^2}{n^2}\right\}$$

From (3.108b) we have

$$\sum_{n=1}^{\infty}\frac{1}{2^n}\left\{6\frac{H_n^{(1)}}{n^3} + 6\frac{H_n^{(2)}}{n^2} - 3\frac{\left(H_n^{(1)}\right)^2}{n^2}\right\} =$$

$$\log^4 2 + 3\log^2 2Li_2(1/2) + 6\log 2Li_3(1/2) + 12Li_4(1/2) - 6\varsigma(4)$$

Therefore we obtain

$$2\sum_{n=1}^{\infty}\frac{\left(H_n^{(1)}\right)^3}{n2^n} + 4\sum_{n=1}^{\infty}\frac{H_n^{(3)}}{n2^n} = 6\varsigma_a(4) - 6\varsigma(4) - 6Li_4(1/2) - 3\log^2 2Li_2(1/2) - 6\log 2Li_3(1/2)$$

Upon division of (3.110f) by $x$ and integration we obtain

$$\frac{1}{3}\int_0^t \frac{\log^4(1-x)}{x}dx + 2Li_5(t) = \sum_{n=1}^{\infty}\left\{\frac{1}{3}\frac{\left(H_n^{(1)}\right)^3}{n^2} - \frac{\left(H_n^{(1)}\right)^2}{n^3} + 2\frac{H_n^{(1)}}{n^4} - \frac{H_n^{(1)}H_n^{(2)}}{n^2} + \frac{H_n^{(2)}}{n^3} + \frac{2}{3}\frac{H_n^{(3)}}{n^2}\right\}t^n$$

Using the Wolfram Integrator we obtain

$$\int_0^t \frac{\log^4(1-x)}{x}dx = 24\left[\begin{array}{l}\frac{1}{24}\log^4(1-t)\log t + \frac{1}{6}\log^3(1-t)Li_2(1-t) - \frac{1}{2}\log^2(1-t)Li_3(1-t) \\ \\ +\log(1-t)Li_4(1-t) - Li_5(1-t) + \varsigma(5)\end{array}\right]$$

and accordingly we have

(3.110g) $$\sum_{n=1}^{\infty}\left\{\frac{1}{3}\frac{\left(H_n^{(1)}\right)^3}{n^2} - \frac{\left(H_n^{(1)}\right)^2}{n^3} + 2\frac{H_n^{(1)}}{n^4} - \frac{H_n^{(1)}H_n^{(2)}}{n^2} + \frac{H_n^{(2)}}{n^3} + \frac{2}{3}\frac{H_n^{(3)}}{n^2}\right\}t^n$$

$$= 8\left[\begin{array}{l}\frac{1}{24}\log^4(1-t)\log t + \frac{1}{6}\log^3(1-t)Li_2(1-t) - \frac{1}{2}\log^2(1-t)Li_3(1-t) \\ \\ +\log(1-t)Li_4(1-t) - Li_5(1-t) + \varsigma(5) + \frac{1}{4}Li_5(t)\end{array}\right]$$

With $t = 1$ we have

(3.110h) $$\sum_{n=1}^{\infty}\left\{\frac{1}{3}\frac{\left(H_n^{(1)}\right)^3}{n^2}-\frac{\left(H_n^{(1)}\right)^2}{n^3}+2\frac{H_n^{(1)}}{n^4}-\frac{H_n^{(1)}H_n^{(2)}}{n^2}+\frac{H_n^{(2)}}{n^3}+\frac{2}{3}\frac{H_n^{(3)}}{n^2}\right\}=2\varsigma(5)$$

In 1994 Borwein et al. [28] proved that where

$$\sigma_h(p,q)=\sum_{n=1}^{\infty}\frac{H_{n-1}^{(p)}}{n^q}=\sum_{n=1}^{\infty}\frac{H_n^{(p)}}{n^q}-\varsigma(p+q)$$

and where $m=p+q$ is odd:

(3.110i) $$\sigma_h(p,q)=\frac{1}{2}\left[\binom{p+q}{p}-1\right]\varsigma(p+q)+\varsigma(p)\varsigma(q)$$

$$-\frac{1}{2}\sum_{j=1}^{n}\left[\binom{2j-2}{p-1}+\binom{2j-2}{q-1}\right]\varsigma(2j-1)\varsigma(p+q-2j+1)$$

Where $m$ is even:

(3.110j) $$\sigma_h(p,q)=-\frac{1}{2}\left[\binom{p+q}{p}+1\right]\varsigma(p+q)$$

$$+\frac{1}{2}\sum_{j=1}^{n}\left[\binom{2j-2}{p-1}+\binom{2j-2}{q-1}\right]\varsigma(2j-1)\varsigma(p+q-2j+1)$$

and where $\varsigma(1)$ should be interpreted as equal to nil wherever it occurs in the summation.

For example, we have

$$\sigma_h(2,3)=-\frac{11}{2}\varsigma(5)+3\varsigma(2)\varsigma(3)$$

which implies that

(3.110k) $$\sum_{n=1}^{\infty}\frac{H_n^{(2)}}{n^3}=-\frac{9}{2}\varsigma(5)+3\varsigma(2)\varsigma(3)$$

In addition we have

$$\sigma_h(3,2)=\frac{9}{2}\varsigma(5)-2\varsigma(2)\varsigma(3)$$

and this implies that

(3.110l) $\sum_{n=1}^{\infty}\frac{H_n^{(3)}}{n^2}=\frac{11}{2}\varsigma(5)-2\varsigma(2)\varsigma(3)$

Therefore we obtain

(3.110m) $\sum_{n=1}^{\infty}\frac{H_n^{(2)}}{n^3}+\sum_{n=1}^{\infty}\frac{H_n^{(3)}}{n^2}=\varsigma(5)+\varsigma(2)\varsigma(3)$

In 1998 Flajolet and Salvy [69], using impressive contour integration techniques, wrote the Borwein et al. result [28] in the following form:

For an odd weight $m = p+q$, the linear sums are reducible to zeta values as follows

(3.110n)

$$\sum_{n=1}^{\infty}\frac{H_n^{(p)}}{n^q}=\varsigma(m)\left\{\frac{1}{2}-\frac{(-1)^p}{2}\binom{m-1}{p}-\frac{(-1)^p}{2}\binom{m-1}{q}\right\}+\frac{1-(-1)^p}{2}\varsigma(p)\varsigma(q)$$

$$+(-1)^p\sum_{k=1}^{\lfloor p/2\rfloor}\binom{m-2k-1}{q-1}\varsigma(2k)\varsigma(m-2k)$$

$$+(-1)^p\sum_{k=1}^{\lfloor q/2\rfloor}\binom{m-2k-1}{p-1}\varsigma(2k)\varsigma(m-2k) \qquad ??$$

where $\varsigma(1)$ should be interpreted as equal to nil wherever it occurs in the summation. However, this formula may contain some misprints since it does not appear to produce the correct results. For example, it suggests that

$$\sum_{n=1}^{\infty}\frac{H_n^{(2)}}{n^3}=\varsigma(5)\left\{\frac{1}{2}-\frac{1}{2}\binom{4}{2}-\frac{1}{2}\binom{4}{3}\right\}+\binom{2}{2}\varsigma(2)\varsigma(3)+\binom{2}{1}\varsigma(2)\varsigma(3)$$

$$=-3\varsigma(5)+2\varsigma(2)\varsigma(3) \quad ??$$

$$\sum_{n=1}^{\infty}\frac{H_n^{(3)}}{n^2}=\varsigma(5)\left\{\frac{1}{2}+\frac{1}{2}\binom{4}{3}+\frac{1}{2}\binom{4}{2}\right\}+\varsigma(2)\varsigma(3)-\binom{2}{1}\varsigma(2)\varsigma(3)-\binom{2}{2}\varsigma(2)\varsigma(3)$$

$$=3\varsigma(5)-\varsigma(2)\varsigma(3) \quad ??$$

both of which disagree with [28], and also do not satisfy the identity (4.4.232a) which shows that

$$\sum_{n=1}^{\infty}\frac{H_n^{(p)}}{n^q}+\sum_{n=1}^{\infty}\frac{H_n^{(q)}}{n^p}=\varsigma(p)\varsigma(q)+\varsigma(p+q)$$

Similarly, the Flajolet and Salvy [69] formula suggests that

$$\sum_{n=1}^{\infty}\frac{H_n^{(1)}}{n^2}=\varsigma(3)\left\{\frac{1}{2}+\frac{1}{2}\binom{2}{1}+\frac{1}{2}\binom{2}{2}\right\}-\sum_{k=1}^{1}\binom{2-2k}{0}\varsigma(2k)\varsigma(3-2k)=\frac{3}{2}\varsigma(3)\ ??$$

whereas the correct result is $2\varsigma(3)$. The Flajolet and Salvy [69] formula is also cited by Freitas in [69a] and it is possible that some of the results quoted in the latter paper may accordingly be questionable.

Georghiou and Philippou [69c] gave the following formula

(3.110o)

$$\sum_{n=1}^{\infty}\frac{H_k^{(2)}}{k^{2n+1}}=\varsigma(2)\varsigma(2n+1)-\frac{(n+2)(2n+1)}{2}\varsigma(2n+3)+2\sum_{j=2}^{n+1}(j-1)\varsigma(2j-1)\varsigma(2n+4-2j)$$

and this gives us for $n=1$

(3.110p) $$\sum_{n=1}^{\infty}\frac{H_n^{(2)}}{n^3}=3\varsigma(2)\varsigma(3)-\frac{9}{2}\varsigma(5)$$

EULER, LANDEN AND SPENCE POLYLOGARITHM IDENTITIES

Let us now revisit Euler's identity (1.6c): after dividing this by $(1-x)$ and integrating, we obtain

$$\int_0^t\frac{\varsigma(2)}{1-x}dx=\int_0^t\frac{\log x\log(1-x)}{1-x}dx+\int_0^t\frac{Li_2(x)}{1-x}dx+\int_0^t\frac{Li_2(1-x)}{1-x}dx$$

Using (3.40) we have

$$\int_0^t\frac{\log x\log(1-x)}{1-x}dx=\log(1-x)Li_2(1-x)-Li_3(1-x)\Big|_0^t$$

$$=\log(1-t)Li_2(1-t)-Li_3(1-t)+\varsigma(3)$$

Using (3.45) we have

$$\int_0^t\frac{Li_2(x)}{1-x}dx=-Li_2(x)\log(1-x)\Big|_0^t-\int_0^t\frac{\log^2(1-x)}{x}dx$$

and, using (3.38), this becomes

(3.110pi)

$$\int_0^t \frac{Li_2(x)}{1-x}dx = -Li_2(x)\log(1-x)\Big|_0^t - \log^2(1-x)\log x\Big|_0^t - 2\int_0^t \frac{\log(1-x)\log x}{1-x}dx$$

$$= -Li_2(t)\log(1-t) - \log^2(1-t)\log t - 2\log(1-t)Li_2(1-t) + 2Li_3(1-t) - 2\varsigma(3)$$

Combining everything together, disappointingly we simply obtain Euler's identity (1.6c) yet again. A more fruitful result is however obtained in (4.4.156) of Volume III by using a slightly different method.

An elementary proof of Euler's dilogarithm identity (1.6c) is set out below.

Let $f(x)$ be defined for $x > 0$ as follows

$$f(x) = \log x \log(1-x)$$

Using L'Hôpital's rule it is clear that $\lim_{x\to 0} f(x) = 0$

Then we have

$$f'(x) = \frac{\log(1-x)}{x} - \frac{\log x}{1-x}$$

$$= -\sum_{n=1}^{\infty}\frac{x^{n-1}}{n} + \sum_{n=1}^{\infty}\frac{(1-x)^{n-1}}{n}$$

and upon integrating this series term by term we get

$$f(t) - f(a) = -\sum_{n=1}^{\infty}\frac{t^n}{n^2} - \sum_{n=1}^{\infty}\frac{(1-t)^n}{n^2} + \sum_{n=1}^{\infty}\frac{a^n}{n^2} + \sum_{n=1}^{\infty}\frac{(1-a)^n}{n^2}$$

Therefore, in the limit as $a \to 0$, we obtain a very straightforward proof of Euler's dilogarithm identity (albeit we still need to prove that term by term integration is valid)

(3.110q) $$\varsigma(2) = \log t \log(1-t) + Li_2(t) + Li_2(1-t)$$

Using partial fractions it is obvious that

$$\frac{\log(1-x)}{x(1-x)} = \frac{\log(1-x)}{1-x} + \frac{\log(1-x)}{x}$$

and integrating we obtain

$$\int_0^t \frac{\log(1-x)}{x(1-x)}dx = -\frac{1}{2}\log^2(1-t) - Li_2(t)$$

Employing the series definition of the dilogarithm it is easily seen that

$$Li_2'\left(\frac{-x}{1-x}\right) = \frac{\log(1-x)}{x(1-x)} \quad \text{and hence} \quad \int_0^t \frac{\log(1-x)}{x(1-x)}dx = Li_2\left(\frac{-t}{1-t}\right)$$

Consequently, we have derived the functional equation originally obtained by Landen [94b] some 227 years ago in 1780.

(3.111) $$Li_2\left(\frac{-t}{1-t}\right) = -\frac{1}{2}\log^2(1-t) - Li_2(t)$$

Equivalently, with $t \to -t$ we have

(3.111a) $$Li_2\left(\frac{t}{1+t}\right) = -\frac{1}{2}\log^2(1+t) - Li_2(-t)$$

With $t = 1/2$ we thereby obtain a very simple derivation of Euler's formula (3.43a)

$$Li_2(1/2) = \frac{\pi^2}{12} - \frac{1}{2}\log^2 2 = \frac{1}{2}\left[\varsigma(2) - \log^2 2\right]$$

Comparing (3.111) and (3.31) we see that

(3.111b) $$Li_2\left(\frac{-t}{1-t}\right) = -\sum_{n=1}^\infty \frac{H_n}{n}t^n \qquad Li_2\left(\frac{t}{1+t}\right) = \sum_{n=1}^\infty (-1)^{n+1}\frac{H_n}{n}t^n$$

We will also show in (4.4.155k) in Volume IV that

$$Li_3\left(\frac{-x}{1-x}\right) = -\frac{1}{2}\sum_{n=1}^\infty \frac{H_n^{(2)}}{n}x^n - \frac{1}{2}\sum_{n=1}^\infty \frac{\left(H_n^{(1)}\right)^2}{n}x^n$$

The Wolfram Integrator gives us

(3.111c)
$$\int_0^x Li_2\left(\frac{-t}{1-t}\right)\frac{dt}{t} = Li_2(x)\log x + Li_2\left(\frac{-x}{1-x}\right)\log x - \log(1-x)Li_2(1-x)$$
$$+ Li_3(1-x) - Li_3(x) - \varsigma(3)$$

and therefore using (3.111b) we have the following identity for $0 \le x \le 1$

(3.111d) $\sum_{n=1}^{\infty} \frac{H_n}{n^2} x^n = -Li_2(x)\log x - Li_2\left(\frac{-x}{1-x}\right)\log x + \log(1-x) Li_2(1-x)$

$$- Li_3(1-x) + Li_3(x) + \varsigma(3)$$

Using (3.111) and (3.42) we also have

$$-\frac{1}{2}\int_0^t \frac{\log^2(1-t)}{t} dt - \int_0^x Li_2(t)\, dt =$$

(3.111e)

$$-\frac{1}{2}\log x \log^2(1-x) - \log(1-x)\, Li_2(1-x) + Li_3(1-x) - \varsigma(3) - Li_3(x)$$

and upon equating (3.111c) and (3.111e) we get (3.111) again

(3.111f) $\quad Li_2(x) + Li_2\left(\frac{-x}{1-x}\right) = -\frac{1}{2}\log^2(1-x)$

As expected, the same identity arises by equating (3.105d) and (3.111d).

From (3.67a) we have

$$\sum_{n=1}^{\infty} t^n \sum_{k=1}^{n} \binom{n}{k} \frac{x^k}{k^s} = \frac{1}{1-t} Li_s\left(\frac{xt}{1-t}\right)$$

and hence with $x = -1$ we have

$$\sum_{n=1}^{\infty} t^n \sum_{k=1}^{n} \binom{n}{k} \frac{(-1)^k}{k^2} - \sum_{n=1}^{\infty} t^{n+1} \sum_{k=1}^{n} \binom{n}{k} \frac{(-1)^k}{k^2} = Li_2\left(\frac{-t}{1-t}\right)$$

Dividing by $t$ and integrating gives us

$$\sum_{n=1}^{\infty} \frac{t^n}{n} \sum_{k=1}^{n} \binom{n}{k} \frac{(-1)^k}{k^2} - \sum_{n=1}^{\infty} \frac{t^{n+1}}{n+1} \sum_{k=1}^{n} \binom{n}{k} \frac{(-1)^k}{k^2} = \int_0^t Li_2\left(\frac{-t}{1-t}\right)\frac{dt}{t} = -\sum_{n=1}^{\infty} \frac{H_n}{n^2} t^n$$

and with $t = 1$ we obtain

$$\sum_{n=1}^{\infty} \frac{1}{n(n+1)} \sum_{k=1}^{n} \binom{n}{k} \frac{(-1)^{k+1}}{k^2} = \sum_{n=1}^{\infty} \frac{H_n}{n^2}$$

Alternatively, using (3.114a)

$$\int_0^t \frac{1}{t(1-t)} Li_2\left(\frac{-t}{1-t}\right) dt = Li_3\left(\frac{-t}{1-t}\right)$$

we may obtain

$$\sum_{n=1}^{\infty}\frac{t^n}{n}\sum_{k=1}^{n}\binom{n}{k}\frac{(-1)^k}{k^2}=Li_3\left(\frac{-t}{1-t}\right)$$

and deduce that

$$\sum_{n=1}^{\infty}\frac{t^{n+1}}{n+1}\sum_{k=1}^{n}\binom{n}{k}\frac{(-1)^k}{k^2}=\sum_{n=1}^{\infty}\frac{H_n}{n^2}t^n-Li_3\left(\frac{-t}{1-t}\right)$$

Using (3.111)

$$Li_2\left(\frac{-t}{1-t}\right)=-\frac{1}{2}\log^2(1-t)-Li_2(t)$$

we also see that

$$(1-t)\sum_{n=1}^{\infty}t^n\sum_{k=1}^{n}\binom{n}{k}\frac{(-1)^k}{k^2}=-\frac{1}{2}\log^2(1-t)-Li_2(t)$$

Dividing the above equation by $(1-t)$ and integrating results in

$$\sum_{n=1}^{\infty}\frac{t^{n+1}}{n+1}\sum_{k=1}^{n}\binom{n}{k}\frac{(-1)^k}{k^2}=\frac{1}{6}\log^3(1-t)-\int_0^t\frac{Li_2(t)}{1-t}dt$$

From (3.110pi) we have

$$\int_0^t\frac{Li_2(t)}{1-t}dt$$

$$=-Li_2(t)\log(1-t)-\log^2(1-t)\log t-2\log(1-t)Li_2(1-t)+2Li_3(1-t)-2\varsigma(3)$$

and hence we get

$$\sum_{n=1}^{\infty}\frac{t^{n+1}}{n+1}\sum_{k=1}^{n}\binom{n}{k}\frac{(-1)^k}{k^2}=$$

$$\frac{1}{6}\log^3(1-t)+Li_2(t)\log(1-t)+\log^2(1-t)\log t+2\log(1-t)Li_2(1-t)-2Li_3(1-t)+2\varsigma(3)$$

By a simple integration we see that

$$\sum_{n=1}^{\infty}\frac{t^n}{n^2}\sum_{k=1}^{n}\binom{n}{k}\frac{(-1)^k}{k^2}=\int_0^t Li_3\left(\frac{-x}{1-x}\right)\frac{dx}{x}$$

Looking ahead to (4.4.64d) in Volume III we have

$$\sum_{n=1}^{\infty}\frac{t^n}{n^2}\sum_{k=1}^{n}\binom{n}{k}\frac{(-1)^k x^k}{k^s}=\frac{(-1)^s x}{\Gamma(s+1)}\int_0^1\frac{\log^s u\log\left[1-t(1-xu)\right]}{1-xu}du$$

and therefore we get

$$\sum_{n=1}^{\infty}\frac{t^n}{n^2}\sum_{k=1}^{n}\binom{n}{k}\frac{(-1)^k}{k^2}=\frac{1}{2}\int_0^1\frac{\log^2 u\log\left[1-t(1-u)\right]}{1-u}du$$

Letting $x=u/t$ we obtain

$$\int_0^t Li_3\left(\frac{-x}{1-x}\right)\frac{dx}{x}=\int_0^1 Li_3\left(\frac{-u}{t-u}\right)\frac{du}{u}$$

and we see that

$$\int_0^1 Li_3\left(\frac{-u}{t-u}\right)\frac{du}{u}=\frac{1}{2}\int_0^1\frac{\log^2 u\log\left[1-t(1-u)\right]}{1-u}du$$

Dividing (3.111d) by $x$ and integrating results in

$$-\sum_{n=1}^{\infty}\frac{H_n}{n^3}t^n=\int_0^t\frac{Li_2(x)\log x}{x}dx+\int_0^t Li_2\left(\frac{-x}{1-x}\right)\log x\frac{dx}{x}$$

$$-\int_0^t\frac{\log(1-x)Li_2(1-x)}{x}dx+\int_0^t\frac{Li_3(1-x)-\varsigma(3)}{x}dx-Li_4(t)$$

We may write this as

$$\text{(3.111g)}\qquad \sum_{n=1}^{\infty}\frac{H_n}{n^3}t^n=-Li_3(t)\log t+2Li_4(t)-K+L-M$$

where

$$K=\int_0^t Li_2\left(\frac{-x}{1-x}\right)\log x\frac{dx}{x}$$

$$L=\int_0^t\frac{\log(1-x)Li_2(1-x)}{x}dx$$

$$M = \int_0^t \frac{Li_3(1-x) - \varsigma(3)}{x} dx$$

Integration by parts gives us

$$\int_0^t Li_2\left(\frac{-x}{1-x}\right) \log x \frac{dx}{x} =$$

$$\log t \left[ Li_2(t) \log x + Li_2\left(\frac{-t}{1-t}\right) \log t - \log(1-t) Li_2(1-t) + Li_3(1-t) - Li_3(t) - \varsigma(3) \right]$$

$$-\int_0^t \left[ Li_2(x) \log x + Li_2\left(\frac{-x}{1-x}\right) \log x - \log(1-x) Li_2(1-x) + Li_3(1-x) - Li_3(x) - \varsigma(3) \right] \frac{dx}{x}$$

and hence we see that

$$2\int_0^t Li_2\left(\frac{-x}{1-x}\right) \log x \frac{dx}{x} =$$

$$\log t \left[ Li_2(t) \log x + Li_2\left(\frac{-t}{1-t}\right) \log t - \log(1-t) Li_2(1-t) + Li_3(1-t) - Li_3(t) - \varsigma(3) \right]$$

$$-\int_0^t \left[ Li_2(x) \log x - \log(1-x) Li_2(1-x) + Li_3(1-x) - Li_3(x) - \varsigma(3) \right] \frac{dx}{x}$$

We thus have

(3.111h)

$$2K = \log t \left[ Li_2(t) \log x + Li_2\left(\frac{-t}{1-t}\right) \log t - \log(1-t) Li_2(1-t) + Li_3(1-t) - Li_3(t) - \varsigma(3) \right]$$

$$-Li_3(t) \log t + 2Li_4(t) + L - M$$

$$= A - Li_3(t) \log t + 2Li_4(t) + L - M$$

where, for convenience, we have designated

$$A = \log t \left[ Li_2(t) \log t + Li_2\left(\frac{-t}{1-t}\right) \log t - \log(1-t) Li_2(1-t) + Li_3(1-t) - Li_3(t) - \varsigma(3) \right]$$

Therefore from (3.111g) and (3.111h) we have

$$K = \sum_{n=1}^{\infty} \frac{H_n}{n^3} t^n + A$$

$$L = 2\sum_{n=1}^{\infty} \frac{H_n}{n^3} t^n + A + Li_3(t)\log t - 2Li_4(t) + M$$

We note from (3.121b) that

$$M = \int_0^t \frac{\varsigma(3) - Li_3(1-x)}{x} dx = \left[\varsigma(3) - Li_3(1-t)\right]\log t - \frac{1}{2}\left[Li_2(1-t)\right]^2 + \frac{1}{2}\varsigma^2(2)$$

and we accordingly obtain

(3.111i) $$K = \int_0^t Li_2\left(\frac{-x}{1-x}\right)\log x \frac{dx}{x} = \sum_{n=1}^{\infty} \frac{H_n}{n^3} t^n$$

$$+\log t\left[Li_2(t)\log t + Li_2\left(\frac{-t}{1-t}\right)\log t - \log(1-t)Li_2(1-t) + Li_3(1-t) - Li_3(t) - \varsigma(3)\right]$$

In particular we get (which may also be obtained by using (3.11c))

(3.111j) $$\int_0^1 Li_2\left(\frac{-x}{1-x}\right)\log x \frac{dx}{x} = \sum_{n=1}^{\infty} \frac{H_n}{n^3} = \frac{5}{4}\varsigma(4)$$

We also see that

$$L = \int_0^t \frac{\log(1-x)Li_2(1-x)}{x} dx =$$

$$2\sum_{n=1}^{\infty} \frac{H_n}{n^3} t^n + \log t\left[Li_2(t)\log t + Li_2\left(\frac{-t}{1-t}\right)\log t - \log(1-t)Li_2(1-t) + Li_3(1-t) - Li_3(t) - \varsigma(3)\right]$$

$$+Li_3(t)\log t - 2Li_4(t) + \left[\varsigma(3) - Li_3(1-t)\right]\log t - \frac{1}{2}\left[Li_2(1-t)\right]^2 + \frac{1}{2}\varsigma^2(2)$$

and this simplifies to

(3.111k)

$$L = \int_0^t \frac{\log(1-x)Li_2(1-x)}{x} dx$$

$$= 2\sum_{n=1}^{\infty}\frac{H_n}{n^3}t^n + \log t\left[Li_2(t)\log t + Li_2\left(\frac{-t}{1-t}\right)\log t - \log(1-t)Li_2(1-t)\right]$$

$$-2Li_4(t) - \frac{1}{2}\left[Li_2(1-t)\right]^2 + \frac{1}{2}\varsigma^2(2)$$

In particular we get

(3.111l)

$$\int_0^1 \frac{\log(1-x)Li_2(1-x)}{x}dx = 2\sum_{n=1}^{\infty}\frac{H_n}{n^3} - 2\varsigma(4) + \frac{1}{2}\varsigma^2(2) = \frac{11}{4}\varsigma(4)$$

The Wolfram Integrator gives us a result containing a negative logarithm

$$\int_0^x Li_2\left(\frac{t}{1+t}\right)\frac{dt}{t} = \frac{1}{2}\left[\log x - \log(-x)\right]\log^2(1+x) - Li_2(1+x)\log(1+x) + \log x\, Li_2(-x)$$

$$+\log x\, Li_2\left(\frac{x}{1+x}\right) - Li_3(-x) + Li_3(1+x)$$

$$= \sum_{n=1}^{\infty}(-1)^{n+1}\frac{H_n}{n^2}x^n$$

where in the final part we used (3.111b).

From [126, p.107] we have for $t > 0$

$$Li_2(t) + Li_2\left(\frac{t}{t-1}\right) = \frac{\pi^2}{2} - \frac{1}{2}\log^2(t-1) - 2i\pi\log t + i\pi\log(t-1)$$

and letting $t = 1+x$ for $x > -1$ we get

$$Li_2(1+x) + Li_2\left(\frac{1+x}{x}\right) = \frac{\pi^2}{2} - \frac{1}{2}\log^2 x - 2i\pi\log(1+x) + i\pi\log x$$

However it is not clear to me if anything particularly meaningful can be derived from this.

Dividing Euler's dilogarithm identity by $(1-x)$ we obtain (provided $x \neq 1$)

(3.112) $$-\frac{\varsigma(2)}{1-x} + \log x\frac{\log(1-x)}{1-x} + \frac{Li_2(1-x)}{1-x} + \frac{Li_2(x)}{1-x} = 0$$

Using Landen's identity (3.111) to substitute for the last term in (3.112) we have

$$(3.113) \qquad -\frac{\varsigma(2)}{1-x}+\log x\frac{\log(1-x)}{1-x}+\frac{Li_2(1-x)}{1-x}-\frac{1}{1-x}Li_2\left(\frac{-x}{1-x}\right)-\frac{1}{2}\frac{\log^2(1-x)}{1-x}=0$$

Judiciously adding and subtracting the term $\frac{1}{2x}\log^2(1-x)$ to (3.113) we obtain

$$-\frac{\varsigma(2)}{1-x}+\left[\log x\frac{\log(1-x)}{1-x}-\frac{1}{2x}\log^2(1-x)\right]+\frac{Li_2(1-x)}{1-x}-\frac{1}{2}\frac{\log^2(1-x)}{1-x}$$

$$-\frac{1}{1-x}Li_2\left(\frac{-x}{1-x}\right)+\frac{1}{2x}\log^2(1-x)=0$$

Integrating the above equation we get

$$(3.114) \qquad \varsigma(2)\log(1-t)-\frac{1}{2}\log t\log^2(1-t)-Li_3(1-t)+\varsigma(3)+\frac{1}{6}\log^3(1-t)=$$

$$\int_0^t\left[\frac{1}{1-x}Li_2\left(\frac{-x}{1-x}\right)-\frac{1}{2x}\log^2(1-x)\right]dx$$

With partial fractions we see that

$$\frac{1}{1-x}=\frac{1}{x(1-x)}-\frac{1}{x}$$

and hence

$$\frac{1}{1-x}Li_2\left(\frac{-x}{1-x}\right)=\frac{1}{x(1-x)}Li_2\left(\frac{-x}{1-x}\right)-\frac{1}{x}Li_2\left(\frac{-x}{1-x}\right)$$

Using (3.111) this becomes

$$\frac{1}{1-x}Li_2\left(\frac{-x}{1-x}\right)=\frac{1}{x(1-x)}Li_2\left(\frac{-x}{1-x}\right)+\frac{1}{2x}\log^2(1-x)+\frac{Li_2(x)}{x}$$

and therefore we have

$$\frac{1}{1-x}Li_2\left(\frac{-x}{1-x}\right)-\frac{1}{2x}\log^2(1-x)=\frac{1}{x(1-x)}Li_2\left(\frac{-x}{1-x}\right)+\frac{Li_2(x)}{x}$$

Hence we deduce that

$$\int_0^t \left[\frac{1}{1-x}Li_2\left(\frac{-x}{1-x}\right)-\frac{1}{2x}\log^2(1-x)\right]dx=\int_0^t \frac{1}{x(1-x)}Li_2\left(\frac{-x}{1-x}\right)dx+Li_3(t)$$

Employing the series definition of the trilogarithm it is easily seen that

$$Li_3'\left(\frac{-x}{1-x}\right)=\frac{1}{x(1-x)}Li_2\left(\frac{-x}{1-x}\right)$$

and accordingly we have

(3.114a) $$\int_0^t \frac{1}{x(1-x)}Li_2\left(\frac{-x}{1-x}\right)dx=Li_3\left(\frac{-t}{1-t}\right)$$

We therefore end up with Landen's functional equation for the trilogarithm which is valid for $0 \le t < 1$

(3.115)

$$Li_3\left(\frac{-t}{1-t}\right)=\varsigma(2)\log(1-t)-\frac{1}{2}\log t\log^2(1-t)-Li_3(1-t)+\varsigma(3)+\frac{1}{6}\log^3(1-t)-Li_3(t)$$

(and this corrects the misprint in, of all places, Lewin's survey [101, p.2]). The identity is also referred to as Spence's formula [123b], and G.N. Watson gave an elementary proof of this in 1928 using a version of the integral representation of the trilogarithm given in (4.4.37).

With $t=1/2$ and using (4.4.67a) we can readily derive Landen's formula (3.43b) for the trilogarithm $Li_3(1/2)$

(3.115a) $$Li_3(1/2)=\frac{7}{8}\varsigma(3)-\frac{\pi^2}{12}\log 2+\frac{1}{6}\log^3 2$$

We now recall (3.67a)

$$\sum_{n=1}^\infty t^n\sum_{k=1}^n \binom{n}{k}\frac{x^k}{k^s}=\frac{1}{1-t}Li_s\left(\frac{xt}{1-t}\right)$$

and with $x=-1$ we have

$$\sum_{n=1}^\infty t^n\sum_{k=1}^n \binom{n}{k}\frac{(-1)^k}{k^s}=\frac{1}{1-t}Li_s\left(\frac{-t}{1-t}\right)$$

Dividing this by $t$ and integrating results in

$$\sum_{n=1}^{\infty}\frac{u^n}{n}\sum_{k=1}^{n}\binom{n}{k}\frac{(-1)^k}{k^s}=\int_0^u\frac{1}{t(1-t)}Li_s\left(\frac{-t}{1-t}\right)dt$$

Comparing this with (3.114a) results in

(3.115b) $$\sum_{n=1}^{\infty}\frac{u^n}{n}\sum_{k=1}^{n}\binom{n}{k}\frac{(-1)^k}{k^2}=Li_3\left(\frac{-u}{1-u}\right)=-\frac{1}{2}\sum_{n=1}^{\infty}\frac{u^n}{n}\left[\left(H_n^{(1)}\right)^2+H_n^{(2)}\right]$$

Using (3.16b) we also see that

(3.115c) $$\sum_{n=1}^{\infty}\frac{u^n}{n^2}\sum_{k=1}^{n}\binom{n}{k}\frac{(-1)^k}{k^2}=-\frac{1}{2}\sum_{n=1}^{\infty}\frac{u^n}{n^2}\left[\left(H_n^{(1)}\right)^2+H_n^{(2)}\right]=\int_0^u Li_3\left(\frac{-u}{1-u}\right)\frac{du}{u}$$

The Wolfram Integrator was not able to compute the above integral. See also (4.4.64d) and (4.4.155l) in Volume IV.

We also see from (3.123) that

$$\sum_{n=1}^{\infty}\frac{u^n}{n}\sum_{k=1}^{n}\binom{n}{k}\frac{(-1)^k}{k^s}=\int_0^u\frac{1}{t(1-t)}Li_s\left(\frac{-t}{1-t}\right)dt=Li_{s+1}\left(\frac{-u}{1-u}\right)$$

and hence we have

(3.115d) $$\sum_{n=1}^{\infty}\frac{u^n}{n}\sum_{k=1}^{n}\binom{n}{k}\frac{(-1)^k}{k^s}=(1-u)\sum_{n=1}^{\infty}u^n\sum_{k=1}^{n}\binom{n}{k}\frac{(-1)^k}{k^{s+1}}$$

Dividing the following by $x$ and integrating with respect to $x$

$$\sum_{n=1}^{\infty}t^n\sum_{k=1}^{n}\binom{n}{k}\frac{x^k}{k^s}=\frac{1}{1-t}Li_s\left(\frac{xt}{1-t}\right)$$

gives us

$$(1-t)\sum_{n=1}^{\infty}t^n\sum_{k=1}^{n}\binom{n}{k}\frac{u^k}{k^{s+1}}=\int_0^u Li_s\left(\frac{xt}{1-t}\right)dx=Li_{s+1}\left(\frac{ut}{1-t}\right)$$

With $u=-1$ we get

$$(1-t)\sum_{n=1}^{\infty}t^n\sum_{k=1}^{n}\binom{n}{k}\frac{(-1)^k}{k^{s+1}}=Li_{s+1}\left(\frac{-t}{1-t}\right)$$

and thereby obtain another derivation of (3.115d).

Dividing (3.115) by $1-t$ and integrating we get

$$\int_0^x \frac{1}{1-t} Li_3\left(\frac{-t}{1-t}\right) dt =$$

$$\int_0^x \left[\varsigma(2)\log(1-t) - \frac{1}{2}\log t \log^2(1-t) \; - Li_3(1-t) + \varsigma(3) + \frac{1}{6}\log^3(1-t) - Li_3(t)\right] \frac{dt}{1-t}$$

As shown below we obtain

(3.116)

$$\int_0^x \frac{1}{1-t} Li_3\left(\frac{-t}{1-t}\right) dt =$$

$$= -\frac{1}{2}\varsigma(2)\log^2(1-x) - \frac{1}{2}\log^2(1-x) Li_2(1-x) + \log(1-x) Li_3(1-x) - \varsigma(3)\log(1-x) - \frac{1}{24}\log^4(1-x)$$

$$+ Li_3(x)\log(1-x) + \frac{1}{2}\left[Li_2(x)\right]^2$$

In evaluating the above, we first of all used integration by parts to determine the following integrals.

$$\int_0^x \frac{\log t}{1-t} \log^2(1-t)\, dt = \log^2(1-x) Li_2(1-x) + 2\int_0^x Li_2(1-t) \frac{\log(1-t)}{1-t} dt$$

In the same way we get

$$\int_0^x \frac{Li_2(1-t)}{1-t} \log(1-t)\, dt = -Li_3(1-x)\log(1-x) - \int_0^x \frac{Li_3(1-t)}{1-t} dt$$

$$= -Li_3(1-x)\log(1-x) + Li_4(1-x) - \varsigma(4)$$

and therefore we have

(3.117)

$$\int_0^x \frac{\log t \log^2(1-t)}{1-t} dt = 2\left[\frac{1}{2}\log^2(1-x) Li_2(1-x) - \log(1-x) Li_3(1-x) + Li_4(1-x) - \varsigma(4)\right]$$

Using integration by parts we have

$$\int_0^x \frac{Li_3(t)}{1-t} dt = -Li_3(x)\log(1-x) + \int_0^x \frac{\log(1-t) Li_2(t)}{t} dt$$

Furthermore we get

$$\int_0^x \frac{\log(1-t)Li_2(t)}{t}dt = -[Li_2(x)]^2 - \int_0^x \frac{\log(1-t)Li_2(t)}{t}dt$$

Hence we have

(3.118) $$\int_0^x \frac{\log(1-t)Li_2(t)}{t}dt = -\frac{1}{2}[Li_2(x)]^2$$

and therefore we get

(3.118a) $$\int_0^x \frac{Li_3(t)}{1-t}dt = -Li_3(x)\log(1-x) - \frac{1}{2}[Li_2(x)]^2$$

This time we now divide (3.115) by $t$ and integrate to get

(3.119)

$$\int_0^x \frac{1}{t}Li_3\left(\frac{-t}{1-t}\right)dt =$$

$$\int_0^x \left[\varsigma(2)\log(1-t) - \frac{1}{2}\log t\log^2(1-t) - Li_3(1-t) + \varsigma(3) + \frac{1}{6}\log^3(1-t) - Li_3(t)\right]\frac{dt}{t}$$

$$= -\varsigma(2)Li_2(x) - \frac{1}{2}\int_0^x \frac{\log t\log^2(1-t)}{t}dt + \int_0^x \frac{\varsigma(3) - Li_3(1-t)}{t}dt + \frac{1}{6}\int_0^x \frac{\log^3(1-t)}{t}dt - Li_4(x)$$

We deal with each of the three integrals in turn. Let us consider the most difficult one first: reference to (4.4.167p) shows that

(3.120)

$$\int_0^x \frac{\log t\log^2(1-t)}{t}dt = \frac{1}{2}\log^2(1-x)\log^2 x + \frac{1}{12}\log^4 x - \log^2(1-x)\log^2 x + \frac{2}{3}\log(1-x)\log^3 x$$

$$-\left[\log(1-x) + \frac{1}{3}\log x\right]\log^2 x\log\left(\frac{x}{1-x}\right) + \frac{1}{2}\log^2 x\left[\log\left(\frac{x}{1-x}\right)\right]^2$$

$$-\frac{1}{4}\left[\log\left(\frac{x}{1-x}\right)\right]^4 + \log^2(1-x)Li_2(1-x) - \log^2 x\, Li_2(x)$$

$$-\left[\log\left(\frac{x}{1-x}\right)\right]^2 Li_2\left(\frac{-x}{1-x}\right) - 2\log(1-x)Li_3(1-x) + 2\log x\, Li_3(x)$$

$$+2\log\left(\frac{x}{1-x}\right)Li_3\left(\frac{-x}{1-x}\right)+2\left[Li_4(1-x)-Li_4(x)-Li_4\left(\frac{-x}{1-x}\right)\right]-2\varsigma(4)$$

We have using integration by parts

(3.121) $$\int_0^x \frac{\varsigma(3)-Li_3(1-t)}{t}\,dt=\left[\varsigma(3)-Li_3(1-x)\right]\log x-\int_0^x \frac{\log t\, Li_2(1-t)}{1-t}\,dt$$

The latter integral reduces to

$$\int_0^x \frac{\log t\, Li_2(1-t)}{1-t}\,dt=-\left[Li_2(1-x)\right]^2+\varsigma^2(2)-\int_0^x \frac{\log t\, Li_2(1-t)}{1-t}\,dt$$

and a further integration by parts gives us

(3.121a) $$\int_0^x \frac{\log t\, Li_2(1-t)}{1-t}\,dt=-\frac{1}{2}\left[Li_2(1-x)\right]^2+\frac{1}{2}\varsigma^2(2)$$

Accordingly we obtain

(3.121b) $$\int_0^x \frac{\varsigma(3)-Li_3(1-t)}{t}\,dt=\left[\varsigma(3)-Li_3(1-x)\right]\log x-\frac{1}{2}\left[Li_2(1-x)\right]^2+\frac{1}{2}\varsigma^2(2)$$

Using integration by parts we obtain

$$\int \frac{\log^3(1-x)}{x}\,dx=\log^3(1-x)\log x+3\int \frac{\log^2(1-x)\log x}{1-x}\,dx$$

$$\int \frac{\log^2(1-x)\log x}{1-x}\,dx=\log^2(1-x)\,Li_2(1-x)+2\int \frac{\log(1-x)Li_2(1-x)}{1-x}\,dx$$

$$\int \frac{\log(1-x)Li_2(1-x)}{1-x}\,dx=-\log(1-x)\,Li_3(1-x)-\int \frac{Li_3(1-x)}{1-x}\,dx$$

$$=-\log(1-x)\,Li_3(1-x)+Li_4(1-x)$$

and hence we get

(3.122)

$$\int_0^x \frac{\log^3(1-t)}{t}\,dt=6\left[\frac{1}{6}\log^3(1-x)\log x+\frac{1}{2}\log^2(1-x)\,Li_2(1-x)-\log(1-x)\,Li_3(1-x)+Li_4(1-x)-\varsigma(4)\right]$$

As before, with partial fractions we see that

$$\frac{1}{1-t}Li_3\left(\frac{-t}{1-t}\right)+\frac{1}{t}Li_3\left(\frac{-t}{1-t}\right)=\frac{1}{t(1-t)}Li_3\left(\frac{-t}{1-t}\right)$$

Therefore, upon integration we obtain

(3.123)
$$\int_0^x\left[\frac{1}{1-t}Li_3\left(\frac{-t}{1-t}\right)+\frac{1}{t}Li_3\left(\frac{-t}{1-t}\right)\right]dt=\int_0^x\frac{1}{t(1-t)}Li_3\left(\frac{-t}{1-t}\right)dt=Li_4\left(\frac{-x}{1-x}\right)$$

because, using the series definition of a polylogarithm, it is easy to prove that

$$Li_4'\left(\frac{-x}{1-x}\right)=\frac{1}{x(1-x)}Li_3\left(\frac{-x}{1-x}\right).$$

Therefore, using (3.116), (3.119) and (3.123), we obtain a functional equation involving several polylogarithms (note that the terms involving the tetralogarithm $Li_4\left(\frac{-x}{1-x}\right)$ cancel out).

(3.124)

$$-\frac{1}{2}\varsigma(2)\log^2(1-x)-\frac{1}{2}\log^2(1-x)Li_2(1-x)+\log(1-x)Li_3(1-x)-\varsigma(3)\log(1-x)-\frac{1}{24}\log^4(1-x)$$

$$+Li_3(x)\log(1-x)+\frac{1}{2}\left[Li_2(x)\right]^2-\varsigma(2)Li_2(x)$$

$$-\frac{1}{4}\log^2(1-x)\log^2 x-\frac{1}{24}\log^4 x+\frac{1}{2}\log^2(1-x)\log^2 x-\frac{1}{3}\log(1-x)\log^3 x$$

$$+\frac{1}{2}\left[\log(1-x)+\frac{1}{3}\log x\right]\log^2 x\log\left(\frac{x}{1-x}\right)-\frac{1}{4}\log^2 x\left[\log\left(\frac{x}{1-x}\right)\right]^2$$

$$+\frac{1}{8}\left[\log\left(\frac{x}{1-x}\right)\right]^4-\frac{1}{2}\log^2(1-x)Li_2(1-x)+\frac{1}{2}\log^2 x\,Li_2(x)$$

$$+\frac{1}{2}\left[\log\left(\frac{x}{1-x}\right)\right]^2 Li_2\left(\frac{-x}{1-x}\right)+\log(1-x)Li_3(1-x)-\log x\,Li_3(x)$$

$$-\log\left(\frac{x}{1-x}\right)Li_3\left(\frac{-x}{1-x}\right)-\left[Li_4(1-x)-Li_4(x)\right]+\varsigma(4)$$

$$+\left[\varsigma(3)-Li_3(1-x)\right]\log x-\frac{1}{2}\left[Li_2(1-x)\right]^2+\frac{1}{2}\varsigma^2(2)$$

$$+\frac{1}{6}\log^3(1-x)\log x+\frac{1}{2}\log^2(1-x)\,Li_2(1-x)-\log(1-x)\,Li_3(1-x)+Li_4(1-x)-\varsigma(4)-Li_4(x)=0$$

Making the obvious cancellations this becomes

(3.124a)

$$\log\left(\frac{x}{1-x}\right)Li_3\left(\frac{-x}{1-x}\right)=$$

$$-\frac{1}{12}\log^4 x+\frac{1}{4}\log^2(1-x)\log^2 x-\frac{1}{3}\log(1-x)\log^3 x+\frac{1}{6}\log^3(1-x)\log x$$

$$+\frac{1}{2}\left[\log(1-x)+\frac{1}{3}\log x\right]\log^2 x\log\left(\frac{x}{1-x}\right)-\frac{1}{4}\log^2 x\left[\log\left(\frac{x}{1-x}\right)\right]^2+\frac{1}{8}\left[\log\left(\frac{x}{1-x}\right)\right]^4$$

$$+\frac{1}{2}\left\{\left[Li_2(x)\right]^2-\left[Li_2(1-x)\right]^2\right\}+\frac{1}{2}\varsigma(2)\left[\varsigma(2)-2Li_2(x)-\log^2(1-x)\right]$$

$$+\frac{1}{2}\left[\log^2 x\,Li_2(x)-\log^2(1-x)Li_2(1-x)\right]+\frac{1}{2}\left[\log\left(\frac{x}{1-x}\right)\right]^2 Li_2\left(\frac{-x}{1-x}\right)$$

$$+\left[2Li_3(x)\log(1-x)-Li_3(1-x)\log x-Li_3(x)\log x\right]$$

$$+\varsigma(3)\left[\log x-\log(1-x)\right]$$

From (3.115) we have

$$Li_3\left(\frac{-x}{1-x}\right)=\varsigma(2)\log(1-x)-\frac{1}{2}\log x\log^2(1-x)+\frac{1}{6}\log^3(1-x)-Li_3(x)-Li_3(1-x)+\varsigma(3)$$

and hence we get for the first term in (3.124a)

(3.124b) $$\log\left(\frac{x}{1-x}\right)Li_3\left(\frac{-x}{1-x}\right)=$$

$$\left\{\log x-\log(1-x)\right\}\left[\varsigma(2)\log(1-x)-\frac{1}{2}\log x\log^2(1-x)+\frac{1}{6}\log^3(1-x)-Li_3(x)-Li_3(1-x)+\varsigma(3)\right]$$

Equating (3.124a) and (3.124b) gives us

$$\{\log x - \log(1-x)\}\left[\varsigma(2)\log(1-x) - \frac{1}{2}\log x \log^2(1-x) + \frac{1}{6}\log^3(1-x) - Li_3(x) - Li_3(1-x) + \varsigma(3)\right] =$$

$$-\frac{1}{12}\log^4 x + \frac{1}{4}\log^2(1-x)\log^2 x - \frac{1}{3}\log(1-x)\log^3 x + \frac{1}{6}\log^3(1-x)\log x$$

$$+\frac{1}{2}\left[\log(1-x) + \frac{1}{3}\log x\right]\log^2 x \log\left(\frac{x}{1-x}\right) - \frac{1}{4}\log^2 x\left[\log\left(\frac{x}{1-x}\right)\right]^2 + \frac{1}{8}\left[\log\left(\frac{x}{1-x}\right)\right]^4$$

$$+\frac{1}{2}\left\{[Li_2(x)]^2 - [Li_2(1-x)]^2\right\} + \frac{1}{2}\varsigma(2)\left\{\varsigma(2) - 2Li_2(x) - \log^2(1-x)\right\}$$

$$+\frac{1}{2}\left[\log^2 x\, Li_2(x) - \log^2(1-x)Li_2(1-x)\right] + \frac{1}{2}\left[\log\left(\frac{x}{1-x}\right)\right]^2 Li_2\left(\frac{-x}{1-x}\right)$$

$$+\left[2Li_3(x)\log(1-x) - Li_3(1-x)\log x - Li_3(x)\log x\right]$$

$$+\varsigma(3)\left[\log x - \log(1-x)\right]$$

This then becomes

(3.124c)

$$\left[Li_3(1-x) - Li_3(x)\right]\log(1-x) =$$

$$-\frac{1}{12}\log^4 x + \frac{3}{4}\log^2(1-x)\log^2 x - \frac{1}{3}\log(1-x)\log^3 x - \frac{1}{2}\log^3(1-x)\log x + \frac{1}{6}\log^4(1-x)$$

$$+\frac{1}{2}\left[\log(1-x) + \frac{1}{3}\log x\right]\log^2 x \log\left(\frac{x}{1-x}\right) - \frac{1}{4}\log^2 x\left[\log\left(\frac{x}{1-x}\right)\right]^2 + \frac{1}{8}\left[\log\left(\frac{x}{1-x}\right)\right]^4$$

$$+\frac{1}{2}\left\{[Li_2(x)]^2 - [Li_2(1-x)]^2\right\} + \frac{1}{2}\varsigma(2)\left[\varsigma(2) - 2Li_2(x) - \log^2(1-x)\right]$$

$$-\left[\log x - \log(1-x)\right]\varsigma(2)\log(1-x) + \frac{1}{2}\left[\log^2 x\, Li_2(x) - \log^2(1-x)Li_2(1-x)\right]$$

$$+\frac{1}{2}\left[\log\left(\frac{x}{1-x}\right)\right]^2 Li_2\left(\frac{-x}{1-x}\right)$$

Using (3.111) and (3.110f), the above identity may be expressed using only the polylogarithms $Li_3(1-x), Li_3(x)$ and $Li_2(x)$.

Let us now take stock of our armoury: we have the following identities involving the trilogarithm $Li_3(x)$ (see for example [126, p.106]). We now consider the irrational number $\rho = \frac{\sqrt{5}-1}{2}$ which is the reciprocal of the golden mean and is the positive root of the quadratic equation

$$\rho^2 + \rho - 1 = 0$$

We have from (3.111)

(3.125a) $$Li_3(\rho^2) = \frac{4}{5}\varsigma(3) + \frac{4}{5}\varsigma(2)\log\rho - \frac{2}{3}\log^3\rho$$

In addition we have the functional equation which is derived in [126, p.113]

(3.125b) $$Li_3(-\rho) - Li_3(-1/\rho) = -\frac{1}{6}\pi^2\log\rho - \frac{1}{6}\log^3\rho$$

From (3.124c) above we get

(3.125c) $$\left[Li_3(1-\rho) - Li_3(\rho)\right]\log(1-\rho) = k$$

where $k$ is equal to the right-hand side of (3.124c) evaluated at $x = \rho$.

Simple algebra shows that $\frac{\rho}{1-\rho} = \frac{1}{\rho}$ and hence we have

(3.125d) $$Li_3\left(-\frac{\rho}{1-\rho}\right) = Li_3(-1/\rho)$$

and

(3.125e) $$Li_3(\rho^2) = Li_3(1-\rho)$$

From (3.115) we get Landen's identity

(3.125f)

$$Li_3\left(\frac{-\rho}{1-\rho}\right) + Li_3(1-\rho) + Li_3(\rho) = \varsigma(3) + \varsigma(2)\log(1-\rho) - \frac{1}{2}\log\rho\log^2(1-\rho) + \frac{1}{6}\log^3(1-\rho)$$

Finally, from (4.4.67) we obtain

(3.125g) $$Li_3(\rho) + Li_3(-\rho) = \frac{1}{4}Li_3(\rho^2) = \frac{1}{5}\varsigma(3) + \frac{1}{5}\varsigma(2)\log\rho - \frac{1}{6}\log^3\rho$$

Regarding $\varsigma(3)$, $Li_3(\rho^2)$, $Li_3\left(-\frac{\rho}{1-\rho}\right)$, $Li_3(-\rho)$, $Li_3(-1/\rho)$, $Li_3(1-\rho)$ and $Li_3(\rho)$ as unknowns, I initially thought that I had "The Magnificent Seven", i.e. 7 simultaneous equations in 7 unknowns. Unfortunately, this system of simultaneous equations is not linearly independent and therefore not solvable with the current input. If we put $x = 1-\rho$ in (3.124c) do we get another linearly independent equation? This area may merit further consideration.

We also have the Euler/Landen/Spence identities for the dilogarithm: equations (3.126c) and (3.126e) are required when evaluating (3.124c) at $x = \rho$.

(3.126a) $$Li_2(\rho) + Li_2\left(\frac{-\rho}{1-\rho}\right) = -\frac{1}{2}\log^2(1-\rho)$$

(3.126b) $$Li_2(\rho) = \frac{\pi^2}{10} - \frac{1}{2}\log^2\rho$$

Hence we get

(3.126c) $$Li_2\left(\frac{-\rho}{1-\rho}\right) = -\frac{\pi^2}{10} + \frac{1}{2}\log^2\rho - \frac{1}{2}\log^2(1-\rho)$$

We also have from (3.110f)

(3.126d) $$Li_2(\rho) + Li_2(1-\rho) = \varsigma(2) - \log\rho\log(1-\rho)$$

and therefore

(3.126e) $$Li_2(1-\rho) = \frac{\pi^2}{15} + \frac{1}{2}\log^2\rho - \log\rho\log(1-\rho)$$

In (3.115c) we saw that

$$\sum_{n=1}^{\infty}\frac{u^n}{n^2}\sum_{k=1}^{n}\binom{n}{k}\frac{(-1)^k}{k^2} = -\frac{1}{2}\sum_{n=1}^{\infty}\frac{u^n}{n^2}\left[\left(H_n^{(1)}\right)^2 + H_n^{(2)}\right] = \int_0^u Li_3\left(\frac{-t}{1-t}\right)\frac{dt}{t}$$

and using Landen's identity (3.115) we have

$$\int_0^u Li_3\left(\frac{-t}{1-t}\right)\frac{dt}{t} =$$

$$\varsigma(2)\int_0^u \frac{\log(1-t)}{t}dt - \frac{1}{2}\int_0^u \frac{\log t\log^2(1-t)}{t}dt + \int_0^u \frac{\varsigma(3) - Li_3(1-t)}{t}dt + \frac{1}{6}\int_0^u \frac{\log^3(1-t)}{t}dt - \int_0^u \frac{Li_3(t)}{t}dt$$

Using (3.121b) and (3.122) this becomes

$$= -\varsigma(2) Li_2(u) - \frac{1}{2}\int_0^u \frac{\log t \log^2(1-t)}{t} dt + \left[\varsigma(3) - Li_3(1-u)\right]\log u - \frac{1}{2}\left[Li_2(1-u)\right]^2 + \frac{1}{2}\varsigma^2(2)$$

$$+\frac{1}{6}\log^3(1-u)\log u + \frac{1}{2}\log^2(1-u) Li_2(1-u) - \log(1-u) Li_3(1-u) + Li_4(1-u) - \varsigma(4) - Li_4(u)$$

Referring back to (3.115c) we then see that (instead of (3.120)) we may write

$$\int_0^u \frac{\log t \log^2(1-t)}{t} dt$$

$$= -2\varsigma(2) Li_2(u) + 2\left[\varsigma(3) - Li_3(1-u)\right]\log u - \left[Li_2(1-u)\right]^2 + \varsigma^2(2)$$

$$+\frac{1}{3}\log^3(1-u)\log u + \log^2(1-u) Li_2(1-u) - 2\log(1-u) Li_3(1-u) + 2Li_4(1-u) - 2\varsigma(4) - 2Li_4(u)$$

$$+\sum_{n=1}^{\infty} \frac{u^n}{n^2}\left[\left(H_n^{(1)}\right)^2 + H_n^{(2)}\right]$$

Equating this with (3.120) we get

$$-2\varsigma(2) Li_2(u) + 2\left[\varsigma(3) - Li_3(1-u)\right]\log u - \left[Li_2(1-u)\right]^2 + \varsigma^2(2) + \sum_{n=1}^{\infty} \frac{u^n}{n^2}\left[\left(H_n^{(1)}\right)^2 + H_n^{(2)}\right]$$

$$= \frac{1}{2}\log^2(1-u)\log^2 u + \frac{1}{12}\log^4 u - \log^2(1-u)\log^2 u + \frac{1}{3}\log(1-u)\log^3 u$$

$$-\left[\log(1-u) + \frac{1}{3}\log u\right]\log^2 u \log\left(\frac{u}{1-u}\right) + \frac{1}{2}\log^2 u \left[\log\left(\frac{u}{1-u}\right)\right]^2$$

$$-\frac{1}{4}\left[\log\left(\frac{u}{1-u}\right)\right]^4 - \log^2 u\, Li_2(u) - \left[\log\left(\frac{u}{1-u}\right)\right]^2 Li_2\left(\frac{-u}{1-u}\right) + 2\log u\, Li_3(u)$$

$$+2\log\left(\frac{u}{1-u}\right) Li_3\left(\frac{-u}{1-u}\right) - 2Li_4\left(\frac{-u}{1-u}\right)$$

Letting $u = 1/2$ we obtain

$$-2\varsigma(2) Li_2(1/2) - 2\left[\varsigma(3) - Li_3(1/2)\right]\log 2 - \left[Li_2(1/2)\right]^2 + \varsigma^2(2) + \sum_{n=1}^{\infty} \frac{1}{2^n n^2}\left[\left(H_n^{(1)}\right)^2 + H_n^{(2)}\right]$$

$$=\frac{1}{2}\log^4 2+\frac{1}{12}\log^4 2-\log^4 2+\frac{1}{3}\log^4 2-\log^2 2\,Li_2(1/2)-2\log 2\,Li_3(1/2)-2Li_4\left(-1\right)$$

$$=-\frac{1}{12}\log^4 2-\log^2 2\,Li_2(1/2)-2\log 2\,Li_3(1/2)-2Li_4\left(-1\right)$$

This may be written as

$$\sum_{n=1}^{\infty}\frac{1}{2^n n^2}\left[\left(H_n^{(1)}\right)^2+H_n^{(2)}\right]$$

$$=-\frac{1}{12}\log^4 2+[2\varsigma(2)-\log^2 2]Li_2(1/2)-\varsigma^2(2)+\left[Li_2(1/2)\right]^2$$

$$+2\left[\varsigma(3)-Li_3(1/2)\right]\log 2-2Li_4\left(-1\right)$$

By differentiation we can easily verify that

$$\int\frac{\log(1+y^2)}{y}dy=-\frac{1}{2}Li_2(-y^2)$$

and using integration by parts we have

$$\int_0^t\frac{\log^2(1+y^2)}{y}dy=-\frac{1}{2}\log(1+y^2)Li_2(-y^2)\Big|_0^t+\int_0^t\frac{Li_2(-y^2)y}{1+y^2}dy$$

$$=-\frac{1}{2}\log(1+t^2)Li_2(-t^2)+\int_0^t\frac{Li_2(-y^2)y}{1+y^2}dy$$

By an obvious substitution the last integral becomes

$$=\frac{1}{2}\int_0^{t^2}\frac{Li_2(-u)}{1+u}du$$

We now recall a trick which will be used in (4.4.146) and modify it slightly: we have

$$\frac{1-(-u)^n}{1+u}-\sum_{k=1}^{n}(-u)^{k-1}=0$$

Hence the integral of the above is identically equal to zero.

$$\int \left[ \frac{1-(-u)^n}{1+u} - \sum_{k=1}^{n} (-u)^{k-1} \right] du = 0$$

Therefore we have the simple result (certainly simpler than the Wolfram Integrator which displays a hypergeometric function as its output!)

$$\log(1+u) - \sum_{k=1}^{n} \frac{(-1)^{k+1} u^k}{k} = \int \frac{(-u)^n}{1+u} du$$

Hence we get

$$\log 2 - \sum_{k=1}^{n} \frac{(-1)^{k+1}}{k} = \int_0^1 \frac{(-u)^n}{1+u} du$$

We have for $q \geq 2$

$$\int_0^t \frac{Li_q(-u)}{1+u} du = \sum_{n=1}^{\infty} \frac{1}{n^q} \int_0^t \frac{(-u)^n}{1+u} du$$

$$= \sum_{n=1}^{\infty} \frac{1}{n^q} \left[ \log(1+t) - \sum_{k=1}^{n} \frac{(-1)^{k+1} t^k}{k} \right]$$

and hence

(3.127) $$\int_0^t \frac{Li_q(-u)}{1+u} du = \varsigma(q) \log(1+t) - \sum_{n=1}^{\infty} \frac{1}{n^q} \sum_{k=1}^{n} \frac{(-1)^{k+1} t^k}{k}$$

We therefore have for $t = 1$

(3.127a) $$\int_0^1 \frac{Li_2(-u)}{1+u} du = \varsigma(2) \log 2 - \sum_{n=1}^{\infty} \frac{1}{n^2} \sum_{k=1}^{n} \frac{(-1)^{k+1}}{k}$$

As an aside, we note that the Wolfram Integrator gives the following output

(3.127ai)
$$\int \frac{Li_2(-u)}{1+u} du = \log(1+u) Li_2(-u) + \log(-u) \log^2(1+u) + 2\log(1+u) Li_2(1+u) - 2Li_3(1+u)$$

Using integration by parts we have

$$\int_0^t \frac{Li_2(-u)}{1+u} du = Li_2(-t) \log(1+t) + \int_0^t \frac{\log^2(1+u)}{u} du$$

As shown by Lewin [100, p.310] we have

(3.127aii)

$$\int_0^t \frac{\log^2(1+u)}{u} du = \log t \log^2(1+t) - \frac{2}{3}\log^3(1+t) - 2\log(1+t)\, Li_2\left(\frac{1}{1+t}\right) - 2Li_3\left(\frac{1}{1+t}\right) + 2\varsigma(3)$$

Hence we get

(3.127b) $$\int_0^t \frac{Li_2(-u)}{1+u} du = Li_2(-t)\log(1+t) + \log t \log^2(1+t) - \frac{2}{3}\log^3(1+t)$$

$$-2\log(1+t)\, Li_2\left(\frac{1}{1+t}\right) - 2Li_3\left(\frac{1}{1+t}\right) + 2\varsigma(3)$$

and we can therefore derive the following integral which is recorded by Devoto and Duke [53a] in their table of integrals for Feynman diagram calculations

$$\int_0^1 \frac{Li_2(-u)}{1+u} du = \frac{1}{4}\varsigma(3) - \frac{1}{2}\varsigma(2)\log 2$$

We therefore deduce that

(3.128) $$\sum_{n=1}^{\infty} \frac{1}{n^2} \sum_{k=1}^{n} \frac{(-1)^{k+1}}{k} = \frac{3}{2}\varsigma(2)\log 2 - \frac{1}{4}\varsigma(3)$$

We have using integration by parts

$$I_3 = \int_0^t \frac{Li_3(-u)}{1+u} du = Li_3(-t)\log(1+t) - \int_0^t \log(1+u)\frac{Li_2(-u)}{u} du$$

$$\int_0^t \log(1+u)\frac{Li_2(-u)}{u} du = -\left[Li_2(-t)\right]^2 + \int_0^t \log(1+u)\frac{Li_2(-u)}{u} du$$

$$\int_0^t \log(1+u)\frac{Li_2(-u)}{u} du = -\frac{1}{2}\left[Li_2(-t)\right]^2$$

and thus we get

(3.128a) $$\int_0^t \frac{Li_3(-u)}{1+u} du = Li_3(-t)\log(1+t) + \frac{1}{2}\left[Li_2(-t)\right]^2$$

Therefore, using (3.127) we obtain

(3.128b) $$\sum_{n=1}^{\infty}\frac{1}{n^3}\sum_{k=1}^{n}\frac{(-1)^{k+1}t^k}{k}=\varsigma(3)\log(1+t)-Li_3(-t)\log(1+t)-\frac{1}{2}\left[Li_2(-t)\right]^2$$

and in particular we have

(3.128c) $$\sum_{n=1}^{\infty}\frac{1}{n^3}\sum_{k=1}^{n}\frac{(-1)^{k+1}}{k}=\varsigma(3)\log 2-Li_3(-1)\log 2-\frac{1}{2}\left[Li_2(-1)\right]^2$$

$$=\frac{7}{4}\varsigma(3)\log 2-\frac{1}{8}\varsigma^2(2)$$

Similarly we have

$$I_4=\int_0^t\frac{Li_4(-u)}{1+u}du=Li_4(-t)\log(1+t)-\int_0^t\log(1+u)\frac{Li_3(-u)}{u}du$$

$$\int_0^t\log(1+u)\frac{Li_3(-u)}{u}du=-Li_2(-t)Li_3(-t)+\int_0^t Li_2(-u)\frac{Li_2(-u)}{u}du$$

$$\int_0^t Li_2(-u)\frac{Li_2(-u)}{u}du=Li_3(-u)Li_2(-u)+\int_0^t Li_3(-u)\frac{\log(1+u)}{u}du$$

$$\int_0^t Li_3(-u)\frac{\log(1+u)}{u}du=Li_4(-t)\log(1+t)-\int_0^t\frac{Li_4(-u)}{1+u}du$$

However, in this case everything cancels out and we simply get $0=0$ !

Upon integrating (3.128b) we see that

$$\sum_{n=1}^{\infty}\frac{1}{n^3}\sum_{k=1}^{n}\frac{(-1)^{k+1}t^k}{k^2}=\varsigma(3)Li_2(-t)-\int_0^t\left(Li_3(-t)\log(1+t)+\frac{1}{2}\left[Li_2(-t)\right]^2\right)\frac{dt}{t}$$

Equations (3.128) and (3.128c) are particular cases of a more general formula derived by Sitaramachandrarao [120a] in 1987 (see also the paper by Flajolet and Salvy [69]).

(3.129)

$$2\sum_{n=1}^{\infty}\frac{1}{n^q}\sum_{k=1}^{n}\frac{(-1)^{k+1}}{k}=2\varsigma(q)\log 2-q\varsigma(q+1)+2\varsigma_a(q+1)+\sum_{k=1}^{q}\varsigma_a(k)\varsigma_a(q-k+1)$$

where $\varsigma_a(1)=\log 2$.

We also see that

$$\frac{1}{x}\log(1+ux)-\sum_{k=1}^{n}\frac{(-1)^{k+1}(ux)^k}{k}=\int\frac{(-ux)^n}{1+ux}du$$

and therefore we have

$$\int_0^t\frac{Li_q(-ux)}{1+ux}du=\frac{1}{x}\varsigma(q)\log(1+tx)-\sum_{n=1}^{\infty}(-1)^n\frac{x^n}{n^q}\sum_{k=1}^{n}\frac{(-1)^{k+1}t^k}{k}$$

Using integration by parts we have

$$\int_0^t\frac{Li_2(-ux)}{1+ux}du=\frac{1}{x}Li_2(-tx)\log(1+tx)+\frac{1}{x^2}\int_0^t\frac{\log^2(1+ux)}{u}du$$

Referring to formula (A6f) in Devoto and Duke [53a] we see that

$$\int_0^1\frac{\log^2(1+ux)}{u}du=2\sum_{n=2}^{\infty}(-1)^n\frac{x^n}{n^2}\sum_{k=1}^{n-1}\frac{1}{k}$$

and hence we get

$$\frac{1}{x}\varsigma(2)\log(1+x)-\sum_{n=1}^{\infty}(-1)^n\frac{x^n}{n^2}\sum_{k=1}^{n}\frac{(-1)^{k+1}}{k}=\frac{1}{x}Li_2(-x)\log(1+x)+\frac{2}{x^2}\sum_{n=2}^{\infty}(-1)^n\frac{x^n}{n^2}\sum_{k=1}^{n-1}\frac{1}{k}$$

This may be written as

$$x\log(1+x)\left[\varsigma(2)-Li_2(-x)\right]-x^2\sum_{n=1}^{\infty}(-1)^n\frac{x^n}{n^2}\sum_{k=1}^{n}\frac{(-1)^{k+1}}{k}=2\sum_{n=1}^{\infty}(-1)^n\frac{x^n}{n^2}\sum_{k=1}^{n}\frac{1}{k}-2\sum_{n=1}^{\infty}(-1)^n\frac{x^n}{n^3}$$

and hence we obtain

$$x\log(1+x)\left[\varsigma(2)-Li_2(-x)\right]+2Li_3(-x)-x^2\sum_{n=1}^{\infty}(-1)^n\frac{x^n}{n^2}\sum_{k=1}^{n}\frac{(-1)^{k+1}}{k}=2\sum_{n=1}^{\infty}(-1)^n\frac{H_n^{(1)}}{n^2}x^n$$

Letting $x\to -x$ we have

$$-x\log(1-x)\left[\varsigma(2)-Li_2(x)\right]+2Li_3(x)-x^2\sum_{n=1}^{\infty}\frac{x^n}{n^2}\sum_{k=1}^{n}\frac{(-1)^{k+1}}{k}=2\sum_{n=1}^{\infty}\frac{H_n^{(1)}}{n^2}x^n$$

Alternatively we may employ Landen's functional equation (3.111a): dividing that by $1/(1+u)$ and $u$ respectively we obtain

$$\int_0^t \frac{1}{1+u} Li_2\left(\frac{u}{1+u}\right) du = -\int_0^t \frac{Li_2(-u)}{1+u} du - \frac{1}{2}\int_0^t \frac{\log^2(1+u)}{1+u} du$$

and

$$\int_0^t \frac{1}{u} Li_2\left(\frac{u}{1+u}\right) du = -\int_0^t \frac{Li_2(-u)}{u} du - \frac{1}{2}\int_0^t \frac{\log^2(1+u)}{u} du$$

It is obvious using partial fractions that

(3.130) $$\int_0^t \frac{1}{u(1+u)} Li_2\left(\frac{u}{1+u}\right) du = -\int_0^t \frac{1}{1+u} Li_2\left(\frac{u}{1+u}\right) du + \int_0^t \frac{1}{u} Li_2\left(\frac{u}{1+u}\right) du$$

$$= \int_0^t \frac{Li_2(-u)}{1+u} du + \frac{1}{2}\int_0^t \frac{\log^2(1+u)}{1+u} du - \int_0^t \frac{Li_2(-u)}{u} du - \frac{1}{2}\int_0^t \frac{\log^2(1+u)}{u} du$$

As before, employing the series definition of the trilogarithm it is easily seen that

$$Li'_{n+1}\left(\frac{u}{1+u}\right) = \frac{1}{u(1+u)} Li_n\left(\frac{u}{1+u}\right)$$

and, in particular, we have for $n=1$

$$Li'_2\left(\frac{u}{1+u}\right) = \frac{1}{u(1+u)} \sum_{k=1}^{\infty} \frac{1}{k}\left(\frac{u}{1+u}\right)^k = \frac{\log(1+u)}{u(1+u)}$$

Hence, it is clear that

(3.131) $$\int_0^t \frac{1}{u(1+u)} Li_2\left(\frac{u}{1+u}\right) du = Li_3\left(\frac{t}{1+t}\right)$$

Also, it is easily seen that

$$\int_0^t \frac{Li_2(-u)}{u} du = Li_3(-t)$$

and therefore, from (3.130) and (3.127a) we obtain

(3.132)

$$Li_3\left(\frac{t}{1+t}\right) = \varsigma(2)\log(1+t) - \sum_{n=1}^{\infty} \frac{1}{n^2} \sum_{k=1}^{n} \frac{(-1)^{k+1} t^k}{k} + \frac{1}{6}\log^3(1+t) - Li_3(-t) - \frac{1}{2}\int_0^t \frac{\log^2(1+u)}{u} du$$

This is to be compared with the puzzling output which we get from the Wolfram Integrator for the following integral

(3.133)

$$\int_0^t \frac{\log^2(1+u)}{u}du = \log(-t)\log^2(1+t) + 2\log(1+t)Li_2(1+t) - 2Li_3(1+t) + 2\varsigma(3)$$

where we have imaginary numbers arising from (i) the logarithm of a negative number and (ii) the argument of the polylogarithm greater than one (and hence initially I thought this was not convergent). Compare the above representation with (3.42).

With $t = 1$ in (3.132) we have

$$\frac{1}{2}\int_0^1 \frac{\log^2(1+u)}{u}du = \varsigma(2)\log 2 - \sum_{n=1}^{\infty}\frac{1}{n^2}\sum_{k=1}^{n}\frac{(-1)^{k+1}}{k} + \frac{1}{6}\log^3 2 - Li_3(-1) - Li_3(1/2)$$

Devoto and Duke [53a] record the following integral

(3.134) $$\int_0^1 \frac{\log^2(1+u)}{u}du = \frac{1}{4}\varsigma(3)$$

and we accordingly obtain another proof of (3.128)

(3.135) $$\sum_{n=1}^{\infty}\frac{1}{n^2}\sum_{k=1}^{n}\frac{(-1)^{k+1}}{k} = \varsigma(2)\log 2 - \frac{1}{8}\varsigma(3) + \frac{1}{6}\log^3 2 - Li_3(-1) - Li_3(1/2)$$

$$= \frac{3}{2}\varsigma(2)\log 2 - \frac{1}{4}\varsigma(3)$$

By differentiation we can readily verify that

$$\int \frac{\log(1+x)}{x}dx = -Li_2(-x)$$

and using integration by parts we have

$$\int_0^t \frac{\log^2(1+x)}{x}dx = -\log(1+x)Li_2(-x)\Big|_0^t + \int_0^t \frac{Li_2(-x)}{1+x}dx$$

$$= -\log(1+t)Li_2(-t) + \int_0^t \frac{Li_2(-x)}{1+x}dx$$

$$= \log(1+t)\left[\varsigma(2) - Li_2(-t)\right] - \sum_{n=1}^{\infty} \frac{1}{n^2} \sum_{k=1}^{n} \frac{(-1)^{k+1} t^k}{k}$$

Hence we have as before

$$\int_0^1 \frac{\log^2(1+x)}{x} dx = \frac{3}{2}\varsigma(2)\log 2 - \sum_{n=1}^{\infty} \frac{1}{n^2} \sum_{k=1}^{n} \frac{(-1)^{k+1}}{k}$$

Christmas came early this year because on 2 December 2005 Jonathan Sondow kindly forwarded me a copy of Lewin's book on polylogarithms [100] and accordingly I have been able to reproduce Lewin's evaluation of the above integral in (4.4.100g) of Volume III.

We have from (3.127aii)

(3.136)

$$\int_0^t \frac{\log^2(1+x)}{x} dx = \log t \log^2(1+t) - \frac{2}{3}\log^3(1+t) - 2\log(1+t) Li_2\left(\frac{1}{1+t}\right) - 2Li_3\left(\frac{1}{1+t}\right) + 2\varsigma(3)$$

and therefore we obtain from (3.132)

(3.137) $$\sum_{n=1}^{\infty} \frac{1}{n^2} \sum_{k=1}^{n} \frac{(-1)^{k+1} t^k}{k} = \varsigma(2)\log(1+t) - Li_3\left(\frac{t}{1+t}\right) + \frac{1}{6}\log^3(1+t) - Li_3(-t)$$

$$-\frac{1}{2}\log t \log^2(1+t) + \frac{1}{3}\log^3(1+t) + \log(1+t) Li_2\left(\frac{1}{1+t}\right) + Li_3\left(\frac{1}{1+t}\right) - \varsigma(3)$$

With $t = 1$ we obtain (3.128).

From (3.30) we have

$$\frac{1}{2}\log^2(1+x) = \sum_{n=1}^{\infty} (-1)^{n+1} \frac{H_n^{(1)}}{n+1} x^{n+1}$$

and integration gives us

$$\frac{1}{2}\int_0^t \frac{\log^2(1+x)}{x} dx = \sum_{n=1}^{\infty} (-1)^{n+1} \frac{H_n^{(1)}}{(n+1)^2} t^{n+1}$$

Now compare with (3.26)

$$\frac{1}{2}\int_0^t \frac{\log^2(1+x)}{x} dx = \sum_{n=1}^{\infty} (-1)^n \frac{H_n^{(1)}}{n^2} t^n - \int_0^t \frac{Li_2(-x)}{x} dx$$

$$= \sum_{n=1}^{\infty} (-1)^n \frac{H_n^{(1)}}{n^2} t^n - Li_3(-t)$$

With $t = 1$ we obtain

$$\sum_{n=1}^{\infty} (-1)^n \frac{H_n^{(1)}}{n^2} = \frac{1}{2} \int_0^1 \frac{\log^2(1+x)}{x} dx + Li_3(-1)$$

Using (4.4.67) we have

$$Li_3(1) + Li_3(-1) = \frac{1}{4} Li_3(1)$$

and we therefore get (see also (4.4.100e) and (4.4.167k))

(3.138) $$\sum_{n=1}^{\infty} (-1)^n \frac{H_n^{(1)}}{n^2} = -\frac{5}{8} \varsigma(3)$$

Similarly we have

$$\frac{1}{2} \log^2(1-x) = \sum_{n=1}^{\infty} \frac{H_n^{(1)}}{n+1} x^{n+1}$$

Integration results in

$$\frac{1}{2} \int_0^t \frac{\log^2(1-x)}{x} dx = \sum_{n=1}^{\infty} \frac{H_n^{(1)}}{(n+1)^2} t^{n+1}$$

Now compare with (3.26)

$$\frac{1}{2} \int_0^t \frac{\log^2(1-x)}{x} dx = \sum_{n=1}^{\infty} \frac{H_n^{(1)}}{n^2} t^n - Li_3(t)$$

and we obtain

(3.139) $$\sum_{n=1}^{\infty} \frac{H_n^{(1)}}{(n+1)^2} t^{n+1} = \sum_{n=1}^{\infty} \frac{H_n^{(1)}}{n^2} t^n - Li_3(t)$$

Therefore with $t = 1$ we have the well-known result

(3.140) $$\sum_{n=1}^{\infty} \frac{H_n^{(1)}}{(n+1)^2} = \varsigma(3)$$

Using integration by parts it is easily shown that

$$\frac{1}{2}\int_0^t \frac{\log^2(1-x)}{x}dx = \frac{1}{2}\log^2(1-t)\log t + \log(1-t)Li_2(1-t) - Li_3(1-t) + \varsigma(3)$$

and we accordingly get as before

$$\sum_{n=1}^{\infty} \frac{H_n^{(1)}}{n^2} t^n = \frac{1}{2}\log^2(1-t)\log t + \log(1-t)Li_2(1-t) - Li_3(1-t) + Li_3(t) + \varsigma(3)$$

Reference to (3.127ai) and the use of the Wolfram Integrator gives us

$$\int \frac{Li_2(-u)}{1+u}du = \log(1+u)Li_2(-u) + \log(-u)\log^2(1+u) + 2\log(1+u)Li_2(1+u) - 2Li_3(1+u)$$

and hence we have

$$\int_0^{-t} \frac{Li_2(-u)}{1+u}du = \log(1+x)\left[Li_2(1+u) + \varsigma(2)\right] - 2Li_3(1+u)\Big|_0^{-t}$$

$$= \log(1-t)\left[Li_2(1-t) + \varsigma(2)\right] - 2Li_3(1-t) + 2\varsigma(3)$$

Therefore we have using (3.127)

(3.140a) $$\sum_{n=1}^{\infty} \frac{1}{n^2}\sum_{k=1}^{n} \frac{t^k}{k} = 2\varsigma(3) - 2Li_3(1-t) + Li_2(1-t)\log(1-t)$$

In (4.4.64a) in Volume III we shall see that

$$\sum_{n=1}^{\infty} \frac{1}{n^2}\sum_{k=1}^{n} \frac{(1-x)^k}{k} = 2\varsigma(3) - 2Li_3(x) + \log x\, Li_2(x)$$

and with $t = 1-x$ we recover (3.140a).

Letting $t \to -t$ in (3.140a) we have

$$\sum_{n=1}^{\infty} \frac{1}{n^2}\sum_{k=1}^{n} \frac{(-1)^{k+1}}{k}t^k = -2\varsigma(3) + 2Li_3(1+t) - \log(1+t)\,Li_2(1+t)$$

and comparing this with (3.137) **suggests** that

(3.141)

$$2Li_3(1+t) - \log(1+t)\,Li_2(1+t) = \varsigma(2)\log(1+t) - Li_3\left(\frac{t}{1+t}\right) + \frac{1}{6}\log^3(1+t) - Li_3(-t)$$

$$-\frac{1}{2}\log t\log^2(1+t)+\frac{1}{3}\log^3(1+t)+\log(1+t)\,Li_2\left(\frac{1}{1+t}\right)+Li_3\left(\frac{1}{1+t}\right)+\varsigma(3)$$

Note: Equation references from (3.142) to (3.199) have been deliberately omitted.

## AN APPLICATION OF THE BINOMIAL THEOREM

Using the binomial theorem we have for $|z|<1$

(3.200) $$\frac{1}{(1-z)^x}=\sum_{n=0}^{\infty}\frac{(x)_n}{n!}z^n$$

where $(x)_n$ is the ascending factorial symbol (also known as the Pochhamer symbol) defined by [25, p.16] as

(3.201) $$(x)_n = x(x+1)(x+2)...(x+n-1) \text{ if } n>0 \text{ and}$$

$$(x)_0 = 1$$

Shen [120] employed equation (3.200) in his 1995 paper "Remarks on some integrals and series involving the Stirling numbers and $\varsigma(n)$", but here I have adopted a different and more elementary approach to obtain some additional identities.

Differentiating (3.200) with respect to $x$ we get

(3.202) $$-\frac{\log(1-z)}{(1-z)^x}=\sum_{n=1}^{\infty}\frac{f_n'(x)}{n!}z^n$$

where we have designated $f_n(x)=(x)_n$ as a more convenient notation. We have

$$\log f_n(x)=\log x+\log(x+1)...+\log(x+n-1)$$

and logarithmic differentiation results in

$$\frac{f_n'(x)}{f_n(x)}=\frac{1}{x}+\frac{1}{x+1}...+\frac{1}{x+n-1}=H_n^{(1)}(x)$$

where $H_n^{(1)}(1)=H_n^{(1)}$.

Therefore we have

$$f_n'(x)=x(x+1)(x+2)...(x+n-1)\left[\frac{1}{x}+\frac{1}{x+1}...+\frac{1}{x+n-1}\right]$$

With $x=1$ we obtain $f_n'(1)=n!H_n^{(1)}$ and hence we have another derivation of (3.28)

(3.203) $$-\frac{\log(1-z)}{(1-z)}=\sum_{n=1}^{\infty}H_n^{(1)}\,z^n$$

We easily see that

$$\log(1-z)=-\sum_{n=1}^{\infty}H_n^{(1)}(1-z)\,z^n$$

$$=-\sum_{n=1}^{\infty}(H_n^{(1)}-H_{n-1}^{(1)})\,z^n$$

$$=-\sum_{n=1}^{\infty}\frac{z^n}{n}$$

We have

$$-\int_0^t\frac{\log(1-z)}{(1-z)}dz=\frac{1}{2}\log^2(1-t)=\sum_{n=1}^{\infty}\frac{H_n^{(1)}}{n+1}t^{n+1}=\sum_{n=1}^{\infty}\left[H_{n+1}^{(1)}-\frac{1}{n+1}\right]\frac{t^{n+1}}{n+1}$$

$$=\sum_{n=1}^{\infty}\frac{H_n^{(1)}}{n}t^n-Li_2(t)$$

We may also easily obtain

(3.204) $$-\int_0^t\frac{\log(1-z)}{z(1-z)}dz=\sum_{n=1}^{\infty}\frac{H_n^{(1)}}{n}t^n$$

and integration by parts shows that

(3.205) $$-\int_0^t\frac{\log(1-z)}{z(1-z)}dz=\frac{1}{2}\log^2(1-t)+Li_2(t)$$

We therefore obtain a further proof of (3.31)

(3.206) $$\sum_{n=1}^{\infty}\frac{H_n^{(1)}}{n}t^n=\frac{1}{2}\log^2(1-t)+Li_2(t)$$

Multiplying (3.206) across by $\log t/t$ and integrating with the use of the elementary integral

(3.207) $$\int_0^x t^{n-1}\log t\,dt=\left[-\frac{1}{n^2}+\frac{\log x}{n}\right]x^n$$

we obtain

$$\sum_{n=1}^{\infty}\frac{H_n^{(1)}}{n}\left[-\frac{1}{n^2}+\frac{\log x}{n}\right]x^n=\frac{1}{2}\int_0^x\frac{\log^2(1-t)\log t}{t}dt+\int_0^x\frac{Li_2(t)\log t}{t}dt$$

Integration by parts results in

$$\int_0^x\frac{Li_2(t)\log t}{t}dt=\frac{1}{2}\log^2 x\,Li_2(x)+\frac{1}{2}\int_0^x\frac{\log^2 t\log(1-t)}{t}dt$$

$$\int_0^x\frac{\log^2 t\log(1-t)}{t}dt=\frac{1}{3}\log^3 x\log(1-x)+\frac{1}{3}\int_0^x\frac{\log^3 t}{(1-t)}dt$$

(3.208) $$\int_0^x\frac{\log^3 t}{(1-t)}dt=-\log^3 x\log(1-x)-3\log^2 x\,Li_2(x)+6\log x\,Li_3(x)-6Li_4(x)$$

(3.208a) $$\int_0^x\frac{\log^2 t\log(1-t)}{t}dt=-\log^2 x\,Li_2(x)+2\log x\,Li_3(x)-2Li_4(x)$$

Therefore we have

(3.209) $$\int_0^x\frac{Li_2(t)\log t}{t}dt=\log x\,Li_3(x)-Li_4(x)$$

It should be noted that this integral may be evaluated much more directly as follows

$$\int_0^x\frac{Li_2(t)}{t}\log t\,dt=\log x\,Li_3(x)-\int_0^x\frac{Li_3(t)}{t}dt$$

Hence we obtain

(3.210) $$\int_0^x\frac{\log^2(1-t)\log t}{t}dt=2Li_4(x)-2\log x\,Li_3(x)-2\sum_{n=1}^{\infty}\frac{H_n^{(1)}}{n^3}x^n+2\log x\sum_{n=1}^{\infty}\frac{H_n^{(1)}}{n^2}x^n$$

A few pages later we shall show in (3.226) that

$$\int_0^x\frac{\log^2(1-t)\log t}{t}dt=\log x\,Li_3(x)-Li_4(x)+Li_2(x)Li_2(1-x)-\varsigma(2)Li_2(x)$$

$$+\frac{1}{2}\left[Li_2(x)\right]^2+\sum_{n=1}^{\infty}\frac{H_n^{(2)}}{n^2}x^n-\log x\sum_{n=1}^{\infty}\frac{H_n^{(2)}}{n}x^n$$

and we therefore get

(3.211) $$3Li_4(x)-3\log x\,Li_3(x)-Li_2(x)Li_2(1-x)+\varsigma(2)Li_2(x)-\frac{1}{2}\left[Li_2(x)\right]^2=$$

$$2\sum_{n=1}^{\infty}\frac{H_n^{(1)}}{n^3}x^n-2\log x\sum_{n=1}^{\infty}\frac{H_n^{(1)}}{n^2}x^n+\sum_{n=1}^{\infty}\frac{H_n^{(2)}}{n^2}x^n-\log x\sum_{n=1}^{\infty}\frac{H_n^{(2)}}{n}x^n$$

A quick check shows that the arithmetic works with $x=1$ whereby we obtain

$$2\sum_{n=1}^{\infty}\frac{H_n^{(1)}}{n^3}+\sum_{n=1}^{\infty}\frac{H_n^{(2)}}{n^2}=3\varsigma(4)+\frac{1}{2}\varsigma^2(2)$$

in agreement with (3.211b).

We have from (4.4.43zh) in Volume II(b) for $0\le x<1$

(3.211i)

$$\sum_{n=1}^{\infty}\frac{H_n^{(2)}}{n}x^n=Li_3(x)-\log(1-x)Li_2(x)-\log x\log^2(1-x)-2\log(1-x)\,Li_2(1-x)+2Li_3(1-x)-2\varsigma(3)$$

and from (3.105d) for $0\le x\le 1$

(3.211ii)

$$\sum_{n=1}^{\infty}\frac{H_n}{n^2}x^n=\frac{1}{2}\log^2(1-x)\log x+\log(1-x)Li_2(1-x)+Li_3(x)-Li_3(1-x)+\varsigma(3)$$

Substituting these identities in (3.211) we get

(3.211a) $$3Li_4(x)-3\log x\,Li_3(x)-Li_2(x)Li_2(1-x)+\varsigma(2)Li_2(x)-\frac{1}{2}\left[Li_2(x)\right]^2=$$

$$2\sum_{n=1}^{\infty}\frac{H_n^{(1)}}{n^3}x^n+\sum_{n=1}^{\infty}\frac{H_n^{(2)}}{n^2}x^n+\log x\log(1-x)Li_2(x)-3\log x\,Li_3(x)-4\log x\,Li_3(1-x)$$

Using (3.108b)

$$2\sum_{n=1}^{\infty}\frac{H_n^{(1)}}{n^3}x^n+\sum_{n=1}^{\infty}\frac{H_n^{(2)}}{n^2}x^n-\sum_{n=1}^{\infty}\frac{\left(H_n^{(1)}\right)^2}{n^2}x^n=$$

$$\frac{1}{3}\log^3(1-x)\log x+\log^2(1-x)Li_2(1-x)-2\log(1-x)Li_3(1-x)+2Li_4(1-x)+2Li_4(x)-2\varsigma(4)$$

and (3.110ea)

$$\sum_{n=1}^{\infty}\frac{\left(H_n^{(1)}\right)^2}{n^2}x^n = -\frac{1}{3}\log^3(1-x)\log x - \log^2(1-x)Li_2(1-x) + 2\log(1-x)Li_3(1-x)$$

$$-2Li_4(1-x) + 2\varsigma(4) + Li_4(x) + \frac{1}{2}\left[Li_2(x)\right]^2$$

we see that for $0 \le x \le 1$

(3.211b) $$2\sum_{n=1}^{\infty}\frac{H_n^{(1)}}{n^3}x^n + \sum_{n=1}^{\infty}\frac{H_n^{(2)}}{n^2}x^n = 3Li_4(x) + \frac{1}{2}\left[Li_2(x)\right]^2$$

Combining (3.211), (3.211i) and (3.211ii), it was disappointing to find that the tedious algebra simply reproduced Euler's dilogarithm identity.

In 2005 Choi and Srivastava [45aci] showed that

(3.211bi) $$\sum_{n=1}^{\infty}\frac{H_n^{(2)}}{2^n n^2} = \frac{1}{16}\varsigma(4) + \frac{1}{4}\varsigma(3)\log 2 - \frac{1}{4}\varsigma(2)\log^2 2 + \frac{1}{24}\log^4 2 + Li_4\left(\frac{1}{2}\right)$$

and hence from (3.211b) we may also determine similar expressions for $\sum_{n=1}^{\infty}\frac{H_n^{(1)}}{2^n n^3}$ and $\sum_{n=1}^{\infty}\frac{\left(H_n^{(1)}\right)^2}{2^n n^2}$.

As reported by Berndt [21, Part I, p.253], Ramanujan proved that

(3.211bii) $$h(1-x) - h\left(1-\frac{1}{x}\right) =$$

$$-\frac{1}{24}\log^4 x + \frac{1}{6}\log^3 x\log(1-x) + \varsigma(3)\log x - 2Li_4(x) + Li_3(x)\log x + 2\varsigma(4)$$

where $h(x) = \sum_{n=1}^{\infty}\frac{H_n^{(1)}}{(n+1)^3}x^{n+1}$

We see that $$\sum_{n=1}^{\infty}\frac{H_n^{(1)}}{(n+1)^3}x^{n+1} = \sum_{n=1}^{\infty}\frac{H_n^{(1)} + \frac{1}{n+1}}{(n+1)^3}x^{n+1} - \sum_{n=1}^{\infty}\frac{x^{n+1}}{(n+1)^4}$$

$$= \sum_{n=1}^{\infty}\frac{H_{n+1}^{(1)}}{(n+1)^3}x^{n+1} - \sum_{n=1}^{\infty}\frac{x^{n+1}}{(n+1)^4} = \sum_{n=1}^{\infty}\frac{H_n^{(1)}}{n^3}x^n - Li_4(t)$$

and therefore we have $h(x)=\sum_{n=1}^{\infty}\frac{H_n^{(1)}}{n^3}x^n-Li_4(x)$. This gives us

$$\sum_{n=1}^{\infty}\frac{H_n^{(1)}}{n^3}(1-x)^n-Li_4(1-x)-\sum_{n=1}^{\infty}\frac{H_n^{(1)}}{n^3}\left(1-\frac{1}{x}\right)^n+Li_4\left(1-\frac{1}{x}\right)=$$

$$-\frac{1}{24}\log^4 x+\frac{1}{6}\log^3 x\log(1-x)+\varsigma(3)\log x-2Li_4(x)+Li_3(x)\log x+2\varsigma(4)$$

Differentiating (3.211b) we obtain

(3.211c) $$2\sum_{n=1}^{\infty}\frac{H_n^{(1)}}{n^2}x^n+\sum_{n=1}^{\infty}\frac{H_n^{(2)}}{n}x^n=3Li_3(x)-Li_2(x)\log(1-x)$$

and this may also be derived by combining (3.211i) and (3.211ii).

Substituting Euler's dilogarithm identity (1.6c) in (3.211) we get

$$3Li_4(x)-3\log x\,Li_3(x)+Li_2(x)\log x\log(1-x)+\frac{1}{2}\left[Li_2(x)\right]^2=$$

$$2\sum_{n=1}^{\infty}\frac{H_n^{(1)}}{n^3}x^n-2\log x\sum_{n=1}^{\infty}\frac{H_n^{(1)}}{n^2}x^n+\sum_{n=1}^{\infty}\frac{H_n^{(2)}}{n^2}x^n-\log x\sum_{n=1}^{\infty}\frac{H_n^{(2)}}{n}x^n$$

Dividing this by $x$ and integrating we have

$$3Li_5(x)-3\int_0^x\frac{\log x\,Li_3(x)}{x}dx+\int_0^x\frac{Li_2(x)\log x\log(1-x)}{x}dx+\frac{1}{2}\int_0^x\frac{\left[Li_2(x)\right]^2}{x}dx=$$

$$2\sum_{n=1}^{\infty}\frac{H_n^{(1)}}{n^4}x^n-2\int_0^x\log x\sum_{n=1}^{\infty}\frac{H_n^{(1)}}{n^2}x^{n-1}dx+\sum_{n=1}^{\infty}\frac{H_n^{(2)}}{n^3}x^n-\int_0^x\log x\sum_{n=1}^{\infty}\frac{H_n^{(2)}}{n}x^{n-1}dx$$

We note that

$$\int_0^x\frac{\log x\,Li_3(x)}{x}dx=\log x\,Li_4(x)-Li_5(x)$$

$$\int_0^x\frac{Li_2(x)\log(1-x)}{x}\log x\,dx=-\frac{1}{2}\left[Li_2(x)\right]^2\log x+\frac{1}{2}\int_0^x\frac{\left[Li_2(x)\right]^2}{x}dx$$

Using (3.207) we have

$$\int_0^x \log x \sum_{n=1}^{\infty} \frac{H_n^{(1)}}{n^2} x^{n-1} dx = -\sum_{n=1}^{\infty} \frac{H_n^{(1)}}{n^4} x^n + \log x \sum_{n=1}^{\infty} \frac{H_n^{(1)}}{n^3} x^n$$

$$\int_0^x \log x \sum_{n=1}^{\infty} \frac{H_n^{(2)}}{n} x^{n-1} dx = -\sum_{n=1}^{\infty} \frac{H_n^{(2)}}{n^3} x^n + \log x \sum_{n=1}^{\infty} \frac{H_n^{(2)}}{n^2} x^n$$

and we then obtain

(3.211d) $$6Li_5(x) - 3\log x\, Li_4(x) - \frac{1}{2}\left[Li_2(x)\right]^2 \log x + \int_0^x \frac{\left[Li_2(x)\right]^2}{x} dx =$$

$$4\sum_{n=1}^{\infty} \frac{H_n^{(1)}}{n^4} x^n - 2\log x \sum_{n=1}^{\infty} \frac{H_n^{(1)}}{n^3} x^n + 2\sum_{n=1}^{\infty} \frac{H_n^{(2)}}{n^3} x^n - \log x \sum_{n=1}^{\infty} \frac{H_n^{(2)}}{n^2} x^n$$

Integrating (3.211b) we get

(3.211e) $$\int_0^t \frac{\left[Li_2(x)\right]^2}{x} dx = 4\sum_{n=1}^{\infty} \frac{H_n^{(1)}}{n^4} t^n + 2\sum_{n=1}^{\infty} \frac{H_n^{(2)}}{n^3} t^n - 6Li_5(t)$$

(the Wolfram Integrator cannot evaluate this integral) and with $t = 1$ we have

$$\int_0^1 \frac{\left[Li_2(x)\right]^2}{x} dx = 4\sum_{n=1}^{\infty} \frac{H_n^{(1)}}{n^4} + 2\sum_{n=1}^{\infty} \frac{H_n^{(2)}}{n^3} - 6\varsigma(5)$$

Inserting (3.211e) in (3.211d) we obtain

(3.211ei) $$6Li_5(x) - 3\log x\, Li_4(x) - \frac{1}{2}\left[Li_2(x)\right]^2 \log x + 4\sum_{n=1}^{\infty} \frac{H_n^{(1)}}{n^4} + 2\sum_{n=1}^{\infty} \frac{H_n^{(2)}}{n^3} - 6\varsigma(5) =$$

$$4\sum_{n=1}^{\infty} \frac{H_n^{(1)}}{n^4} x^n - 2\log x \sum_{n=1}^{\infty} \frac{H_n^{(1)}}{n^3} x^n + 2\sum_{n=1}^{\infty} \frac{H_n^{(2)}}{n^3} x^n - \log x \sum_{n=1}^{\infty} \frac{H_n^{(2)}}{n^2} x^n$$

Combining (3.211c) and (3.211d) we simply recover (3.211b).

Using (3.110i) and (3.110p)

$$\sum_{n=1}^{\infty} \frac{H_n^{(1)}}{n^4} = 3\varsigma(5) - \varsigma(2)\varsigma(3)$$

$$\sum_{n=1}^{\infty} \frac{H_n^{(2)}}{n^3} = 3\varsigma(2)\varsigma(3) - \frac{9}{2}\varsigma(5)$$

we see that

(3.211f) $$\int_0^1 \frac{[Li_2(x)]^2}{x}\,dx = 2\varsigma(2)\varsigma(3) - 3\varsigma(5)$$

which is in accordance with the result previously obtained by Freitas [69a].

We easily see using integration by parts that

(3.212) $$\int_0^t \frac{\log(1-z)\log z}{(1-z)z}\,dz = \log(1-t)Li_2(1-t) - \log t\, Li_2(t) - Li_3(1-t) + Li_3(t) + \varsigma(3)$$

and, alternatively, by multiplying (3.203) across by $\log z / z$, and then integrating we have

$$\int_0^t \frac{\log(1-z)\log z}{(1-z)z}\,dz = -\sum_{n=1}^\infty H_n^{(1)} \int_0^t z^{n-1}\log z\,dz$$

Using the integral again

$$\int_0^t z^{n-1}\log z\,dz = \left[-\frac{1}{n^2} + \frac{\log t}{n}\right] t^n$$

we therefore get

$$\int_0^t \frac{\log(1-z)\log z}{(1-z)z}\,dz = \sum_{n=1}^\infty H_n^{(1)} \left[\frac{1}{n^2} - \frac{\log t}{n}\right] t^n$$

In conjunction with (3.212), this then results in the identity

(3.213)

$$\sum_{n=1}^\infty \frac{H_n^{(1)}}{n^2} t^n - \log t \sum_{n=1}^\infty \frac{H_n^{(1)}}{n} t^n = \log(1-t)Li_2(1-t) - \log t\, Li_2(t) - Li_3(1-t) + Li_3(t) + \varsigma(3)$$

and using (3.206) we get yet another proof of (3.105c), namely

(3.213a) $$\sum_{n=1}^\infty \frac{H_n^{(1)}}{n^2} t^n = \frac{1}{2}\log^2(1-t)\log t + \log(1-t)Li_2(1-t) - Li_3(1-t) + Li_3(t) + \varsigma(3)$$

From (3.206) we can write

$$\sum_{n=1}^{\infty}(-1)^n \frac{H_n^{(1)}}{n} t^n = \frac{1}{2}\log^2(1+t) + Li_2(-t)$$

and accordingly we have

(3.213b) $$\sum_{n=1}^{\infty}(-1)^n \frac{H_n^{(1)}}{n}\int_0^x t^{n-1}\log t\, dt = \frac{1}{2}\int_0^x \frac{\log^2(1+t)}{t}\log t\, dt + \int_0^x \frac{Li_2(-t)}{t}\log t\, dt$$

which gives us

(3.213c)

$$\sum_{n=1}^{\infty}(-1)^{n+1} \frac{H_n^{(1)}}{n^3} x^n - \log x \sum_{n=1}^{\infty}(-1)^{n+1} \frac{H_n^{(1)}}{n^2} x^n = \frac{1}{2}\int_0^x \frac{\log^2(1+t)}{t}\log t\, dt + \int_0^x \frac{Li_2(-t)}{t}\log t\, dt$$

This reminded me of the result obtained by Rutledge and Douglass [116aa] in 1934

(3.214) $$\int_0^1 \frac{\log^2(1+t)}{t}\log t\, dt = \sum_{n=1}^{\infty}(-1)^{n+1}\frac{H_n^{(2)}}{(n+1)^2} - \frac{\pi^4}{288}$$

$$= -\sum_{n=1}^{\infty}(-1)^{n+1}\frac{H_n^{(2)}}{n^2} + \varsigma_a(4) - \frac{\pi^4}{288}$$

$$= -\sum_{n=1}^{\infty}(-1)^{n+1}\frac{H_n^{(2)}}{n^2} + \frac{9}{16}\varsigma(4)$$

The other integral in (3.213c) may be written as

$$\int_0^x \frac{Li_2(-t)}{t}\log t\, dt = \frac{1}{2}\log^2 x\, Li_2(-x) + \frac{1}{2}\int_0^x \frac{\log^2 t \log(1+t)}{t} dt$$

since, using the series definition of the dilogarithm, we have $Li_2'(-t) = -\dfrac{\log(1+t)}{t}$.

We now have from (3.213c)

(3.215) $$\sum_{n=1}^{\infty}(-1)^{n+1} \frac{H_n^{(1)}}{n^3} x^n - \log x \sum_{n=1}^{\infty}(-1)^{n+1} \frac{H_n^{(1)}}{n^2} x^n =$$

$$\frac{1}{2}\log^2 x\, Li_2(-x) + \frac{1}{2}\int_0^x \frac{\log^2(1+t)\log t}{t} dt + \frac{1}{2}\int_0^x \frac{\log^2 t \log(1+t)}{t} dt$$

With integration by parts we have

$$\int_0^x \frac{\log^2 t \log(1+t)}{t} dt = -\log^2 x\, Li_2(-x) + 2\int_0^x \frac{\log t\, Li_2(-t)}{t} dt$$

$$\int_0^x \frac{\log t\, Li_2(-t)}{t} dt = \log x\, Li_3(-x) - \int_0^x \frac{Li_3(-t)}{t} dt$$

and therefore we get

$$\int_0^x \frac{\log^2 t \log(1+t)}{t} dt = -\log^2 x\, Li_2(-x) + 2\log x\, Li_3(-x) - 2Li_4(-x)$$

The above integral could also be obtained by noting that

$$\int_0^x \frac{Li_2(-t)}{t} \log t\, dt = \log x\, Li_3(-x) - Li_4(-x)$$

We obtain from (3.215)

(3.215a)

$$\sum_{n=1}^\infty (-1)^{n+1} \frac{H_n^{(1)}}{n^3} x^n - \log x \sum_{n=1}^\infty (-1)^{n+1} \frac{H_n^{(1)}}{n^2} x^n = \log x\, Li_3(-x) - Li_4(-x) + \frac{1}{2}\int_0^x \frac{\log^2(1+t)\log t}{t} dt$$

With $x=1$ we get

$$\sum_{n=1}^\infty (-1)^{n+1} \frac{H_n^{(1)}}{n^3} = -Li_4(-1) + \frac{1}{2}\int_0^1 \frac{\log^2(1+t)\log t}{t} dt$$

and using (3.214) this becomes

$$= \frac{7}{8}\varsigma(4) - \frac{1}{2}\sum_{n=1}^\infty (-1)^{n+1} \frac{H_n^{(2)}}{n^2} + \frac{9}{32}\varsigma(4)$$

Therefore we get

(3.215b) $$2\sum_{n=1}^\infty (-1)^{n+1} \frac{H_n^{(1)}}{n^3} + \sum_{n=1}^\infty (-1)^{n+1} \frac{H_n^{(2)}}{n^2} = \frac{37}{16}\varsigma(4)$$

Using (3.220b) we see that

(3.216)

$$\mu_1 = \sum_{n=1}^{\infty} (-1)^{n+1} \frac{H_n^{(1)}}{n^3} = \frac{11}{4}\varsigma(4) + \frac{1}{2}\varsigma(2)\log^2 2 - \frac{1}{12}\log^4 2 - \frac{7}{4}\varsigma(3)\log 2 - 2Li_4(1/2)$$

and therefore we get

(3.216a)

$$\sum_{n=1}^{\infty} (-1)^{n+1} \frac{H_n^{(2)}}{n^2} = -\frac{51}{16}\varsigma(4) - \varsigma(2)\log^2 2 + \frac{1}{6}\log^4 2 + \frac{7}{2}\varsigma(3)\log 2 + 4Li_4(1/2)$$

We note that

$$\sum_{n=1}^{\infty} \frac{H_n^{(1)}}{n} \int_0^x t^{n-1} \log t \, dt = -\sum_{n=1}^{\infty} \frac{H_n^{(1)}}{n^3} x^n + \log x \sum_{n=1}^{\infty} \frac{H_n^{(1)}}{n^2} x^n$$

and dividing (3.213) by $t$ and carrying out a (by now near compulsory) integration gives us

(3.216b)

$$2\sum_{n=1}^{\infty} \frac{H_n^{(1)}}{n^3} x^n - \log x \sum_{n=1}^{\infty} \frac{H_n^{(1)}}{n^2} x^n = \int_0^x \left[\log(1-t)Li_2(1-t) - \log t \, Li_2(t) + \varsigma(3) - Li_3(1-t) + Li_3(t)\right] \frac{dt}{t}$$

$$= \int_0^x \frac{\log(1-t)Li_2(1-t)}{t} dt - \int_0^x \frac{\log t \, Li_2(t)}{t} dt + \int_0^x \frac{\varsigma(3) - Li_3(1-t)}{t} dt + \int_0^x \frac{Li_3(t)}{t} dt$$

We have

$$\int_0^x \frac{\log(1-t)Li_2(1-t)}{t} dt = \int_{1-x}^1 \frac{\log t \, Li_2(t)}{1-t} dt = \int_0^1 \frac{\log t \, Li_2(t)}{1-t} dt - \int_0^{1-x} \frac{\log t \, Li_2(t)}{1-t} dt$$

In (4.4.233c) we show that

$$\int_0^y \frac{\log t \, Li_2(t)}{1-t} dt = \sum_{n=1}^{\infty} \frac{1}{n^2} \sum_{k=1}^{n} \frac{y^k}{k^2} + \varsigma(2)Li_2(1-y) - \varsigma^2(2)$$

and therefore we get

$$\int_{1-x}^1 \frac{\log t \, Li_2(t)}{1-t} dt = -\sum_{n=1}^{\infty} \frac{1}{n^2} \sum_{k=1}^{n} \frac{(1-x)^k}{k^2} - \varsigma(2)Li_2(x) + \sum_{n=1}^{\infty} \frac{1}{n^2} \sum_{k=1}^{n} \frac{1}{k^2}$$

From (4.4.167s) we have $\sum_{n=1}^{\infty}\frac{1}{n^2}\sum_{k=1}^{n}\frac{1}{k^2}=\frac{7}{4}\varsigma(4)$ and hence we obtain

(3.217a) $$\int_{1-x}^{1}\frac{\log t\, Li_2(t)}{1-t}dt=-\sum_{n=1}^{\infty}\frac{1}{n^2}\sum_{k=1}^{n}\frac{(1-x)^k}{k^2}-\varsigma(2)Li_2(x)+\frac{7}{4}\varsigma(4)$$

We have from (3.209)

(3.217b) $$\int_{0}^{x}\frac{\log t\, Li_2(t)}{t}dt=\log x\, Li_3(x)-Li_4(x)$$

and from (4.4.168j)

(3.217c)

$$\int_{0}^{x}\frac{\varsigma(3)-Li_3(1-t)}{t}dt=\left[\varsigma(3)-Li_3(1-x)\right]\log x-\frac{1}{2}\left[Li_2(1-x)\right]^2+\frac{1}{2}\varsigma^2(2)$$

We therefore obtain from (3.216b)

(3.217d)

$$2\sum_{n=1}^{\infty}\frac{H_n^{(1)}}{n^3}x^n-\log x\sum_{n=1}^{\infty}\frac{H_n^{(1)}}{n^2}x^n=-\sum_{n=1}^{\infty}\frac{1}{n^2}\sum_{k=1}^{n}\frac{(1-x)^k}{k^2}-\varsigma(2)Li_2(x)+\frac{7}{4}\varsigma(4)$$

$$-\log x\, Li_3(x)+\left[\varsigma(3)-Li_3(1-x)\right]\log x-\frac{1}{2}\left[Li_2(1-x)\right]^2+\frac{1}{2}\varsigma^2(2)+2Li_4(x)$$

and using (3.213a) we get

(3.218)

$$2\sum_{n=1}^{\infty}\frac{H_n^{(1)}}{n^3}x^n=-\sum_{n=1}^{\infty}\frac{1}{n^2}\sum_{k=1}^{n}\frac{(1-x)^k}{k^2}-\varsigma(2)Li_2(x)+\frac{7}{4}\varsigma(4)-\frac{1}{2}\left[Li_2(1-x)\right]^2+\frac{1}{2}\varsigma^2(2)+2Li_4(x)$$

$$+\frac{1}{2}\log^2(1-x)\log^2 x+\log(1-x)\log x\, Li_2(1-x)+2\log x\,\varsigma(3)-2\log x\, Li_2(1-x)$$

With $x=1/2$ we get

$$2\sum_{n=1}^{\infty}\frac{H_n^{(1)}}{2^n n^3}=-\sum_{n=1}^{\infty}\frac{1}{n^2}\sum_{k=1}^{n}\frac{1}{2^k k^2}-\varsigma(2)Li_2(1/2)+\frac{7}{4}\varsigma(4)-\frac{1}{2}\left[Li_2(1/2)\right]^2+\frac{1}{2}\varsigma^2(2)+2Li_4(1/2)$$

$$+\frac{1}{2}\log^4 2+\log^2 2\,Li_2(1/2)-2\log 2\,\varsigma(3)+2\log 2\,Li_2(1/2)$$

Using (3.23) we see that

$$\sum_{n=1}^{\infty}\frac{1}{n^2}\sum_{k=1}^{n}\frac{1}{2^k k^2}=\sum_{n=1}^{\infty}\frac{1}{2^n n^2}\sum_{k=n}^{\infty}\frac{1}{k^2}=\sum_{n=1}^{\infty}\frac{1}{2^n n^2}\left[\varsigma(2)-H_{n-1}^{(2)}\right]$$

$$=\varsigma(2)\sum_{n=1}^{\infty}\frac{1}{2^n n^2}-\sum_{n=1}^{\infty}\frac{H_n^{(2)}}{2^n n^2}+\sum_{n=1}^{\infty}\frac{1}{2^n n^4}$$

$$=\varsigma(2)Li_2(1/2)-\sum_{n=1}^{\infty}\frac{H_n^{(2)}}{2^n n^2}+Li_4(1/2)$$

$$=\frac{1}{2}\varsigma^2(2)-\frac{1}{2}\varsigma(2)\log^2 2-\sum_{n=1}^{\infty}\frac{H_n^{(2)}}{2^n n^2}+Li_4(1/2)$$

and therefore we obtain

$$2\sum_{n=1}^{\infty}\frac{H_n^{(1)}}{2^n n^3}-\sum_{n=1}^{\infty}\frac{H_n^{(2)}}{2^n n^2}=\frac{1}{2}\varsigma(2)\log^2 2-\varsigma(2)Li_2(1/2)+\frac{7}{4}\varsigma(4)-\frac{1}{2}\left[Li_2(1/2)\right]^2+Li_4(1/2)$$

$$+\frac{1}{2}\log^4 2+\log^2 2\,Li_2(1/2)-2\log 2\,\varsigma(3)+2\log 2\,Li_2(1/2)$$

We saw in (3.217a) that

$$\int_0^x\frac{\log(1-t)Li_2(1-t)}{t}dt=-\sum_{n=1}^{\infty}\frac{1}{n^2}\sum_{k=1}^{n}\frac{(1-x)^k}{k^2}-\varsigma(2)Li_2(x)+\frac{7}{4}\varsigma(4)$$

and referring back to (3.111k)

$$\int_0^x\frac{\log(1-t)Li_2(1-t)}{t}dt=$$

$$2\sum_{n=1}^{\infty}\frac{H_n}{n^3}x^n+\log x\left[Li_2(x)\log x+Li_2\left(\frac{-x}{1-x}\right)\log x-\log(1-x)Li_2(1-x)\right]$$

$$-2Li_4(x)-\frac{1}{2}\left[Li_2(1-x)\right]^2+\frac{1}{2}\varsigma^2(2)$$

We therefore have

$$-\sum_{n=1}^{\infty}\frac{1}{n^2}\sum_{k=1}^{n}\frac{(1-x)^k}{k^2}-\varsigma(2)Li_2(x)+\frac{7}{4}\varsigma(4)=$$

$$2\sum_{n=1}^{\infty}\frac{H_n}{n^3}x^n+\log x\left[Li_2(x)\log x+Li_2\left(\frac{-x}{1-x}\right)\log x-\log(1-x)Li_2(1-x)\right]$$

$$-2Li_4(x)-\frac{1}{2}\left[Li_2(1-x)\right]^2+\frac{1}{2}\varsigma^2(2)$$

in agreement with (3.218) above.

We may obtain a further identity by differentiating (3.218).

We have from (3.203) by letting $z\to -z$

(3.218a) $$-\frac{\log(1+z)}{(1+z)}=\sum_{n=1}^{\infty}(-1)^n H_n^{(1)} z^n$$

and hence

$$\int_0^t \frac{\log^2 z\log(1+z)}{z(1+z)}dz=\sum_{n=1}^{\infty}(-1)^{n+1}H_n^{(1)}\int_0^t z^{n-1}\log^2 z\, dz$$

It is easily shown using integration by parts that

(3.218b) $$\int_0^t z^{n-1}\log^2 z\, dz=t^n\left[\frac{2}{n^3}-\frac{2\log t}{n^2}+\frac{\log^2 t}{n}\right]$$

and hence we get

(3.218c) $$\int_0^1 z^{n-1}\log^2 z\, dz=\frac{2}{n^3}$$

Therefore we have

(3.219)

$$\int_0^t \frac{\log^2 z\log(1+z)}{z(1+z)}dz=2\sum_{n=1}^{\infty}(-1)^{n+1}\frac{H_n^{(1)}}{n^3}t^n-2\log t\sum_{n=1}^{\infty}(-1)^{n+1}\frac{H_n^{(1)}}{n^2}t^n+\log^2 t\sum_{n=1}^{\infty}(-1)^{n+1}\frac{H_n^{(1)}}{n}t^n$$

This is a particular case of the following identity referred to by Flajolet and Salvy in [69]

(3.220) $$\mu_q = \sum_{n=1}^{\infty} (-1)^{n+1} \frac{H_n^{(1)}}{n^{2q+1}} = \frac{1}{(2q)!} \int_0^1 \frac{\log^{2q} z \log(1+z)}{z(1+z)} dz$$

where we have from De Doelder [55] and Sitaramachandrarao [120a]

(3.220a) $$\mu_0 = \frac{1}{2}\varsigma(2) - \frac{1}{2}\log^2 2$$

(3.220b) $$\mu_1 = \frac{11}{4}\varsigma(4) + \frac{1}{2}\varsigma(2)\log^2 2 - \frac{1}{12}\log^4 2 - \frac{7}{4}\varsigma(3)\log 2 - 2Li_4(1/2)$$

Equation (3.220b) is De Doelder's result [55] as corrected by Coffey in [45d] (it is correctly stated in [69]).

The above result (3.220) may be generalised by parametric differentiation of $\int_0^t z^{\alpha-1}\, dz = \frac{t^\alpha}{\alpha}$ so that we have

$$\frac{d^p}{d\alpha^p} \int_0^t z^{\alpha-1}\, dz = \int_0^t z^{\alpha-1} \log^p z\, dz = \frac{d^p}{d\alpha^p}\left(\frac{t^\alpha}{\alpha}\right)$$

Using Leibniz's theorem we have

$$\frac{d^p}{d\alpha^p}\left(\frac{t^\alpha}{\alpha}\right) = \sum_{j=0}^{p} \binom{p}{j} \frac{d^j}{d\alpha^j}\left(t^\alpha\right) \frac{d^{p-j}}{d\alpha^{p-j}}\left(\frac{1}{\alpha}\right)$$

$$= \sum_{j=0}^{p} \binom{p}{j} (-1)^{p-j} \frac{t^\alpha}{\alpha^{p-j+1}} (p-j)!\, \log^j t$$

Hence we have from (3.218a)

(3.220c) $$\int_0^t \frac{\log^p z \log(1+z)}{z(1+z)} dz = \sum_{n=1}^{\infty} (-1)^{n+1} H_n^{(1)} \int_0^t z^{n-1} \log^p z\, dz$$

$$= \sum_{n=1}^{\infty} (-1)^{n+1} H_n^{(1)} \sum_{j=0}^{p} \binom{p}{j} (-1)^{p-j} \frac{t^n}{n^{p-j+1}} (p-j)!\, \log^j t$$

and we therefore get from (3.219)

(3.220d) $$2\sum_{n=1}^{\infty} (-1)^{n+1} \frac{H_n^{(1)}}{n^3} t^n - 2\log t \sum_{n=1}^{\infty} (-1)^{n+1} \frac{H_n^{(1)}}{n^2} t^n + \log^2 t \sum_{n=1}^{\infty} (-1)^{n+1} \frac{H_n^{(1)}}{n} t^n$$

$$= \sum_{n=1}^{\infty} (-1)^{n+1} \frac{H_n^{(1)}}{n^3} t^n \sum_{j=0}^{2} \binom{p}{j} (-1)^{p-j} n^j (p-j)! \log^j t$$

In the case where $t = 1$ we have

$$\int_0^1 z^{n-1} \log^p z \, dz = (-1)^p \frac{p!}{n^{p+1}}$$

Hence we have from (3.220c)

(3.220e) $$\int_0^1 \frac{\log^p z \log(1+z)}{z(1+z)} dz = \sum_{n=1}^{\infty} (-1)^{n+1} H_n^{(1)} \int_0^1 z^{n-1} \log^p z \, dz$$

$$= \frac{(-1)^p}{p!} \sum_{n=1}^{\infty} (-1)^{n+1} \frac{H_n^{(1)}}{n^{p+1}}$$

Differentiating (3.202) we get

(3.220f) $$\frac{\log^2(1-z)}{(1-z)^x} = \sum_{n=1}^{\infty} \frac{f_n''(x)}{n!} z^n$$

and we have

$$f_n''(x) = -x(x+1)(x+2)...(x+n-1)\left[\frac{1}{x^2} + \frac{1}{(x+1)^2} ... + \frac{1}{(x+n-1)^2}\right]$$

$$+x(x+1)(x+2)...(x+n-1)\left[\frac{1}{x} + \frac{1}{x+1} ... + \frac{1}{x+n-1}\right]^2$$

Thus we have $f_n''(1) = -n! H_n^{(2)} + n!\left[H_n^{(1)}\right]^2$ and therefore

(3.221) $$\frac{\log^2(1-z)}{(1-z)} = \sum_{n=1}^{\infty} \left(\left[H_n^{(1)}\right]^2 - H_n^{(2)}\right) z^n$$

This is in agreement with (3.34) and (3.35) which are reproduced below

(3.221a) $$\frac{Li_2(x)}{1-x} = \sum_{n=1}^{\infty} H_n^{(2)} x^n \qquad , x \in [0,1)$$

(3.221b) $$\frac{\log^2(1-x) + Li_2(x)}{1-x} = \sum_{n=1}^{\infty} \left(H_n^{(1)}\right)^2 x^n \qquad , x \in [0,1)$$

From (3.221a) we get

$$\int_0^t \frac{Li_2(x)\log x}{(1-x)x}dx = \sum_{n=1}^{\infty} H_n^{(2)} \int_0^t x^{n-1}\log x dx$$

Using the integral $\int_0^t x^{n-1}\log x\, dx = \left[-\frac{1}{n^2}+\frac{\log t}{n}\right]t^n$ we obtain

(3.222) $$\int_0^t \frac{Li_2(x)\log x}{(1-x)x}dx = -\sum_{n=1}^{\infty}\frac{H_n^{(2)}}{n^2}t^n + \log x\sum_{n=1}^{\infty}\frac{H_n^{(2)}}{n}t^n$$

We have using partial fractions

(3.223) $$\int_0^t \frac{Li_2(x)\log x}{(1-x)x}dx = \int_0^t \frac{Li_2(x)\log x}{x}dx + \int_0^t \frac{Li_2(x)\log x}{(1-x)}dx$$

and we have already seen from (3.209) that

$$\int_0^t \frac{Li_2(x)\log x}{x}dx = \log t\, Li_3(t) - Li_4(t)$$

We now address the third integral in (3.223)

$$\int_0^t Li_2(x)\frac{\log x}{(1-x)}dx = Li_2(t)Li_2(1-t) + \int_0^t \frac{Li_2(1-x)\log(1-x)}{x}dx$$

We have from (4.4.167ma)

$$\int_0^t \frac{\log(1-x)Li_2(1-x)}{x}dx + \int_0^t \frac{\log x\log^2(1-x)}{x}dx = -\varsigma(2)Li_2(t) + \frac{1}{2}\left[Li_2(t)\right]^2$$

and therefore we obtain

(3.224)

$$\int_0^t \frac{Li_2(x)\log x}{(1-x)}dx = Li_2(t)Li_2(1-t) - \varsigma(2)Li_2(t) + \frac{1}{2}\left[Li_2(t)\right]^2 - \int_0^t \frac{\log x\log^2(1-x)}{x}dx$$

Hence we get

(3.225) $$\int_0^t \frac{Li_2(x)\log x}{(1-x)x}dx =$$

$$= \log t\, Li_3(t) - Li_4(t) + Li_2(t)Li_2(1-t) - \varsigma(2)Li_2(t) + \frac{1}{2}\left[Li_2(t)\right]^2 - \int_0^t \frac{\log x \log^2(1-x)}{x}dx$$

Therefore we have obtained the following integral using (3.222) and (3.225)

(3.226) $$\int_0^t \frac{\log x \log^2(1-x)}{x}dx =$$

$$= \log t\, Li_3(t) - Li_4(t) + Li_2(t)Li_2(1-t) - \varsigma(2)Li_2(t) + \frac{1}{2}\left[Li_2(t)\right]^2 + \sum_{n=1}^{\infty}\frac{H_n^{(2)}}{n^2}t^n - \log t \sum_{n=1}^{\infty}\frac{H_n^{(2)}}{n}t^n$$

and this is valid for $t \in [0,1]$.

Letting $t = 1$ we obtain a proof of (4.4.167q)

$$\int_0^1 \frac{\log x \log^2(1-x)}{x}dx = = -\varsigma(4) - \frac{1}{2}\varsigma^2(2) + \sum_{n=1}^{\infty}\frac{H_n^{(2)}}{n^2} = -\frac{1}{2}\varsigma(4)$$

where $\sum_{n=1}^{\infty}\frac{H_n^{(2)}}{n^2} = \frac{7}{4}\varsigma(4)$ and $\varsigma^2(2) = \frac{5}{2}\varsigma(4)$ .

We also have from (4.4.168f)

$$\int_0^t \frac{\log x \log^2(1-x)}{x}dx = \frac{1}{2}\log t\left[\log^2(1-t)\log t + 2\log(1-t)Li_2(1-t) - 2Li_3(1-t) + 2\varsigma(3)\log t\right]$$

$$-\int_0^t \frac{Li_2(1-x)\log(1-x)}{x}dx - \int_0^t \frac{\varsigma(3) - Li_3(1-x)}{x}dx$$

and using (3.217a) and (3.217c) this becomes

(3.226a)

$$\int_0^t \frac{\log x \log^2(1-x)}{x}dx =$$

$$\frac{1}{2}\log t\left[\log^2(1-t)\log t + 2\log(1-t)Li_2(1-t) - 2Li_3(1-t) + 2\varsigma(3)\log t\right]$$

$$+\sum_{n=1}^{\infty}\frac{1}{n^2}\sum_{k=1}^{n}\frac{(1-t)^k}{k^2}+\varsigma(2)Li_2(t)-\frac{7}{4}\varsigma(4)-\left[\varsigma(3)-Li_3(1-t)\right]\log t+\frac{1}{2}\left[Li_2(1-t)\right]^2-\frac{1}{2}\varsigma^2(2)$$

$$=\frac{1}{2}\log^2(1-t)\log^2 t+\log t\log(1-t)Li_2(1-t)$$

$$+\sum_{n=1}^{\infty}\frac{1}{n^2}\sum_{k=1}^{n}\frac{(1-t)^k}{k^2}+\varsigma(2)Li_2(t)-\frac{7}{4}\varsigma(4)+\frac{1}{2}\left[Li_2(1-t)\right]^2-\frac{1}{2}\varsigma^2(2)$$

We may write

$$\sum_{n=1}^{\infty}\frac{1}{n^2}\sum_{k=1}^{n}\frac{(1-t)^k}{k^2}=\sum_{n=1}^{\infty}\frac{(1-t)^n}{n^2}\sum_{k=n}^{\infty}\frac{1}{k^2}=\sum_{n=1}^{\infty}\frac{(1-t)^n}{n^2}\left[\varsigma(2)-H_{n-1}^{(2)}\right]$$

$$=\varsigma(2)Li_2(1-t)-\sum_{n=1}^{\infty}\frac{H_n^{(2)}}{n^2}(1-t)^n+\sum_{n=1}^{\infty}\frac{1}{n^3}(1-t)^n$$

An alternative derivation of (4.4.168f) is shown below. We have from (4.4.100gii)

$$\int_0^t\frac{\log^2(1-x)}{x}dx=\log t\log^2(1-t)+2\log(1-t)\,Li_2(1-t)-2Li_3(1-t)+2\varsigma(3)$$

and using integration by parts we get

$$\int_0^t\frac{\log x\log^2(1-x)}{x}dx=\log t\left[\log t\log^2(1-t)+2\log(1-t)\,Li_2(1-t)-2Li_3(1-t)+2\varsigma(3)\right]$$

$$-\int_0^t\left[\log t\log^2(1-t)+2\log(1-t)\,Li_2(1-t)-2Li_3(1-t)+2\varsigma(3)\right]\frac{dx}{x}$$

Therefore

$$\int_0^t\frac{\log x\log^2(1-x)}{x}dx=\frac{1}{2}\log t\left[\log t\log^2(1-t)+2\log(1-t)\,Li_2(1-t)-2Li_3(1-t)+2\varsigma(3)\right]$$

$$-\int_0^t\frac{\log(1-t)\,Li_2(1-t)}{x}dx-\int_0^t\frac{\varsigma(3)-Li_3(1-t)}{x}dx$$

Equating (3.226) and (3.226a) we get

$$\log t\, Li_3(t) - Li_4(t) + Li_2(t)Li_2(1-t) - \varsigma(2)Li_2(t) + \frac{1}{2}\left[Li_2(t)\right]^2 + \sum_{n=1}^{\infty}\frac{H_n^{(2)}}{n^2}t^n - \log t\sum_{n=1}^{\infty}\frac{H_n^{(2)}}{n}t^n =$$

$$\frac{1}{2}\log t\left[\log^2(1-t)\log t + 2\log(1-t)Li_2(1-t) - 2Li_3(1-t) + 2\varsigma(3)\right]$$

$$+\sum_{n=1}^{\infty}\frac{1}{n^2}\sum_{k=1}^{n}\frac{(1-t)^k}{k^2} + \varsigma(2)Li_2(t) - \frac{7}{4}\varsigma(4) - \left[\varsigma(3) - Li_3(1-t)\right]\log t + \frac{1}{2}\left[Li_2(1-t)\right]^2 - \frac{1}{2}\varsigma^2(2)$$

and this may be simplified to

(3.226b) $$\sum_{n=1}^{\infty}\frac{H_n^{(2)}}{n^2}t^n - \log t\sum_{n=1}^{\infty}\frac{H_n^{(2)}}{n}t^n =$$

$$\frac{1}{2}\log^2(1-t)\log^2 t + \log(1-t)\log t\, Li_2(1-t) - \log t\, Li_3(t) + Li_4(t) - Li_2(t)Li_2(1-t)$$

$$+\sum_{n=1}^{\infty}\frac{1}{n^2}\sum_{k=1}^{n}\frac{(1-t)^k}{k^2} + 2\varsigma(2)Li_2(t) - \frac{7}{4}\varsigma(4) + \frac{1}{2}\left[Li_2(1-t)\right]^2 - \frac{1}{2}\left[Li_2(t)\right]^2 - \frac{1}{2}\varsigma^2(2)$$

As we will see later in (4.4.167p) the Wolfram Integrator miraculously produces the result

(3.226c) $$\int_0^t \frac{\log x\log^2(1-x)}{x}dx =$$

$$\frac{1}{2}\log^2(1-t)\log^2 t + \frac{1}{12}\log^4 t - \log^2(1-t)\log^2 t + \frac{2}{3}\log(1-t)\log^3 t$$

$$-\left[\log(1-t) + \frac{1}{3}\log t\right]\log^2 t\log\left(\frac{t}{1-t}\right) + \frac{1}{2}\log^2 t\left[\log\left(\frac{t}{1-t}\right)\right]^2$$

$$-\frac{1}{4}\left[\log\left(\frac{t}{1-t}\right)\right]^4 + \log^2(1-t)Li_2(1-t) - \log^2 t\, Li_2(t)$$

$$-\left[\log\left(\frac{x}{1-x}\right)\right]^2 Li_2\left(\frac{-t}{1-t}\right) - 2\log(1-t)Li_3(1-t) + 2\log t\, Li_3(t)$$

$$+2\log\left(\frac{t}{1-t}\right)Li_3\left(\frac{-t}{1-t}\right) + 2\left[Li_4(1-t) - Li_4(t) - Li_4\left(\frac{-t}{1-t}\right)\right] - 2\varsigma(4)$$

(the machine generated proof may of course be easily verified by differentiating).

The presence of so many terms involving $\log x$ initially made me think that the integral was not convergent at $x = 0$ but closer inspection reveals the net sum of the terms only involving $\log^4 x$ is zero (the factor of $\log(1-x)$ elsewhere ensures convergence for the other powers involving $\log x$). However, it is not immediately apparent to me how the integrated part behaves as $x \to 1$, having regard to the polylogarithmic terms involving $(-x/1-x)$: perhaps the Wolfram Integrator output is only valid for $x < 1$ because the Maclaurin series expansion for $\log(1-x)$ is not convergent at $x = 1$.

Equating (3.226a) and (3.226c) gives us a rather complex expression for $\sum_{n=1}^{\infty} \frac{1}{n^2} \sum_{k=1}^{n} \frac{(1-t)^k}{k^2}$.

From (3.226c) we get much cancellation with $t = 1/2$ and end up with a concise result

(3.227) $$\int_0^{1/2} \frac{\log x \log^2(1-x)}{x} dx = \frac{1}{4}\log^4 2 + 2Li_4(-1) - 2\varsigma(4) = \frac{1}{4}\log^4 2 - \frac{15}{4}\varsigma(4)$$

where we have used (4.4.67) to show that $Li_4(-1) = -\frac{7}{8}\varsigma(4)$ .

By differentiation we may easily prove that

(3.228)

$$\int_0^t \frac{Li_2(x)}{(1-x)x} dx = \log t \log^2(1-t) + Li_2(t)\log(1-t) - 2\varsigma(2)\log(1-t) + 2Li_3(1-t) + Li_3(t) - 2\varsigma(3)$$

Alternatively, using (3.34) we have

$$\int_0^t \frac{Li_2(x)}{(1-x)x} dx = \sum_{n=1}^{\infty} \frac{H_n^{(2)}}{n} t^n$$

and hence we get

(3.229)

$$\sum_{n=1}^{\infty} \frac{H_n^{(2)}}{n} t^n = \log t \log^2(1-t) + Li_2(t)\log(1-t) - 2\varsigma(2)\log(1-t) + 2Li_3(1-t) + Li_3(t) - 2\varsigma(3)$$

Using Euler's dilogarithm identity we note that (3.229) is equivalent to (3.211i).

Substituting (3.229) in (3.226) we obtain

(3.230) $$\int_0^t \frac{\log x \log^2(1-x)}{x}dx =$$

$$= \log t\, Li_3(t) - Li_4(t) + Li_2(t)Li_2(1-t) - \varsigma(2)Li_2(t) + \frac{1}{2}\left[Li_2(t)\right]^2 + \sum_{n=1}^{\infty}\frac{H_n^{(2)}}{n^2}t^n$$

$$-\log^2 t \log^2(1-t) - Li_2(t)\log t \log(1-t) + 2\varsigma(2)\log t \log(1-t) - 2\log t\, Li_3(1-t)$$

$$-Li_3(t)\log t + 2\varsigma(3)\log t$$

$$= -Li_4(t) + Li_2(t)Li_2(1-t) - \varsigma(2)Li_2(t) + \frac{1}{2}\left[Li_2(t)\right]^2 + \sum_{n=1}^{\infty}\frac{H_n^{(2)}}{n^2}t^n$$

$$-\log^2 t \log^2(1-t) - Li_2(t)\log t \log(1-t) + 2\varsigma(2)\log t \log(1-t) - 2\log t\, Li_3(1-t) + 2\varsigma(3)\log t$$

Therefore we get

(3.231) $$\int_0^1 \frac{\log x \log^2(1-x)}{x}dx = -\varsigma(4) - \frac{1}{2}\varsigma^2(2) + \sum_{n=1}^{\infty}\frac{H_n^{(2)}}{n^2} = -\frac{1}{2}\varsigma(4)$$

The above analysis is equivalent to (3.226) as shown below

$$\int_0^t \frac{\log x \log^2(1-x)}{x}dx =$$

$$= \log t\, Li_3(t) - Li_4(t) + Li_2(t)Li_2(1-t) - \varsigma(2)Li_2(t) + \frac{1}{2}\left[Li_2(t)\right]^2 + \sum_{n=1}^{\infty}\frac{H_n^{(2)}}{n^2}t^n - \log t\sum_{n=1}^{\infty}\frac{H_n^{(2)}}{n}t^n$$

We have

$$\sum_{n=1}^{\infty}\frac{1}{n^2}\sum_{k=1}^{n}\frac{(1-t)^k}{k^2} = \sum_{n=1}^{\infty}\frac{(1-t)^n}{n^2}\sum_{k=n}^{\infty}\frac{1}{k^2} = \sum_{n=1}^{\infty}\frac{(1-t)^n}{n^2}\left[\varsigma(2) - H_{n-1}^{(2)}\right]$$

$$= \varsigma(2)Li_2(1-t) - \sum_{n=1}^{\infty}\frac{(1-t)^n H_{n-1}^{(2)}}{n^2}$$

$$= \varsigma(2)Li_2(1-t) - \sum_{n=1}^{\infty}\frac{(1-t)^n H_n^{(2)}}{n^2} + \sum_{n=1}^{\infty}\frac{(1-t)^n}{n^4}$$

(3.232) $$\sum_{n=1}^{\infty}\frac{1}{n^2}\sum_{k=1}^{n}\frac{(1-t)^k}{k^2} = \varsigma(2)Li_2(1-t) - \sum_{n=1}^{\infty}\frac{H_n^{(2)}}{n^2}(1-t)^n + Li_4(1-t)$$

We have in particular for $t = 1/2$

$$\sum_{n=1}^{\infty}\frac{1}{n^2}\sum_{k=1}^{n}\frac{1}{2^k k^2}=\varsigma(2)Li_2(1/2)-\sum_{n=1}^{\infty}\frac{H_n^{(2)}}{2^n n^2}+Li_4(1/2)$$

and for $t=0$ we rather easily obtain

$$2\sum_{n=1}^{\infty}\frac{H_n^{(2)}}{n^2}=\varsigma^2(2)+\varsigma(4)$$

as required by (4.4.232a).

Integrating (3.221) we get

$$-\frac{1}{3}\log^3(1-t)=\frac{1}{t}\sum_{n=1}^{\infty}\frac{\left(-H_n^{(2)}+\left[H_n^{(1)}\right]^2\right)}{n+1}t^n$$

We see that

$$\sum_{n=1}^{\infty}\frac{\left(-H_n^{(2)}+\left[H_n^{(1)}\right]^2\right)}{n+1}t^n=\sum_{n=1}^{\infty}\frac{\left[H_{n+1}^{(1)}-\frac{1}{n+1}\right]^2}{n+1}t^n-\sum_{n=1}^{\infty}\frac{\left[H_{n+1}^{(2)}-\frac{1}{(n+1)^2}\right]}{n+1}t^n$$

$$=\sum_{n=1}^{\infty}\frac{\left[H_{n+1}^{(1)}\right]^2}{n+1}t^n-2\sum_{n=1}^{\infty}\frac{H_{n+1}^{(1)}}{(n+1)^2}t^n+\sum_{n=1}^{\infty}\frac{t^n}{(n+1)^3}-\sum_{n=1}^{\infty}\frac{H_{n+1}^{(2)}}{n+1}t^n+\sum_{n=1}^{\infty}\frac{t^n}{(n+1)^3}$$

$$=\frac{1}{t}\left[\sum_{n=1}^{\infty}\frac{\left[H_n^{(1)}\right]^2}{n}t^n-2\sum_{n=1}^{\infty}\frac{H_n^{(1)}}{n^2}t^n+2Li_3(t)-\sum_{n=1}^{\infty}\frac{H_n^{(2)}}{n}t^n\right]$$

We therefore get

(3.233) $$-\frac{1}{3}\log^3(1-t)=\sum_{n=1}^{\infty}\frac{\left[H_n^{(1)}\right]^2}{n}t^n-2\sum_{n=1}^{\infty}\frac{H_n^{(1)}}{n^2}t^n+2Li_3(t)-\sum_{n=1}^{\infty}\frac{H_n^{(2)}}{n}t^n$$

which we recall from (3.106). We may easily deduce the value of the Stirling number $s(n,3)$ from (3.233).

Dividing (3.233) by $t$ and integrating results in

(3.234) $$-\frac{1}{3}\int_0^x\frac{\log^3(1-t)}{t}dt=\sum_{n=1}^{\infty}\frac{\left[H_n^{(1)}\right]^2}{n^2}x^n-2\sum_{n=1}^{\infty}\frac{H_n^{(1)}}{n^3}x^n+2Li_4(x)-\sum_{n=1}^{\infty}\frac{H_n^{(2)}}{n^2}x^n$$

We have

$$\int_0^x \frac{\log^3(1-t)}{t}\,dt = \log x \log^3(1-x) + 3Li_2(1-x)\log^2(1-x)$$

$$-6Li_3(1-x)\log(1-x) + 6Li_4(1-x) - 6\varsigma(4)$$

and therefore we obtain

(3.235) $$6\sum_{n=1}^\infty \frac{H_n^{(1)}}{n^3}x^n - 3\sum_{n=1}^\infty \frac{\left[H_n^{(1)}\right]^2}{n^2}x^n - 6Li_4(x) + 3\sum_{n=1}^\infty \frac{H_n^{(2)}}{n^2}x^n =$$

$$\log x \log^3(1-x) + 3Li_2(1-x)\log^2(1-x) - 6Li_3(1-x)\log(1-x) + 6Li_4(1-x) - 6\varsigma(4)$$

We now try something different: multiplying (3.233) by $\log t/t$ , and then integrating, we obtain (on the assumption that it is valid to interchange the order of integration and summation)

(3.236) $$-\frac{1}{3}\int_0^x \frac{\log^3(1-t)\log t}{t}\,dt =$$

$$\sum_{n=1}^\infty \frac{\left[H_n^{(1)}\right]^2}{n}\int_0^x t^{n-1}\log t dt - 2\sum_{n=1}^\infty \frac{H_n^{(1)}}{n^2}\int_0^x t^{n-1}\log t dt + 2\int_0^x \frac{\log t}{t}Li_3(t)dt - \sum_{n=1}^\infty \frac{H_n^{(2)}}{n}\int_0^x t^{n-1}\log t dt$$

Using the following integrals (which may be easily derived by the parametric differentiation of $\int_0^x t^{\alpha-1}dt$ )

(3.237i) $$\int_0^x t^{n-1}\log t\,dt = \left[-\frac{1}{n^2} + \frac{\log x}{n}\right]x^n$$

(3.237ii) $$\int_0^x t^{n-1}\log^2 t\,dt = \left[\frac{2}{n^3} - \frac{2\log x}{n^2} + \frac{\log^2 x}{n}\right]x^n$$

(3.237iii) $$\int_0^x t^{n-1}\log^3 t\,dt = \left[-\frac{6}{n^4} + \frac{6\log x}{n^3} - \frac{3\log^2 x}{n^2} + \frac{\log^3 x}{n}\right]x^n$$

we have

(3.238) $$-\frac{1}{3}\int_0^x \frac{\log^3(1-t)\log t}{t}\,dt =$$

$$\sum_{n=1}^{\infty}\left[\frac{\left[H_n^{(1)}\right]^2}{n}-2\frac{H_n^{(1)}}{n^2}+\frac{H_n^{(2)}}{n}\right]\left[-\frac{1}{n^2}+\frac{\log x}{n}\right]x^n+2\int_0^x\frac{\log t}{t}Li_3(t)dt$$

Integration by parts gives us

$$\int_0^x\frac{\log t}{t}Li_3(t)dt=\frac{1}{2}\log^2 x\,Li_3(x)-\frac{1}{2}\int_0^x\log^2 t\frac{Li_2(t)}{t}dt$$

$$\int_0^x\frac{\log^2 t}{t}Li_2(t)dt=\frac{1}{3}\log^3 x\,Li_2(x)+\frac{1}{3}\int_0^x\frac{\log^3 t\log(1-t)}{t}dt$$

and hence we get

$$\int_0^x\frac{\log t}{t}Li_3(t)dt=\frac{1}{2}\log^2 x\,Li_3(x)-\frac{1}{6}\log^3 x\,Li_2(x)-\frac{1}{6}\int_0^x\frac{\log^3 t\log(1-t)}{t}dt$$

Substituting the above in (3.238) results in the cancellation of the term involving $\int_0^x\frac{\log^3 t\log(1-t)}{t}dt$ and we obtain

(3.239)

$$2\sum_{n=1}^{\infty}\frac{H_n^{(1)}}{n^4}x^n-\sum_{n=1}^{\infty}\frac{H_n^{(2)}}{n^3}x^n-\sum_{n=1}^{\infty}\frac{\left[H_n^{(1)}\right]^2}{n^3}x^n+\log x\left[\sum_{n=1}^{\infty}\frac{\left[H_n^{(1)}\right]^2}{n^2}x^n-2\sum_{n=1}^{\infty}\frac{H_n^{(1)}}{n^3}x^n+\sum_{n=1}^{\infty}\frac{H_n^{(2)}}{n^2}x^n\right]$$

$$+\log^2 x\,Li_3(x)-\frac{1}{3}\log^3 x\,Li_2(x)=0$$

With $x=1$ we have

(3.240) $$\sum_{n=1}^{\infty}\frac{\left[H_n^{(1)}\right]^2}{n^3}-2\sum_{n=1}^{\infty}\frac{H_n^{(1)}}{n^4}+\sum_{n=1}^{\infty}\frac{H_n^{(2)}}{n^3}=0$$

and a more complex expression arises, for example, with $x=1/2$. Further identities may be found by using (3.237ii) and (3.237iii).

Multiplying (3.239) by $1/x$ , and then integrating, we obtain (on the assumption that it is valid to interchange the order of integration and summation)

$$2\sum_{n=1}^{\infty}\frac{H_n^{(1)}}{n^5}x^n-\sum_{n=1}^{\infty}\frac{H_n^{(2)}}{n^4}x^n-\sum_{n=1}^{\infty}\frac{\left[H_n^{(1)}\right]^2}{n^4}x^n-\sum_{n=1}^{\infty}\frac{\left[H_n^{(1)}\right]^2}{n^4}x^n+\log x\sum_{n=1}^{\infty}\frac{\left[H_n^{(1)}\right]^2}{n^3}x^n$$

$$+2\sum_{n=1}^{\infty}\frac{H_n^{(1)}}{n^5}x^n - 2\log x\sum_{n=1}^{\infty}\frac{H_n^{(1)}}{n^4}x^n - \sum_{n=1}^{\infty}\frac{H_n^{(2)}}{n^4}x^n + \log x\sum_{n=1}^{\infty}\frac{H_n^{(2)}}{n^4}x^n$$

$$+\int_0^x \frac{\log^2 t}{t}Li_3(t)dt - \frac{1}{3}\int_0^x \frac{\log^3 t}{t}Li_2(t)dt = 0$$

Integration by parts gives us

$$\int_0^x \frac{\log^2 t}{t}Li_3(t)dt = \frac{1}{3}\log^3 x\, Li_3(x) - \frac{1}{3}\int_0^x \log^3 t\frac{Li_2(t)}{t}dt$$

$$\int_0^x \frac{\log^3 t}{t}Li_2(t)dt = \frac{1}{4}\log^4 x\, Li_2(x) + \frac{1}{4}\int_0^x \frac{\log^4 t\log(1-t)}{t}dt$$

and hence we get

$$\int_0^x \frac{\log^2 t}{t}Li_3(t)dt = \frac{1}{3}\log^3 x\, Li_3(x) - \frac{1}{12}\log^4 x\, Li_2(x) - \frac{1}{12}\int_0^x \frac{\log^4 t\log(1-t)}{t}dt$$

Upon inspection we see that

$$\int_0^x \frac{\log^2 t}{t}Li_3(t)dt - \frac{1}{3}\int_0^x \frac{\log^3 t}{t}Li_2(t)dt = \frac{1}{3}\log^3 x\, Li_3(x) - \frac{1}{6}\log^4 x\, Li_2(x) - \frac{1}{6}\int_0^x \frac{\log^4 t\log(1-t)}{t}dt$$

Integration by parts also gives us

$$\int_0^x \frac{\log^4 t\log(1-t)}{t}dt = \log^4 x\, Li_2(x) - 4\log^3 x\, Li_3(x) + 12\log^2 x\, Li_4(x) - 24\log x\, Li_5(x) + 24 Li_6(x)$$

and we therefore obtain

(3.240a)

$$4\sum_{n=1}^{\infty}\frac{H_n^{(1)}}{n^5}x^n - 2\sum_{n=1}^{\infty}\frac{H_n^{(2)}}{n^4}x^n - 2\sum_{n=1}^{\infty}\frac{\left[H_n^{(1)}\right]^2}{n^4}x^n + \log x\left[\sum_{n=1}^{\infty}\frac{\left[H_n^{(1)}\right]^2}{n^3}x^n - 2\sum_{n=1}^{\infty}\frac{H_n^{(1)}}{n^4}x^n + \sum_{n=1}^{\infty}\frac{H_n^{(2)}}{n^4}x^n\right]$$

$$+\log^3 x\, Li_3(x) - \frac{1}{3}\log^4 x\, Li_2(x) - 2\log^2 x\, Li_4(x) + 4\log x\, Li_5(x) - 4Li_6(x) = 0$$

With $x = 1$ we get

$$2\sum_{n=1}^{\infty}\frac{H_n^{(1)}}{n^5} - \sum_{n=1}^{\infty}\frac{H_n^{(2)}}{n^4} - \sum_{n=1}^{\infty}\frac{\left[H_n^{(1)}\right]^2}{n^4} = 2\varsigma(6)$$

From Flajolet and Salvy [69] we have

$$\sum_{n=1}^{\infty} \frac{H_n^{(2)}}{n^4} = \varsigma^2(3) - \frac{1}{3}\varsigma(6)$$

Georghiou and Philippou [69c] have shown that

$$\sum_{n=1}^{\infty} \frac{H_n^{(1)}}{n^k} = \left(1+\frac{k}{2}\right)\varsigma(k+1) - \frac{1}{2}\sum_{j=2}^{k-1}\varsigma(j)\varsigma(2k-j+2)$$

and hence we get

$$2\sum_{n=1}^{\infty} \frac{H_n^{(1)}}{n^5} = 2\varsigma(2)\varsigma(4) - \varsigma^2(3)$$

We therefore obtain

$$\text{(3.240b)} \qquad \sum_{n=1}^{\infty} \frac{\left[H_n^{(1)}\right]^2}{n^4} = -2\varsigma^2(3) - \frac{97}{24}\varsigma(6)$$

which is also reported by Flajolet and Salvy [69].

Multiplying (3.239) by $\log x / x$ , and then integrating, we obtain (on the assumption that it is valid to interchange the order of integration and summation)

$$-2\sum_{n=1}^{\infty} \frac{H_n^{(1)}}{n^6} x^n + \log x \sum_{n=1}^{\infty} \frac{H_n^{(1)}}{n^5} x^n + \sum_{n=1}^{\infty} \frac{H_n^{(2)}}{n^5} x^n - \log x \sum_{n=1}^{\infty} \frac{H_n^{(2)}}{n^4} x^n + \sum_{n=1}^{\infty} \frac{\left[H_n^{(1)}\right]^2}{n^5} x^n$$

$$-\log x \sum_{n=1}^{\infty} \frac{\left[H_n^{(1)}\right]^2}{n^4} x^n + 2\sum_{n=1}^{\infty} \frac{\left[H_n^{(1)}\right]^2}{n^5} x^n - 2\log x \sum_{n=1}^{\infty} \frac{\left[H_n^{(1)}\right]^2}{n^4} x^n + \log^2 x \sum_{n=1}^{\infty} \frac{\left[H_n^{(1)}\right]^2}{n^3} x^n$$

$$-4\sum_{n=1}^{\infty} \frac{H_n^{(1)}}{n^6} x^n + 4\log x \sum_{n=1}^{\infty} \frac{H_n^{(1)}}{n^5} x^n - 2\log^2 x \sum_{n=1}^{\infty} \frac{H_n^{(1)}}{n^4} x^n$$

$$+2\sum_{n=1}^{\infty} \frac{H_n^{(2)}}{n^5} x^n - 2\log x \sum_{n=1}^{\infty} \frac{H_n^{(2)}}{n^4} x^n + \log^2 x \sum_{n=1}^{\infty} \frac{H_n^{(2)}}{n^3} x^n$$

$$+\int_0^x \frac{\log^3 t}{t} Li_3(t)\,dt - \frac{1}{3}\int_0^x \frac{\log^4 t}{t} Li_2(t)\,dt = 0$$

and this easily simplifies to

$$-6\sum_{n=1}^{\infty}\frac{H_n^{(1)}}{n^6}x^n+3\sum_{n=1}^{\infty}\frac{H_n^{(2)}}{n^5}x^n+3\sum_{n=1}^{\infty}\frac{\left[H_n^{(1)}\right]^2}{n^5}x^n$$

$$+5\log x\left[\sum_{n=1}^{\infty}\frac{H_n^{(1)}}{n^5}x^n-3\sum_{n=1}^{\infty}\frac{\left[H_n^{(1)}\right]^2}{n^4}x^n-3\sum_{n=1}^{\infty}\frac{H_n^{(2)}}{n^4}x^n\right]$$

$$+\log^2 x\left[\sum_{n=1}^{\infty}\frac{\left[H_n^{(1)}\right]^2}{n^3}x^n-2\sum_{n=1}^{\infty}\frac{H_n^{(1)}}{n^4}x^n+\sum_{n=1}^{\infty}\frac{H_n^{(2)}}{n^3}x^n\right]$$

$$+\int_0^x\frac{\log^3 t}{t}Li_3(t)dt-\frac{1}{3}\int_0^x\frac{\log^4 t}{t}Li_2(t)dt=0$$

Integration by parts gives us

$$\int_0^x\frac{\log^3 t}{t}Li_3(t)dt=\frac{1}{4}\log^4 x\,Li_3(x)-\frac{1}{4}\int_0^x\log^4 t\frac{Li_2(t)}{t}dt$$

$$\int_0^x\frac{\log^4 t}{t}Li_2(t)dt=\frac{1}{5}\log^5 x\,Li_2(x)+\frac{1}{5}\int_0^x\frac{\log^4 t\log(1-t)}{t}dt$$

$$\int_0^x\frac{\log^4 t\,Li_2(t)}{t}dt=\frac{1}{5}\log^5 x\,Li_2(x)+\frac{1}{5}\int_0^x\frac{\log^5 t\log(1-t)}{t}dt$$

and hence we get

$$\int_0^x\frac{\log^3 t}{t}Li_3(t)dt=\frac{1}{4}\log^4 x\,Li_3(x)-\frac{1}{20}\log^5 x\,Li_2(x)-\frac{1}{20}\int_0^x\frac{\log^5 t\log(1-t)}{t}dt$$

Integration by parts also gives us

$$\int_0^x\frac{\log^5 t\log(1-t)}{t}dt=\log^5 x\,Li_2(x)-5\log^4 x\,Li_3(x)+20\log^3 x\,Li_4(x)-60\log^2 x\,Li_5(x)$$

$$+120\log x\,Li_6(x)-120Li_7(x)$$

Therefore we see that

$$\int_0^x\frac{\log^3 t}{t}Li_3(t)dt=-\frac{1}{10}\log^5 x\,Li_2(x)+\frac{1}{2}\log^4 x\,Li_3(x)-\log^3 x\,Li_4(x)+3\log^2 x\,Li_5(x)$$

$$-6\log x\, Li_6(x)+6Li_7(x)$$

$$\int_0^x \frac{\log^4 t\, Li_2(t)}{t}dt=\frac{2}{5}\log^5 x\, Li_2(x)-\log^4 x\, Li_3(x)+4\log^3 x\, Li_4(x)-12\log^2 x\, Li_5(x)$$

$$+24\log x\, Li_6(x)-24Li_7(x)$$

$$\int_0^x \frac{\log^3 t}{t}Li_3(t)dt-\frac{1}{3}\int_0^x \frac{\log^4 t}{t}Li_2(t)dt=$$

$$-\frac{7}{30}\log^5 x\, Li_2(x)+\frac{5}{6}\log^4 x\, Li_3(x)-\frac{7}{3}\log^3 x\, Li_4(x)+7\log^2 x\, Li_5(x)-14\log x\, Li_6(x)+14Li_7(x)$$

Therefore we obtain

(3.240c)
$$-6\sum_{n=1}^\infty \frac{H_n^{(1)}}{n^6}x^n+3\sum_{n=1}^\infty \frac{H_n^{(2)}}{n^5}x^n+3\sum_{n=1}^\infty \frac{\left[H_n^{(1)}\right]^2}{n^5}x^n$$

$$+5\log x\left[\sum_{n=1}^\infty \frac{H_n^{(1)}}{n^5}x^n-3\sum_{n=1}^\infty \frac{\left[H_n^{(1)}\right]^2}{n^4}x^n-3\sum_{n=1}^\infty \frac{H_n^{(2)}}{n^4}x^n\right]$$

$$+\log^2 x\left[\sum_{n=1}^\infty \frac{\left[H_n^{(1)}\right]^2}{n^3}x^n-2\sum_{n=1}^\infty \frac{H_n^{(1)}}{n^4}x^n+\sum_{n=1}^\infty \frac{H_n^{(2)}}{n^3}x^n\right]$$

$$-\frac{7}{30}\log^5 x\, Li_2(x)+\frac{5}{6}\log^4 x\, Li_3(x)-\frac{7}{3}\log^3 x\, Li_4(x)$$

$$+7\log^2 x\, Li_5(x)-14\log x\, Li_6(x)+14Li_7(x)=0$$

With $x=1$ we get

(3.240d)
$$-6\sum_{n=1}^\infty \frac{H_n^{(1)}}{n^6}+3\sum_{n=1}^\infty \frac{H_n^{(2)}}{n^5}+3\sum_{n=1}^\infty \frac{\left[H_n^{(1)}\right]^2}{n^5}=-14\varsigma(7)$$

The above Euler sums are also reported individually by Flajolet and Salvy [69].

The above operations may be repeated indefinitely to produce similar formulae for higher order Euler sums.

From (3.29) we see that

$$-\int_0^x \frac{\log(1-t)\log t}{(1-t)t}\,dt = \sum_{n=1}^{\infty} H_n^{(1)} \int_0^x t^{n-1}\log t\,dt$$

We have using integration by parts

$$\int_0^x \frac{\log(1-t)\log t}{(1-t)t}\,dt = \log(1-x)Li_2(1-x) + \log x\, Li_2(x) - Li_3(1-x) + Li_3(x) + \varsigma(3)$$

and hence we obtain

$$\sum_{n=1}^{\infty} \frac{H_n^{(1)}}{n^2} x^n - \log x \sum_{n=1}^{\infty} \frac{H_n^{(1)}}{n} x^n = \log(1-x)Li_2(1-x) + \log x\, Li_2(x) - Li_3(1-x) + Li_3(x) + \varsigma(3)$$

From (3.26) it is seen that

$$\log x \sum_{n=1}^{\infty} \frac{H_n^{(1)}}{n} x^n = \frac{1}{2}\log^2(1-x)\log x + Li_2(x)\log x$$

and hence

(3.241)

$$\sum_{n=1}^{\infty} \frac{H_n^{(1)}}{n^2} x^n = \frac{1}{2}\log^2(1-x)\log x + \log(1-x)Li_2(1-x) + 2\log x\, Li_2(x) - Li_3(1-x) + Li_3(x) + \varsigma(3)$$

With $x=1$ we have yet again

$$\sum_{n=1}^{\infty} \frac{H_n^{(1)}}{n^2} = 2\varsigma(3)$$

Using integration by parts it is easily shown that

$$\int_0^t \frac{\log^2(1-z)}{z(1-z)}\,dz = -\frac{1}{3}\log^3(1-t) + \log^2(1-t)\log t + 2\log(1-t)Li_2(1-t) - 2Li_3(1-t) + 2\varsigma(3)$$

and hence, using (3.221), we obtain for $t \in [0,1)$

(3.242)

$$\sum_{n=1}^{\infty} \frac{\left[H_n^{(1)}\right]^2}{n} t^n - \sum_{n=1}^{\infty} \frac{H_n^{(2)}}{n} t^n =$$

$$-\frac{1}{3}\log^3(1-t) + \log^2(1-t)\log t + 2\log(1-t)Li_2(1-t) - 2Li_3(1-t) + 2\varsigma(3)$$

Using (3.233) we obtain (3.313a)

$$\text{(3.243)}\qquad \sum_{n=1}^{\infty}\frac{H_n^{(1)}}{n^2}t^n=\frac{1}{2}\log^2(1-t)\log t+\log(1-t)Li_2(1-t)-Li_3(1-t)+Li_3(t)+\varsigma(3)$$

A further differentiation of (3.221) gives us

$$\text{(3.244)}\qquad -\frac{\log^3(1-z)}{(1-z)^x}=\sum_{n=1}^{\infty}\frac{f_n^{(3)}(x)}{n!}z^n$$

and we have

$$f_n^{(3)}(x)=2x(x+1)(x+2)...(x+n-1)\left[\frac{1}{x^3}+\frac{1}{(x+1)^3}...+\frac{1}{(x+n-1)^3}\right]$$

$$-3x(x+1)(x+2)...(x+n-1)\left[\frac{1}{x}+\frac{1}{x+1}...+\frac{1}{x+n-1}\right]\left[\frac{1}{x^2}+\frac{1}{(x+1)^2}...+\frac{1}{(x+n-1)^2}\right]$$

$$+x(x+1)(x+2)...(x+n-1)\left[\frac{1}{x}+\frac{1}{x+1}...+\frac{1}{x+n-1}\right]^3$$

Thus we have $f_n^{(3)}(1)=n!\left(2H_n^{(3)}-3H_n^{(1)}H_n^{(2)}+\left[H_n^{(1)}\right]^3\right)$ and therefore we see that

$$\text{(3.245)}\qquad -\frac{\log^3(1-z)}{(1-z)}=\sum_{n=1}^{\infty}\left(2H_n^{(3)}-3H_n^{(1)}H_n^{(2)}+\left[H_n^{(1)}\right]^3\right)z^n$$

Another integration by parts shows us that

$$\int_0^t\frac{\log^3(1-z)}{z(1-z)}dz=-\frac{1}{3}\log^4(1-t)+\log^3(1-t)\log t+3\log^2(1-t)Li_2(1-t)$$

$$-6\log(1-t)Li_3(1-t)+6Li_4(1-t)-6\varsigma(4)$$

and hence we obtain for $t\in[0,1)$

$$\text{(3.246)}\qquad \sum_{n=1}^{\infty}\frac{\left(2H_n^{(3)}-3H_n^{(1)}H_n^{(2)}+\left[H_n^{(1)}\right]^3\right)}{n}t^n$$

$$= -\frac{1}{3}\log^4(1-t) + \log^3(1-t)\log t + 3\log^2(1-t)Li_2(1-t) - 6\log(1-t)Li_3(1-t)$$

$$+ 6Li_4(1-t) - 6\varsigma(4)$$

We now divide (3.242) by $t$ and integrate to obtain

$$\sum_{n=1}^{\infty}\frac{\left[H_n^{(1)}\right]^2}{n^2} - \sum_{n=1}^{\infty}\frac{H_n^{(2)}}{n^2} =$$

$$-\frac{1}{3}\int_0^1 \frac{\log^3(1-t)}{t}dt + \int_0^1 \frac{\log^2(1-t)\log t}{t}dt + 2\int_0^1 \frac{\log(1-t)Li_2(1-t)}{t}dt + 2\int_0^1 \frac{\varsigma(3) - Li_3(1-t)}{t}dt$$

We have the following integrals

$$\int_0^1 \frac{\log^3(1-t)}{t}dt = -6\varsigma(4) \qquad \text{from (3.108a)}$$

$$\int_0^1 \frac{\log^2(1-t)\log t}{t}dt = -\frac{1}{2}\varsigma(4) \qquad \text{from (4.4.167q)}$$

$$\int_0^1 \frac{\log(1-t)Li_2(1-t)}{t}dt = -\frac{3}{4}\varsigma(4) \qquad \text{from (4.4.167r)}$$

$$\int_0^1 \frac{\varsigma(3) - Li_3(1-t)}{t}dt = \frac{5}{4}\varsigma(4) \qquad \text{from (4.4.167v)}$$

and hence we get

(3.247) $$\sum_{n=1}^{\infty}\frac{\left[H_n^{(1)}\right]^2}{n^2} - \sum_{n=1}^{\infty}\frac{H_n^{(2)}}{n^2} = \frac{5}{4}\varsigma(4)$$

This concurs with (4.4.168) and (4.4.167s) where we show that

$$\sum_{n=1}^{\infty}\frac{\left[H_n^{(1)}\right]^2}{n^2} = \frac{17}{4}\varsigma(4)$$

$$\sum_{n=1}^{\infty}\frac{H_n^{(2)}}{n^2} = \frac{7}{4}\varsigma(4)$$

and hence we have a useful check on my arithmetic.

We now divide (3.243) by $t$ and integrate to obtain

$$2\sum_{n=1}^{\infty}\frac{H_n^{(1)}}{n^3}=\int_0^1\frac{\log^2(1-t)\log t}{t}dt+2\int_0^1\frac{\log(1-t)Li_2(1-t)}{t}dt+2\int_0^1\frac{\varsigma(3)-Li_3(1-t)}{t}dt+2\int_0^1\frac{Li_3(t)}{t}dt$$

Therefore, employing the above definite integrals, we obtain as in (3.108d)

(3.248) $$\sum_{n=1}^{\infty}\frac{H_n^{(1)}}{n^3}=\frac{5}{4}\varsigma(4)$$

Repeating the same integration exercise again we get from (3.246)

(3.248a) $$\sum_{n=1}^{\infty}\frac{\left(2H_n^{(3)}-3H_n^{(1)}H_n^{(2)}+\left[H_n^{(1)}\right]^3\right)}{n^2}=$$

$$-\frac{1}{3}\int_0^1\frac{\log^4(1-t)}{t}dt+\int_0^1\frac{\log^3(1-t)\log t}{t}dt+3\int_0^1\frac{\log^2(1-t)Li_2(1-t)}{t}dt$$

$$-6\int_0^1\frac{\log(1-t)Li_3(1-t)}{t}dt-6\int_0^1\frac{\varsigma(4)-Li_4(1-t)}{t}dt$$

We have from (4.4.231)

$$\int_0^1\frac{\log^{r-1}t}{1-t}dt=\int_0^1\frac{\log^{r-1}(1-t)}{t}dt=(-1)^{r-1}\varsigma(r)\Gamma(r) \qquad ,r\geq 2$$

and accordingly we get

$$\int_0^1\frac{\log^4(1-t)}{t}dt=24\varsigma(5)$$

The Wolfram Integrator provides the result (which is also easily obtained using integration by parts)

(3.249)

$$\int\frac{\log^3(1-t)}{t}dt=\log^3(1-t)\log t+3\log^2(1-t)Li_2(1-t)-6\log(1-t)Li_3(1-t)+6Li_4(1-t)$$

We then obtain using integration by parts over the interval $[a,1]$

$$\int_a^1 \frac{\log^3(1-t)}{t} \log t \, dt =$$

$$\left(\log^3(1-t)\log t + 3\log^2(1-t)Li_2(1-t) - 6\log(1-t)Li_3(1-t) + 6Li_4(1-t)\right)\log t\Big|_a^1 + 6\varsigma(4)\log a$$

$$-\int_a^1 \frac{\log^3(1-t)\log t}{t} dt - 3\int_a^1 \frac{\log^2(1-t)Li_2(1-t)}{t} dt + 6\int_a^1 \frac{\log(1-t)Li_3(1-t)}{t} dt + 6\int_a^1 \frac{\varsigma(4) - Li_4(1-t)}{t} dt$$

where, in order to achieve convergence, we have judiciously added and subtracted a factor of

$$6\varsigma(4)\log a = 6\int_a^1 \frac{\varsigma(4)}{t} dt$$

Therefore, in the limit as $a \to 0$ we obtain

$$2\int_0^1 \frac{\log^3(1-t)\log t}{t} dt =$$

$$-3\int_a^1 \frac{\log^2(1-t)Li_2(1-t)}{t} dt + 6\int_a^1 \frac{\log(1-t)Li_3(1-t)}{t} dt + 6\int_a^1 \frac{\varsigma(4) - Li_4(1-t)}{t} dt$$

Using an obvious substitution we have

$$\int_0^1 \frac{\log^2(1-t)Li_2(1-t)}{t} dt = \int_0^1 \frac{\log^2 t \, Li_2(t)}{1-t} dt$$

and we will recognise the latter integral from (4.4.229)

$$\sum_{n=1}^\infty \frac{H_n^{(r)}}{n^q} = \varsigma(r)\varsigma(q) - \frac{(-1)^{r-1}}{(r-1)!}\int_0^1 \frac{\log^{r-1} x \, Li_q(x)}{1-x} dx \qquad \text{for } q, r \geq 2$$

Hence we have

$$(3.250) \qquad \sum_{n=1}^\infty \frac{H_n^{(3)}}{n^2} = \varsigma(2)\varsigma(3) - \frac{1}{2}\int_0^1 \frac{\log^2 x \, Li_2(x)}{1-x} dx$$

and therefore

$$(3.251) \qquad \int_0^1 \frac{\log^2(1-t)Li_2(1-t)}{t} dt = 2\varsigma(2)\varsigma(3) - 2\sum_{n=1}^\infty \frac{H_n^{(3)}}{n^2}$$

Similarly we obtain

(3.252) $$\int_0^1 \frac{\log(1-t)Li_3(1-t)}{t}dt = \int_0^1 \frac{\log t\, Li_3(t)}{1-t}dt = \sum_{n=1}^{\infty}\frac{H_n^{(2)}}{n^3} - \varsigma(2)\varsigma(3)$$

With integration by parts we get

$$\int_0^1 \frac{\varsigma(4) - Li_4(1-x)}{x}dx = -\int_0^1 \frac{\log x\, Li_3(1-x)}{1-x}dx$$

$$\int_0^1 \frac{\log x\, Li_3(1-x)}{1-x}dx = -\varsigma(2)\varsigma(3) + \int_0^1 \frac{\left[Li_2(1-x)\right]^2}{1-x}dx$$

The following integral was given by Freitas [69a] and is also derived in (3.211e)

$$\int_0^1 \frac{\left[Li_2(1-x)\right]^2}{1-x}dx = \int_0^1 \frac{\left[Li_2(x)\right]^2}{x}dx = 2\varsigma(2)\varsigma(3) - 3\varsigma(5)$$

and hence we get

(3.253) $$\int_0^1 \frac{\varsigma(4) - Li_4(1-x)}{x}dx = 3\varsigma(5) - \varsigma(2)\varsigma(3)$$

We then deduce that

$$2\int_0^1 \frac{\log^3(1-t)\log t}{t}dt =$$

$$-3\left(2\varsigma(2)\varsigma(3) - 2\sum_{n=1}^{\infty}\frac{H_n^{(3)}}{n^2}\right) + 6\left(\sum_{n=1}^{\infty}\frac{H_n^{(2)}}{n^3} - \varsigma(2)\varsigma(3)\right) + 6\left(3\varsigma(5) - \varsigma(2)\varsigma(3)\right)$$

$$= -18\varsigma(2)\varsigma(3) + 6\left(\sum_{n=1}^{\infty}\frac{H_n^{(2)}}{n^3} + \sum_{n=1}^{\infty}\frac{H_n^{(3)}}{n^2}\right) + 18\varsigma(5)$$

We have from (4.4.232a)

$$\sum_{n=1}^{\infty}\frac{H_n^{(p)}}{n^q} + \sum_{n=1}^{\infty}\frac{H_n^{(q)}}{n^p} = \varsigma(p)\varsigma(q) + \varsigma(p+q)$$

and hence we see that

$$\sum_{n=1}^{\infty}\frac{H_n^{(2)}}{n^3}+\sum_{n=1}^{\infty}\frac{H_n^{(3)}}{n^2}=\varsigma(2)\varsigma(3)+\varsigma(5)$$

Georghiou and Philippou [69c] gave the following formula in 1983

$$\sum_{n=1}^{\infty}\frac{H_k^{(2)}}{k^{2n+1}}=\varsigma(2)\varsigma(2n+1)-\frac{(n+2)(2n+1)}{2}\varsigma(2n+3)+2\sum_{j=2}^{n+1}(j-1)\varsigma(2j-1)\varsigma(2n+4-2j)$$

and this gives us for $n=1$

(3.255) $$\sum_{n=1}^{\infty}\frac{H_n^{(2)}}{n^3}=3\varsigma(2)\varsigma(3)-\frac{9}{2}\varsigma(5)$$

Therefore we immediately get

(3.256) $$\sum_{n=1}^{\infty}\frac{H_n^{(3)}}{n^2}=-2\varsigma(2)\varsigma(3)+\frac{11}{2}\varsigma(5)$$

We therefore obtain

(3.254) $$\int_0^1\frac{\log^3(1-t)\log t}{t}dt=-6\varsigma(2)\varsigma(3)+12\varsigma(5)$$

in agreement with the integral reported by Coffey in [45d]. Reference should also be made to Zheng's recent paper [142b]. The Wolfram Integrator was unable to evaluate this integral.

This gives us a set of relevant integrals

(3.257) $$\int_0^1\frac{\log^4(1-t)}{t}dt=24\varsigma(5)$$

(3.258) $$\int_0^1\frac{\log^3(1-t)\log t}{t}dt=-6\varsigma(2)\varsigma(3)+12\varsigma(5)$$

(3.259) $$\int_0^1\frac{\log^2(1-t)Li_2(1-t)}{t}dt=6\varsigma(2)\varsigma(3)-11\varsigma(5)$$

(3.260) $$\int_0^1\frac{\log(1-t)Li_3(1-t)}{t}dt=2\varsigma(2)\varsigma(3)-\frac{9}{2}\varsigma(5)$$

(3.261) $$\int_0^1\frac{\varsigma(4)-Li_4(1-x)}{x}dx=-\varsigma(2)\varsigma(3)+3\varsigma(5)$$

Now, after substituting the individual integral evaluations in (3.248a) we obtain

$$\sum_{n=1}^{\infty}\frac{\left(2H_n^{(3)}-3H_n^{(1)}H_n^{(2)}+\left[H_n^{(1)}\right]^3\right)}{n^2}=6\varsigma(2)\varsigma(3)-20\varsigma(5)$$

We may safely rearrange the above convergent series to produce

$$\sum_{n=1}^{\infty}\frac{\left(2H_n^{(3)}-3H_n^{(1)}H_n^{(2)}+\left[H_n^{(1)}\right]^3\right)}{n^2}=2\sum_{n=1}^{\infty}\frac{H_n^{(3)}}{n^2}-3\sum_{n=1}^{\infty}\frac{H_n^{(1)}H_n^{(2)}}{n^2}+\sum_{n=1}^{\infty}\frac{\left[H_n^{(1)}\right]^3}{n^2}$$

$$=-4\varsigma(2)\varsigma(3)+7\varsigma(5)-3\sum_{n=1}^{\infty}\frac{H_n^{(1)}H_n^{(2)}}{n^2}+\sum_{n=1}^{\infty}\frac{\left[H_n^{(1)}\right]^3}{n^2}$$

and we then deduce that

(3.262) $$\sum_{n=1}^{\infty}\frac{\left[H_n^{(1)}\right]^3}{n^2}-3\sum_{n=1}^{\infty}\frac{H_n^{(1)}H_n^{(2)}}{n^2}=10\varsigma(2)\varsigma(3)-21\varsigma(5)$$

We now need to carry out some algebraic manipulation of series as follows

$$\sum_{n=1}^{\infty}\frac{\left(H_n^{(1)}\right)^3}{n^2}=\sum_{n=1}^{\infty}\frac{\left(H_{n-1}^{(1)}+\frac{1}{n}\right)^3}{n^2}=\sum_{n=1}^{\infty}\frac{\left(H_{n-1}^{(1)}\right)^3}{n^2}+3\sum_{n=1}^{\infty}\frac{\left(H_{n-1}^{(1)}\right)^2}{n^3}+3\sum_{n=1}^{\infty}\frac{H_{n-1}^{(1)}}{n^4}+\sum_{n=1}^{\infty}\frac{1}{n^5}$$

$$\sum_{n=1}^{\infty}\frac{\left(H_{n-1}^{(1)}\right)^2}{n^3}=\sum_{n=1}^{\infty}\frac{\left(H_n^{(1)}-\frac{1}{n}\right)^2}{n^3}$$

$$=\sum_{n=1}^{\infty}\frac{\left(H_n^{(1)}\right)^2}{n^3}-2\sum_{n=1}^{\infty}\frac{H_n^{(1)}}{n^4}+\sum_{n=1}^{\infty}\frac{1}{n^5}$$

Therefore we get

$$\sum_{n=1}^{\infty}\frac{\left(H_n^{(1)}\right)^3}{n^2}=\sum_{n=1}^{\infty}\frac{\left(H_{n-1}^{(1)}\right)^3}{n^2}+3\left[\sum_{n=1}^{\infty}\frac{\left(H_n^{(1)}\right)^2}{n^3}-2\sum_{n=1}^{\infty}\frac{H_n^{(1)}}{n^4}+\sum_{n=1}^{\infty}\frac{1}{n^5}\right]+3\sum_{n=1}^{\infty}\frac{H_{n-1}^{(1)}}{n^4}+\sum_{n=1}^{\infty}\frac{1}{n^5}$$

$$=\sum_{n=1}^{\infty}\frac{\left(H_{n-1}^{(1)}\right)^3}{n^2}+3\left[\sum_{n=1}^{\infty}\frac{\left(H_n^{(1)}\right)^2}{n^3}-2\sum_{n=1}^{\infty}\frac{H_n^{(1)}}{n^4}+\sum_{n=1}^{\infty}\frac{1}{n^5}\right]+3\left[\sum_{n=1}^{\infty}\frac{H_n^{(1)}}{n^4}-\sum_{n=1}^{\infty}\frac{1}{n^5}\right]+\sum_{n=1}^{\infty}\frac{1}{n^5}$$

and hence we have

$$\sum_{n=1}^{\infty}\frac{\left(H_n^{(1)}\right)^3}{n^2}=\sum_{n=1}^{\infty}\frac{\left(H_{n-1}^{(1)}\right)^3}{n^2}+3\sum_{n=1}^{\infty}\frac{\left(H_n^{(1)}\right)^2}{n^3}-3\sum_{n=1}^{\infty}\frac{H_n^{(1)}}{n^4}+\sum_{n=1}^{\infty}\frac{1}{n^5}$$

By Flajolet and Salvy [69, Cor. 5.2] we have

(3.263) $$\sum_{n=1}^{\infty}\frac{\left(H_{n-1}^{(1)}\right)^3}{n^2}=\sum_{n=1}^{\infty}\frac{\left(H_n^{(1)}\right)^3}{(n+1)^2}=\varsigma(2)\varsigma(3)+\frac{15}{2}\varsigma(5)$$

From [69, Th. 3.1] we have

(3.264) $$\sum_{n=1}^{\infty}\frac{\left(H_n^{(1)}\right)^2}{n^3}=-\varsigma(2)\varsigma(3)+\frac{7}{2}\varsigma(5)$$

and [69, Th. 2.2]

(3.265) $$\sum_{n=1}^{\infty}\frac{H_n^{(1)}}{n^4}=-\varsigma(2)\varsigma(3)+3\varsigma(5)$$

and we may accordingly obtain $\sum_{n=1}^{\infty}\frac{H_n^{(1)}H_n^{(2)}}{n^2}$ from (3.262). See also (4.3.60c).

We note the result obtained by Borwein and Girgensohn [28a] in 1996

(3.266) $$\sum_{n=1}^{\infty}\frac{H_{n-1}^{(1)}H_{n-1}^{(2)}}{n^2}=-\varsigma(2)\varsigma(3)+\frac{7}{2}\varsigma(5)$$

and also see Zheng's recent paper [142b].

## SUMMARY OF HARMONIC NUMBER SERIES IDENTITIES

Some of the series derived in these papers are set out below for ease of reference.

(3.105a) $$\sum_{n=1}^{\infty}\frac{H_n}{n}x^n=\frac{1}{2}\log^2(1-x)+Li_2(x)$$

(3.105d)
$$\sum_{n=1}^{\infty}\frac{H_n}{n^2}x^n=\frac{1}{2}\log^2(1-x)\log x+\log(1-x)Li_2(1-x)+Li_3(x)-Li_3(1-x)+\varsigma(3)$$

(3.106) $$\frac{1}{3}\log^3(1-x)+2Li_3(x)=2\sum_{n=1}^{\infty}\frac{H_n}{n^2}x^n+\sum_{n=1}^{\infty}\frac{H_n^{(2)}}{n}x^n-\sum_{n=1}^{\infty}\frac{\left(H_n\right)^2}{n}x^n$$

(3.211i)

$$\sum_{n=1}^{\infty}\frac{H_n^{(2)}}{n}x^n=Li_3(x)-\log(1-x)Li_2(x)-\log x\log^2(1-x)-2\log(1-x)\,Li_2(1-x)+2Li_3(1-x)-2\varsigma(3)$$

(3.106e)

$$2\sum_{n=1}^{\infty}\frac{\left(H_n\right)^2}{n}x^n=\log^2(1-x)\log x-\frac{1}{3}\log^3(1-x)+2\log(1-x)Li_2(1-x)-2Li_3(1-x)$$

$$+2\varsigma(3)-2Li_3\left(-\frac{x}{1-x}\right)$$

(3.106f) $$\sum_{n=1}^{\infty}\frac{H_n^{(2)}}{n}x^n=-2Li_3\left(-\frac{x}{1-x}\right)$$

(3.242) $$\sum_{n=1}^{\infty}\frac{\left[H_n^{(1)}\right]^2}{n}x^n-\sum_{n=1}^{\infty}\frac{H_n^{(2)}}{n}x^n=$$

$$-\frac{1}{3}\log^3(1-x)+\log^2(1-x)\log x+2\log(1-x)Li_2(1-x)-2Li_3(1-x)+2\varsigma(3)$$

(4.4.43zh)

$$\sum_{n=1}^{\infty}\frac{H_n^{(2)}}{n}x^n=Li_3(x)-\log(1-x)Li_2(x)-\log x\log^2(1-x)-2\log(1-x)\,Li_2(1-x)+2Li_3(1-x)-2\varsigma(3)$$

(3.108b) $$2\sum_{n=1}^{\infty}\frac{H_n^{(1)}}{n^3}x^n+\sum_{n=1}^{\infty}\frac{H_n^{(2)}}{n^2}x^n-\sum_{n=1}^{\infty}\frac{\left(H_n^{(1)}\right)^2}{n^2}x^n=$$

$$\frac{1}{3}\log^3(1-x)\log x+\log^2(1-x)Li_2(1-x)-2\log(1-x)Li_3(1-x)+2Li_4(1-x)+2Li_4(x)-2\varsigma(4)$$

(3.211) $$3Li_4(x)-3\log x\,Li_3(x)-Li_2(x)Li_2(1-x)+\varsigma(2)Li_2(x)-\frac{1}{2}\left[Li_2(x)\right]^2=$$

$$2\sum_{n=1}^{\infty}\frac{H_n^{(1)}}{n^3}x^n-2\log x\sum_{n=1}^{\infty}\frac{H_n^{(1)}}{n^2}x^n+\sum_{n=1}^{\infty}\frac{H_n^{(2)}}{n^2}x^n-\log x\sum_{n=1}^{\infty}\frac{H_n^{(2)}}{n}x^n$$

(3.211b) $$2\sum_{n=1}^{\infty}\frac{H_n^{(1)}}{n^3}x^n+\sum_{n=1}^{\infty}\frac{H_n^{(2)}}{n^2}x^n=3Li_4(x)+\frac{1}{2}\left[Li_2(x)\right]^2$$

(3.109b) $$2\sum_{n=1}^{\infty}\frac{H_n^{(1)}}{n^3}x^n-\int_0^1\frac{\log^2 y\log[1-x(1-y)]}{1-y}dy-2\sum_{n=1}^{\infty}\frac{\left(H_n^{(1)}\right)^2}{n^2}x^n=$$

$$\frac{1}{3}\log^3(1-x)\log x+\log^2(1-x)Li_2(1-x)-2\log(1-x)Li_3(1-x)+2Li_4(1-x)+2Li_4(x)-2\varsigma(4)$$

(3.217d)

$$2\sum_{n=1}^{\infty}\frac{H_n^{(1)}}{n^3}x^n-\log x\sum_{n=1}^{\infty}\frac{H_n^{(1)}}{n^2}x^n=-\sum_{n=1}^{\infty}\frac{1}{n^2}\sum_{k=1}^{n}\frac{(1-x)^k}{k^2}-\varsigma(2)Li_2(x)+\frac{7}{4}\varsigma(4)$$

$$-\log x\,Li_3(x)+\left[\varsigma(3)-Li_3(1-x)\right]\log x-\frac{1}{2}\left[Li_2(1-x)\right]^2+\frac{1}{2}\varsigma^2(2)+2Li_4(x)$$

(3.235) $$6\sum_{n=1}^{\infty}\frac{H_n^{(1)}}{n^3}x^n-3\sum_{n=1}^{\infty}\frac{\left[H_n^{(1)}\right]^2}{n^2}x^n-6Li_4(x)+3\sum_{n=1}^{\infty}\frac{H_n^{(2)}}{n^2}x^n=$$

$$\log x\log^3(1-x)+3Li_2(1-x)\log^2(1-x)-6Li_3(1-x)\log(1-x)+6Li_4(1-x)-6\varsigma(4)$$

(3.110a) $$\sum_{n=1}^{\infty}\frac{\left(H_n^{(1)}\right)^2}{n}x^n=-\frac{1}{3}\log^3(1-x)+Li_3(x)-Li_2(x)\log(1-x)$$

(3.110ea)

$$\sum_{n=1}^{\infty}\frac{\left(H_n^{(1)}\right)^2}{n^2}x^n=Li_4(x)+\frac{1}{2}\left[Li_2(x)\right]^2$$

$$-\frac{1}{3}\left[\log^3(1-x)\log x+3\log^2(1-x)Li_2(1-x)-6\log(1-x)Li_3(1-x)+6Li_4(1-x)-6\varsigma(4)\right]$$

(3.110f)

$$\frac{1}{3}\log^4(1-x)+2Li_4(x)=\sum_{n=1}^{\infty}\left\{\frac{1}{3}\frac{\left(H_n^{(1)}\right)^3}{n}-\frac{\left(H_n^{(1)}\right)^2}{n^2}+2\frac{H_n^{(1)}}{n^3}-\frac{H_n^{(1)}H_n^{(2)}}{n}+\frac{H_n^{(2)}}{n^2}+\frac{2}{3}\frac{H_n^{(3)}}{n}\right\}x^n$$

(3.110g) $$\sum_{n=1}^{\infty}\left\{\frac{1}{3}\frac{\left(H_n^{(1)}\right)^3}{n^2}-\frac{\left(H_n^{(1)}\right)^2}{n^3}+2\frac{H_n^{(1)}}{n^4}-\frac{H_n^{(1)}H_n^{(2)}}{n^2}+\frac{H_n^{(2)}}{n^3}+\frac{2}{3}\frac{H_n^{(3)}}{n^2}\right\}x^n$$

$$=8\left[\begin{array}{l}\frac{1}{24}\log^4(1-x)\log x+\frac{1}{6}\log^3(1-x)Li_2(1-x)-\frac{1}{2}\log^2(1-x)Li_3(1-x)\\ \\ +\log(1-x)Li_4(1-x)-Li_5(1-x)+\varsigma(5)+\frac{1}{4}Li_5(x)\end{array}\right]$$

(3.111d) $$\sum_{n=1}^{\infty}\frac{H_n}{n^2}x^n=-Li_2(x)\log x-Li_2\left(\frac{-x}{1-x}\right)\log x+\log(1-x)Li_2(1-x)$$

$$-Li_3(1-x)+Li_3(x)+\varsigma(3)$$

(3.241)

$$\sum_{n=1}^{\infty}\frac{H_n^{(1)}}{n^2}x^n=\frac{1}{2}\log^2(1-x)\log x+\log(1-x)Li_2(1-x)+2\log x\,Li_2(x)-Li_3(1-x)+Li_3(x)+\varsigma(3)$$

(3.211ei) $$6Li_5(x)-3\log x\,Li_4(x)-\frac{1}{2}\left[Li_2(x)\right]^2\log x+4\sum_{n=1}^{\infty}\frac{H_n^{(1)}}{n^4}+2\sum_{n=1}^{\infty}\frac{H_n^{(2)}}{n^3}-6\varsigma(5)=$$

$$4\sum_{n=1}^{\infty}\frac{H_n^{(1)}}{n^4}x^n-2\log x\sum_{n=1}^{\infty}\frac{H_n^{(1)}}{n^3}x^n+2\sum_{n=1}^{\infty}\frac{H_n^{(2)}}{n^3}x^n-\log x\sum_{n=1}^{\infty}\frac{H_n^{(2)}}{n^2}x^n$$

(3.239)

$$2\sum_{n=1}^{\infty}\frac{H_n^{(1)}}{n^4}x^n-\sum_{n=1}^{\infty}\frac{H_n^{(2)}}{n^3}x^n-\sum_{n=1}^{\infty}\frac{\left[H_n^{(1)}\right]^2}{n^3}x^n+\log x\left[\sum_{n=1}^{\infty}\frac{\left[H_n^{(1)}\right]^2}{n^2}x^n-2\sum_{n=1}^{\infty}\frac{H_n^{(1)}}{n^3}x^n+\sum_{n=1}^{\infty}\frac{H_n^{(2)}}{n^2}x^n\right]$$

$$+\log^2 x\,Li_3(x)-\frac{1}{3}\log^3 x\,Li_2(x)=0$$

(3.240a)

$$4\sum_{n=1}^{\infty}\frac{H_n^{(1)}}{n^5}x^n-2\sum_{n=1}^{\infty}\frac{H_n^{(2)}}{n^4}x^n-2\sum_{n=1}^{\infty}\frac{\left[H_n^{(1)}\right]^2}{n^4}x^n+\log x\left[\sum_{n=1}^{\infty}\frac{\left[H_n^{(1)}\right]^2}{n^3}x^n-2\sum_{n=1}^{\infty}\frac{H_n^{(1)}}{n^4}x^n+\sum_{n=1}^{\infty}\frac{H_n^{(2)}}{n^4}x^n\right]$$

$$+\log^3 x\, Li_3(x)-\frac{1}{3}\log^4 x\, Li_2(x)-2\log^2 x\, Li_4(x)+4\log x\, Li_5(x)-4Li_6(x)=0$$

(3.240c) $$-6\sum_{n=1}^{\infty}\frac{H_n^{(1)}}{n^6}x^n+3\sum_{n=1}^{\infty}\frac{H_n^{(2)}}{n^5}x^n+3\sum_{n=1}^{\infty}\frac{\left[H_n^{(1)}\right]^2}{n^5}x^n$$

$$+5\log x\left[\sum_{n=1}^{\infty}\frac{H_n^{(1)}}{n^5}x^n-3\sum_{n=1}^{\infty}\frac{\left[H_n^{(1)}\right]^2}{n^4}x^n-3\sum_{n=1}^{\infty}\frac{H_n^{(2)}}{n^4}x^n\right]$$

$$+\log^2 x\left[\sum_{n=1}^{\infty}\frac{\left[H_n^{(1)}\right]^2}{n^3}x^n-2\sum_{n=1}^{\infty}\frac{H_n^{(1)}}{n^4}x^n+\sum_{n=1}^{\infty}\frac{H_n^{(2)}}{n^3}x^n\right]$$

$$-\frac{7}{30}\log^5 x\, Li_2(x)+\frac{5}{6}\log^4 x\, Li_3(x)-\frac{7}{3}\log^3 x\, Li_4(x)$$

$$+7\log^2 x\, Li_5(x)-14\log x\, Li_6(x)+14Li_7(x)=0$$

(3.246) $$\sum_{n=1}^{\infty}\frac{\left(2H_n^{(3)}-3H_n^{(1)}H_n^{(2)}+\left[H_n^{(1)}\right]^3\right)}{n}x^n$$

$$=-\frac{1}{3}\log^4(1-x)+\log^3(1-x)\log x+3\log^2(1-x)Li_2(1-x)-6\log(1-x)Li_3(1-x)$$

$$+6Li_4(1-x)-6\varsigma(4)$$

.

Donal F. Connon
Elmhurst
Dundle Road
Matfield
Kent TN12 7HD
dconnon@btopenworld.com